\documentclass[a4paper]{article}

\usepackage{pstricks,pst-plot}
\usepackage{amsfonts}

\usepackage[latin1]{inputenc}
\usepackage[frenchb]{babel}  
\def\C{\mathbb{C}}

\def\g{\mathfrak{g}}

\def\z{\mathfrak{z}}
\def\p{\mathfrak{p}}
\def\q{\mathfrak{q}}
\def\m{\mathfrak{m}}

\def\u{\mathfrak{u}}

\def\n{\mathfrak{n}}

\def\sl2{\mathfrak{sl}_{2}-{\rm triplet}}

\def\dem{{\it D{\'e}monstration. ~}}
\def\rem{{\bf Remarque} - }
\def\qed{\hfill \rule{.2cm}{.2cm} \\}

\let \epsilon=\varepsilon
\def \phi {\varphi}

\usepackage{amsmath}

\begin{document}

\newtheorem{lemme}{Lemme}
\newtheorem{defi}{D{\'e}finition}
\newtheorem{prop}{Proposition}
\newtheorem{pte}{Propri{\'e}t{\'e}}
\newtheorem{thm}{Th{\'e}or{\`e}me}
\newtheorem{cor}{Corollaire}

\def\dv#1#2{\langle {#1},{#2}\rangle}
\def\an#1#2{\def\deux{#2} \ifx\deux\empty {\cal O}_{#1} 
\else {\cal O}_{#1,#2} \fi }
\def\tk#1#2{{#2}\otimes _{#1}}
\def\tens{\raisebox{.3mm}{\scriptsize$\otimes$}}

\large
\begin{center} \LARGE
{\bf Indice du normalisateur du centralisateur d'un {\'e}l{\'e}ment nilpotent
  dans une alg{\`e}bre de Lie semi-simple.}\\
\vspace{1cm}
\large
{\sc Anne Moreau }\\

\vspace{1cm}

\footnotesize
\begin{center}
{\bf Abstract}
\end{center}
\begin{quotation}
\vskip .5em
The index of a complex Lie algebra is the minimal codimension of its
coadjoint orbits. Let us suppose $\g$ semisimple, then its index, ${\rm
ind} \; \g$, is equal to its rank, ${\rm rk \; \g}$. The goal of this paper is to establish a simple general formula for
the index of $\n(\g^{\xi})$, for $\xi$ nilpotent, where $\n(\g^{\xi})$
is the normaliser in $\g$ of the centraliser $\g^{\xi}$ of $\xi$. More precisely, we
have to show the following result, conjectured by D. Panyushev
\cite{Panyushev}\,:
$${\rm ind \;} \n(\g^{\xi}) = {\rm rk \; \g}-\dim \z(\g^{\xi}),$$ 
where $\z(\g^{\xi})$ is the center of
$\g^{\xi}$. D. Panyushev obtained in \cite{Panyushev} the inequality \hbox{${\rm ind
  \;} \n(\g^{\xi}) \geq {\rm rg \; \g}-\dim \z(\g^{\xi})$} and we show
  that the maximality of the rank of a certain matrix with entries in the
  symmetric algebra ${\cal S}(\g^{\xi})$ implies the other
  inequality. The main part of this paper consists of the proof of the
  maximality of the rank
 of this matrix. 
\end{quotation}
\end{center}

\large
\section*{Introduction}

L'indice d'une alg{\`e}bre de Lie complexe est la codimension minimale de
ses orbites coadjointes. Si $\g$ est une alg{\`e}bre de Lie semi-simple
complexe, son indice ${\rm ind} \; \g$ est {\'e}gal {\`a} son rang, ${\rm rg \; \g}$. Plus g{\'e}n{\'e}ralement, l'indice d'une repr{\'e}sentation arbitraire $V$
d'une alg{\`e}bre de Lie complexe $\q$ est la codimension minimale de ses
orbites sous l'action contragr{\'e}diente. On le note ${\rm ind}(\q,V)$.
   
Dans tout ce qui suit, $\g$ est une alg{\`e}bre de Lie semi-simple
complexe de groupe adjoint $G$. On identifie $\g$ {\`a} son image par la repr{\'e}sentation adjointe. Pour $x$ dans $\g$, on
note $\g^{x}$ son centralisateur, $\z(\g^{x})$ le centre de $\g^{x}$ et
$\n(\g^{x})$ le normalisateur dans $\g$ du centralisateur
$\g^{x}$. Le but de cet article est de donner une expression 
simple de l'indice de $\n(\g^{\xi})$, pour $\xi$ {\'e}l{\'e}ment nilpotent de
$\g$. L'alg{\`e}bre $\n(\g^{\xi})$ agit sur le sous-espace
$\g^{\xi}$ par la repr{\'e}sentation adjointe et on {\'e}tablit en outre
une formule pour l'indice, ${\rm ind}(\n(\g^{\xi}),\g^{\xi})$, du
$\n(\g^{\xi})$-module $\g^{\xi}$. Plus pr{\'e}cis{\'e}ment, on se propose de
montrer les deux r{\'e}sultats suivants, conjectur{\'e}s
par D. Panyushev en \cite{Panyushev}, Conjectures 6.1 et 6.2\,:

\begin{thm}\label{intro1} Soit $\xi$ un {\'e}l{\'e}ment nilpotent de $\g$. Alors 
$${\rm ind \;} \n(\g^{\xi}) = {\rm rg \; \g}-\dim \z(\g^{\xi}) \cdot$$
\end{thm}

\begin{thm}\label{intro2} Soit $\xi$ un {\'e}l{\'e}ment nilpotent de $\g$. Alors 
$${\rm ind \;} (\n(\g^{\xi}),\g^{\xi}) = {\rm rg \; \g}-\dim \z(\g^{\xi}) \cdot$$
\end{thm}

Notons que les deux relations pr{\'e}c{\'e}dentes sont ind{\'e}pendantes du choix d'un repr{\'e}sentant
dans l'orbite de $\xi$ sous l'action du groupe adjoint. Dans \cite{Panyushev}, D. Panyushev dresse une liste de cas o{\`u} ces
deux {\'e}galit{\'e}s sont satisfaites. 
Remarquons {\`a} ce propos qu'il obtient seulement les relations\,: 
\hbox{${\rm  ind \;} \n(\g^{\xi}) = {\rm ind \;} \g^{\xi}-\dim \z(\g^{\xi})$} et 
${\rm ind \;} (\n(\g^{\xi}),\g^{\xi}) 
={\rm ind \;} \g^{\xi} -\dim \z(\g^{\xi})$. La relation \hbox{${\rm ind \;}
\g^{\xi}={\rm rg \; \g}$} n'est en effet {\'e}nonc{\'e}e dans \cite{Panyushev}
que sous forme de conjecture, \cite{Panyushev}, Conjecture 3.2 d'Elashvili. Cette {\'e}galit{\'e} est d{\'e}montr{\'e}e en \cite{Charbonnel}, Th{\'e}or{\`e}me 5.5. Dans tous ses exemples, il obtient ces
relations en montrant que le groupe $N_{\g}(\xi)$ a une orbite ouverte
dans ${\z(\g^{\xi})}^{*}$, o{\`u} $N_{\g}(\xi)$ d{\'e}signe le sous-groupe connexe
de $G$ d'alg{\`e}bre de Lie $\n(\g^{\xi})$. Cette condition est, comme on aura
l'occasion de le voir, suffisante mais ne permet pas de traiter tous les
cas. 
 
La premi{\`e}re partie regroupe un certain nombre de r{\'e}sultats autour du normalisateur du centralisateur d'un {\'e}l{\'e}ment
nilpotent. On introduit en outre dans cette partie une propri{\'e}t{\'e} $(P)$ qui
interviendra dans la suite. On montre dans la deuxi{\`e}me
partie que le th{\'e}or{\`e}me \ref{intro1} est en fait une cons{\'e}quence du
th{\'e}or{\`e}me \ref{intro2} et qu'obtenir l'identit{\'e} du th{\'e}or{\`e}me
\ref{intro2} {\'e}quivaut {\`a} montrer qu'une certaine matrice {\`a} coefficients dans
${\cal S}(\g^{\xi})$ est de rang maximal. On consacre les deux parties
qui suivent {\`a} la d{\'e}monstration de ce dernier point dans deux cas particuliers\,:
lorsque $\g$ est une alg{\`e}bre de Lie simple classique (partie 3)
et lorsque l'{\'e}l{\'e}ment $\xi$ v{\'e}rifie la propri{\'e}t{\'e} $(P)$ (partie
4). La troisi{\`e}me partie utilise des
propri{\'e}tes g{\'e}om{\'e}triques des alg{\`e}bres de Lie classiques tandis que la
quatri{\`e}me partie repose pour une large part sur des
r{\'e}sultats expos{\'e}s dans \cite{Charbonnel} par J.Y. Charbonnel. 
On {\'e}tudie dans la derni{\`e}re partie la propri{\'e}t{\'e}
$(P)$ pour achever la
d{\'e}monstration des th{\'e}or{\`e}mes \ref{intro1} et \ref{intro2} dans le cas
exceptionnel. La fin de la d{\'e}monstration consiste {\`a} v{\'e}rifier certaines
conditions sur un nombre fini de cas {\`a} l'aide du logiciel GAP4.

\section{R{\'e}sultats pr{\'e}liminaires.}

 Soit $\xi$ un {\'e}l{\'e}ment nilpotent de $\g$. Puisque $\xi$ est nilpotent,
 le th{\'e}or{\`e}me de Jacobson-Morosov assure l'existence de deux {\'e}l{\'e}ments
 $\rho$ et $\eta$ dans $\g$ pour lesquels $\xi, \rho, \eta$
 satisfont les relations de $\mathfrak{sl}_{2}$-triplet\,:
$$[\rho,\xi]=2 \xi, \hspace{1cm} [\xi,\eta]=\rho,  \hspace{1cm}
 [\rho,\eta]=-2 \eta \cdot$$
Le  normalisateur $\n(\g^{\xi})$ du
 centralisateur $\g^{\xi}$ de $\xi$ est, par d{\'e}finition, l'ensemble
 des $y$ de $\g$ tels qu'on ait l'inclusion\,: $[y,\g^{\xi}] \subset
 \g^{\xi}$. C'est aussi le normalisateur du centre  $\z(\g^{\xi})$ de 
$\g^{\xi}$, comme on le v{\'e}rifie facilement. La proposition suivante rassemble un certain nombre de propri{\'e}t{\'e}s
 du normalisateur. On en trouve une preuve dans \cite{Tauvel},
 Chapitre XVII, Propositions 5.6 et 5.12.
\begin{prop}\label{normalisateur} On a les relations suivantes\,:
\begin{enumerate}
\item[ 1) ] $\n(\g^{\xi}) = \{y \in \g \; | \; [y,\xi] \in
  \z(\g^{\xi}) \}$.
\item[ 2) ] $\n(\g^{\xi}) = \g^{\xi} \oplus [\eta,\z(\g^{\xi})]$.\\
En particulier, \hbox{$\dim \n(\g^{\xi})=\dim \g^{\xi} + \dim
\z(\g^{\xi})$} et on a l'{\'e}galit{\'e} $[ \n(\g^{\xi}),\xi ] = \z(\g^{\xi})$.
\end{enumerate}
\end{prop}

Si $\u$ est un sous-espace de $\g$, on note $\u^{\perp}$ l'orthogonal
de $\u$ pour la forme de Killing $\langle , \rangle$ de $\g$. On va
d{\'e}crire l'orthogonal de certains sous-espaces. On a la proposition bien connue dont la
d{\'e}monstration est rappel{\'e}e en \cite{Charbonnel}, Lemme 5.6\,:
\begin{prop}\label{centralisateur} L'orthogonal de $\g^{\xi}$ est
   le sous-espace $[\xi,\g]$ et on a la d{\'e}composition \hbox{$\g^{\xi} \oplus [\eta,\g]=\g$}.
\end{prop}
 
Puisque ${\g^{\eta}}^{\perp}=[\eta,\g]$, la deuxi{\`e}me relation permet d'identifier le dual de $\g^{\xi}$ {\`a}
$\g^{\eta}$ via la forme de Killing.\\

Soit $\g'$ un sous-espace de $\g^{\xi}$ stable par ${\rm ad}
\rho$. Les sous-espaces $\g^{\xi}$ et $[\eta,\g']$ ont une intersection
nulle d'apr{\`e}s la proposition \ref{centralisateur} et on
s'int{\'e}resse au sous-espace \hbox{$\g^{\xi} \oplus
  [\eta,\g']$} de $\g$. Son 
orthogonal est d{\'e}crit par la
proposition suivante\,:
\begin{prop}\label{perp} Soit $\g'$ un sous-espace de $\g^{\xi}$ stable par 
${\rm ad} \rho$, alors on a\,:
$${(\g^{\xi} \oplus [\eta,\g'])}^{\perp}=[\xi,{\g'}^{\perp}] \cdot$$
\end{prop}
\dem Il est clair que $\g^{\xi}$ est contenu dans l'orthogonal de
$[\xi,\g'^{\perp}]$. Soit $u$ dans $\g'$, alors on a\,:
$$\langle [\eta,u],[\xi,{\g'}^{\perp}] \rangle=\langle
[[\eta,u],\xi],{\g'}^{\perp} \rangle=\langle -[\rho,u], {\g'}^{\perp}
\rangle=\{0 \}, $$ 
car $[\rho,u]$ appartient {\`a} $\g'$, puisque $\g'$ est stable par ${\rm ad} \rho$. On a ainsi montr{\'e} que le sous-espace $\g^{\xi} \oplus
[\eta,\g']$ est contenu dans l'orthogonal de
$[\xi,{\g'}^{\perp}]$. Calculons les dimensions des deux
sous-espaces\,:
\begin{eqnarray*}
\dim {(\g^{\xi} \oplus [\eta,\g'])}^{\perp} & = & 
\dim \g-(\dim \g^{\xi} + \dim [\eta,\g'])\\
& = & \dim \g-(\dim \g^{\xi} + (\dim \g'-\dim(\g' \cap \g^{\eta}))) \\
& = & \dim {\g'}^{\perp} - (\dim \g^{\xi} - \dim(\g' \cap \g^{\eta}))
\end{eqnarray*}
et,
\begin{eqnarray*}
\dim ([\xi,{\g'}^{\perp}]) & = & \dim {\g'}^{\perp}- \dim(
{\g'}^{\perp} \cap \g^{\xi}) \cdot\\
\end{eqnarray*}
De la d{\'e}monstration de la propostion \ref{centralisateur}, il r{\'e}sulte la d{\'e}compositon 
\hbox{$\g'=\g' \cap {\g^{\xi}}^{\perp} \oplus \g' \cap
  \g^{\eta}$}. De plus, on a\,:
\begin{eqnarray*} 
\dim(\g'\cap {\g^{\xi}}^{\perp}) & = & \dim \g-\dim({\g'}^{\perp}+\g^{\xi})\\
& = & \dim \g -(\dim{\g'}^{\perp}+\dim \g^{\xi} -\dim ({\g'}^{\perp}
\cap\g^{\xi}) \\
& = & \dim \g' -\dim \g^{\xi}+\dim ({\g'}^{\perp} \cap\g^{\xi})\cdot
\end{eqnarray*}
De cette {\'e}galit{\'e} et de la d{\'e}compostion
pr{\'e}c{\'e}dente, on d{\'e}duit la relation\,:
$$\dim\g'=(\dim \g' -\dim \g^{\xi}+\dim ({\g'}^{\perp} \cap \g^{\xi}))+\dim(\g' \cap
  \g^{\eta}),$$
ce qui donne\,: $\dim ({\g'}^{\perp} \cap \g^{\xi})=\dim \g^{\xi}-\dim(\g' \cap
  \g^{\eta})$. Par suite les deux sous-espaces $\g^{\xi} \oplus
  [\eta,\g']$ et $[\xi,{\g'}^{\perp}]$ sont de m{\^e}me dimension et la
  proposition s'ensuit.
\qed

{\`A} l'aide de la proposition pr{\'e}c{\'e}dente, on retrouve l'orthogonal de
sous-espaces 
connus. Lorsque le sous-espace $\g'$ est nul, on retrouve l'orthogonal de
  $\g^{\xi}$. D'apr{\`e}s \cite{Panyushev}, Th{\'e}or{\`e}me 2.4, on dispose de la
  d{\'e}composition $\z(\g^{\xi}) \oplus [\g,\g^{\eta}]=\g$. On en d{\'e}duit que l'orthogonal de $\z(\g^{\xi})$ est le sous-espace
  $[\g,\g^{\xi}]$. La proposition pr{\'e}c{\'e}dente
  appliqu{\'e}e {\`a} $\g'=\z(\g^{\xi})$ permet alors de d{\'e}crire
  l'orthogonal de $\n(\g^{\xi})$\,:
$${\n(\g^{\xi})}^{\perp}={(\g^{\xi} \oplus
  [\eta,\z(\g^{\xi})])}^{\perp}=[\xi,[\g^{\xi},\g]] \cdot$$ 
On a utilis{\'e} la proposition \ref{normalisateur} pour la premi{\`e}re
  {\'e}galit{\'e}. Enfin, la proposition
\ref{centralisateur} et la proposition pr{\'e}c{\'e}dente, appliqu{\'e}es {\`a} $\g'=\g^{\xi}$, donnent l'orthogonal du sous-espace 
\hbox{$\g^{\xi} \oplus
  [\eta,\g^{\xi}]$}. Ce dernier sous-espace interviendra {\`a}
plusieurs reprises dans la suite. On a\,:
$${(\g^{\xi} \oplus [\eta,\g^{\xi}])}^{\perp}=[\xi,[\xi,\g]] \cdot$$

On termine cette partie par l'introduction d'une propri{\'e}t{\'e} $(P)$\,: 
\begin{defi}\label{prop_P}
On note $W$ le sous-espace propre de la restriction de ${\rm ad
} \rho$  {\`a} $\z(\g^{\xi})$ relativement {\`a} la plus grande valeur
propre. On dira que $\xi$  {\it v{\'e}rifie la propri{\'e}t{\'e} $(P)$}
si, pour tout {\'e}l{\'e}ment non nul $v$ de $\z(\g^{\xi})$, le sous-espace $W$ est contenu dans le sous-espace $[[\eta,\g^{\xi}],v]$.
 
Il est clair que si $\xi$ v{\'e}rifie la propri{\'e}t{\'e} $(P)$, il en est de
m{\^e}me de tous les {\'e}l{\'e}ments de l'orbite de $\xi$ sous l'action du groupe
adjoint. On dira qu'une orbite nilpotente de $\g$ {\it v{\'e}rifie la
  propri{\'e}t{\'e} $(P)$} si l'un de ses repr{\'e}sentants la v{\'e}rifie. On dira
enfin que {\it $\g$ v{\'e}rifie la propri{\'e}t{\'e} $(P)$} si toutes ses orbites
nilpotentes distingu{\'e}es non r{\'e}guli{\`e}res v{\'e}rifient la propri{\'e}t{\'e} $(P)$.
\end{defi}

\section{Rappels sur l'indice d'une alg{\`e}bre de Lie et premi{\`e}res r{\'e}ductions.}
Soit $\q$ une alg{\`e}bre de Lie complexe et $\phi$ une forme lin{\'e}aire sur $\q$. On
d{\'e}signe par $\q_{\phi}$ l'ensemble des $x$ de $\q$ tels que
$\phi([\q,x])=0$. Autrement dit $\q_{\phi}=\{x \in \q \; | \; ({\rm
  ad}^{*} x) \cdot \phi=0 \}$, o{\`u} ${\rm ad}^{*} : \q \rightarrow \mathfrak{gl}(\q^{*})$ est la
repr{\'e}sentation coadjointe de $\q$. On rappelle que l'{\it indice de
$\q$}, not{\'e} ${\rm ind} \; \q$, est d{\'e}fini par\,:
$${\rm ind} \; \q= \min\limits_{\phi \in \q^{*}} \dim \q_{\phi} \;
\cdot$$  
L'indice d'une alg{\`e}bre de Lie $\q$ ainsi d{\'e}fini est un entier li{\'e} {\`a} la repr{\'e}sentation adjointe de $\q$. Une m{\'e}thode similaire, appliqu{\'e}e
{\`a} une repr{\'e}sentation de $\q$ arbitraire, permet de d{\'e}finir l'{\it
  indice d'une repr{\'e}sentation}. Soit $\rho : \q \rightarrow 
\mathfrak{gl}(V)$ une repr{\'e}sentation de $\q$. On note, de mani{\`e}re
abusive, $s\cdot v$ {\`a} la place de $\rho(s)v$, pour $s$ dans $\q$ et
$v$ dans $V$. De m{\^e}me, pour $\phi$ dans le dual $V^*$ de $V$ et pour
$s$ dans $\q$, on note $s \cdot \phi$ au lieu de $\rho^*(s) \phi$ o{\`u}
$\rho^*$ est la repr{\'e}sentation contragr{\'e}diente {\`a} $\rho$. L'entier $\; \dim V -
\max\limits_{\phi \in V^{*}} (\dim \q \cdot \phi)$ est appel{\'e} l'{\it
  indice de $V$} ou {\it l'indice du $\q$-module $V$}. On le note 
${\rm ind}(\rho,V)$ ou ${\rm ind}(\q,V)$ et il est clair que ${\rm ind}({\rm ad},\q)={\rm ind} \; \q$ au sens pr{\'e}c{\'e}dent.\\

On consid{\`e}re la forme bilin{\'e}aire {\`a} valeurs dans $V$\,:
$${\cal K}(\q,V) \; : \; \q \times V \rightarrow V \; ; \; (s,v) \mapsto
s\cdot v \; \cdot$$
En composant cette application avec un {\'e}l{\'e}ment $\phi$ de $V^{*}$, on
obtient une forme bilin{\'e}aire {\`a} valeurs dans $\C$\,:
$${\cal K}(\q,V)_{\phi} \; : \; \q \times V \rightarrow V
\stackrel{\phi}{\rightarrow} \C \; \cdot$$  
L'application ${\cal K}(\q,V)_{\phi}$ peut {\^e}tre vue comme un {\'e}l{\'e}ment
de ${\rm Hom}(\q,V^{*})$ et on v{\'e}rifie facilement l'{\'e}galit{\'e}\,: 
$${\rm ind \; }(\q,V) =\dim V - \max\limits_{\phi \in V^{*}}({\rm rang \; } {\cal K}(\q,V)_{\phi}) \cdot $$
Soit $n=\dim \q$ et $m=\dim V$. En choisissant une base sur $\q$ et $V$, on 
peut consid{\'e}rer ${\cal
K}(\q,V)$ comme une matrice de taille $m \times n$ {\`a} coefficients dans $V$, o{\`u} $V$
est identifi{\'e} {\`a} la composante de degr{\'e} $1$ de l'alg{\`e}bre sym{\'e}trique
${\cal S}(V)$. Ainsi \hbox{${\rm rang \; } {\cal K}(\q,V)=\max\limits_{\phi \in
  V^{*}}({\rm rang \; } {\cal K}(\q,V)_{\phi})$} et l'on obtient l'{\'e}galit{\'e}\,: 
\begin{eqnarray}\label{indice(q,V)}
{\rm ind \; }(\q,V) =  \dim V -{\rm rang \; } {\cal K}(\q,V).
\end{eqnarray}
Pour $V=\q$, on note ${\cal K}(\q)$ au lieu de ${\cal K}(\q,\q)$ et il
vient\,:
\begin{eqnarray}\label{indice(q)}
{\rm ind \; }(\q)  =  \dim \q -{\rm rang \; } {\cal K}(\q).
\end{eqnarray}
~\\

On s'int{\'e}resse maintenant {\`a} l'alg{\`e}bre de Lie $\n(\g^{\xi})$ et on
commence par examiner l'identit{\'e} du th{\'e}or{\`e}me \ref{intro1}. On
choisit une base $e_{1},\ldots,e_{n}$ de $\g^{\xi}$ telle que
$e_{1},\ldots,e_{m}$ est une base de $\z(\g^{\xi})$, $m \leq
n$. Alors la famille $\{e_{1},\ldots,e_{n},[\eta,e_{1}],\ldots,[\eta,e_{m}]\}$ est
une base de $\n(\g^{\xi})$ d'apr{\`e}s la proposition \ref{normalisateur}. Dans cette base, la matrice ${\cal
  K}(\n(\g^{\xi}))$ est de la forme 
$$\left[
\begin{array}{ccc} 
0 & 0 & \mathfrak{D} \\
0 & \mathfrak{C} & \mathfrak{E} \\
-\mathfrak{D}^{t} & -\mathfrak{E}^{t} & \mathfrak{F} 
\end{array}
\right],$$
o{\`u} les matrices carr{\'e}es $\mathfrak{D}$ et $\mathfrak{C}$ sont
respectivement d'ordre $m$ et $n-m$. On reconnait certains blocs de cette matrice. Ainsi 
${\cal K}(\n(\g^{\xi}),\g^{\xi})=
\left[
\begin{array}{ccc} 
0 & 0 & \mathfrak{D} \\
0 & \mathfrak{C} & \mathfrak{E}
\end{array}
\right]$ et 
${\cal K}(\g^{\xi})=
\left[
\begin{array}{cc} 
0 & 0 \\
0 & \mathfrak{C} 
\end{array}
\right]$. D'apr{\`e}s (\ref{indice(q)}), on a la relation 
$${\rm ind \;} \n(\g^{\xi})=\dim \n(\g^{\xi})-{\rm rang \; } {\cal K}(
\n(\g^{\xi})) \cdot$$
Or la structure de la matrice ${\cal K}(\n(\g^{\xi}))$ montre que
l'on a l'in{\'e}galit{\'e} 
$${\rm rang \; } {\cal K}(\n(\g^{\xi})) \leq {\rm rang \;} \mathfrak{C} +2
\dim \z(\g^{\xi}) \cdot$$ 
Par ailleurs, puisque 
${\cal K}(\g^{\xi})=
\left[
\begin{array}{cc} 
0 & 0 \\
0 & \mathfrak{C} \\
\end{array}
\right]$, on a l'{\'e}galit{\'e}\,:
$${\rm rang \; } \mathfrak{C}={\rm rang \; } {\cal K}(\g^{\xi})=\dim \g^{\xi}-{\rm ind
  \;}\g^{\xi} \cdot$$
Par suite, il vient\,: 
$${\rm ind \;} \n(\g^{\xi}) \geq \dim
  \n(\g^{\xi})-\dim \g^{\xi}+{\rm ind \;}\g^{\xi}-2 \dim
  \z(\g^{\xi}) \cdot$$ 
D'apr{\`e}s \cite{Charbonnel}, Th{\'e}or{\`e}me 5.5, on a la relation \hbox{${\rm
  ind \;}\g^{\xi}={\rm rg \;} \g$}. Par ailleurs, la proposition
  \ref{normalisateur}, 2) donne \hbox{$\dim \n(\g^{\xi})=\dim \g^{\xi}+\dim  \z(\g^{\xi})$}. On obtient finalement l'in{\'e}galit{\'e}\,: 
$${\rm ind \;} \n(\g^{\xi}) \geq {\rm rg \;} \g - \dim
\z(\g^{\xi}) \cdot$$ 
Voyons maintenant pourquoi le th{\'e}or{\`e}me \ref{intro2} implique le
th{\'e}or{\`e}me \ref{intro1}. Il s'agit de prouver la proposition suivante\,:
\begin{prop}\label{implique} Si l'{\'e}galit{\'e} 
${\rm ind \;} (\n(\g^{\xi}),\g^{\xi}) = {\rm rg \;} \g - \dim
\z(\g^{\xi})$
est satisfaite, alors l'{\'e}galit{\'e}
${\rm ind \;} \n(\g^{\xi}) = {\rm rg \;} \g - \dim
\z(\g^{\xi})$
est satisfaite.
\end{prop}
\dem
 On suppose le th{\'e}or{\`e}me \ref{intro2} d{\'e}montr{\'e}. D'apr{\`e}s ce qui pr{\'e}c{\`e}de, il suffit d'obtenir l'in{\'e}galit{\'e}\,:
$${\rm ind \;} \n(\g^{\xi}) \leq {\rm rg \;} \g - \dim
\z(\g^{\xi}) \cdot$$ 
Rappelons le r{\'e}sultat suivant, d{\'e}montr{\'e} dans \cite{Panyushev}, Th{\'e}or{\`e}me 1.4\,:
\begin{lemme}\label{lemme_indice}
Soit $\q$ un id{\'e}al d'une alg{\`e}bre de Lie $\widetilde{\q}$. Alors\\
$${\rm ind \ } \q +{\rm ind \ } \widetilde{\q} \leq
\dim(\displaystyle{\widetilde{\q} / \q }) + 2 {\rm \
  ind}(\widetilde{\q},\q) \cdot$$.
\end{lemme}

En appliquant ce lemme {\`a} l'id{\'e}al $\g^{\xi}$ de $\n(\g^{\xi})$, on
obtient\,:
$${\rm ind \;} \g^{\xi}+{\rm ind \;} \n(\g^{\xi}) \leq \dim
(\displaystyle{\n(\g^{\xi}) /  \g^{\xi}}) + 2 {\rm \ ind \;}
(\n(\g^{\xi}),\g^{\xi}), $$ 
ce qui donne, en utilisant de nouveau les {\'e}galit{\'e}s $\dim
\n(\g^{\xi})=\dim \g^{\xi} + \dim \z(\g^{\xi})$ et \hbox{${\rm ind \;}
\g^{\xi}={\rm rg \ } \g$},
$${\rm ind \;} \n(\g^{\xi}) \leq 2 {\rm \ ind \;} (\n(\g^{\xi}),\g^{\xi}) -
({\rm rg \ } \g - \dim \z(\g^{\xi})) \cdot$$
Le th{\'e}or{\`e}me \ref{intro2} entraine alors l'in{\'e}galit{\'e} souhait{\'e}e.
\qed

On s'int{\'e}resse d{\'e}sormais {\`a} l'identit{\'e} du th{\'e}or{\`e}me \ref{intro2}. On
montre la proposition suivante\,:
\begin{prop}\label{rang} Le th{\'e}or{\`e}me \ref{intro2} est {\'e}quivalent {\`a}
  l'assertion suivante\,: la matrice 
$$\left[
\begin{array}{c} 
\mathfrak{D} \\
\mathfrak{E}
\end{array}
\right]=
\left[
\begin{array}{ccc} 
[[\eta,e_{1}],e_{1}] & \cdots & [[\eta,e_{m}],e_{1}] \\
\vdots &  & \vdots \\
\left[ [\eta,e_{1}],e_{m} \right] & \cdots & [[\eta,e_{m}],e_{m}]\\
\vdots &  & \vdots \\
\left[ [\eta,e_{1}],e_{n} \right] & \cdots & [[\eta,e_{m}],e_{n}]
\end{array}
\right]$$
de taille $n \times m$, {\`a}
  coefficients dans ${\cal S}(\g^{\xi})$,
est de rang maximal {\'e}gal {\`a} $m=\dim \z(\g^{\xi})$.
\end{prop}

Avant de d{\'e}montrer la proposition, on {\'e}tend encore un peu la
d{\'e}finition de l'indice. Ceci permettra d'interpr{\'e}ter g{\'e}om{\'e}triquement 
la matrice $\left[
\begin{array}{c} 
\mathfrak{D} \\
\mathfrak{E}
\end{array}
\right]$. Soit
$\mathfrak{q}$ une alg{\`e}bre de Lie complexe, $V$ et $V'$ deux espaces
vectoriels de dimension finie sur $\C$ et soit $\rho\,: \q
\rightarrow {\rm L}(V,V')$ une application lin{\'e}aire de $\q$ dans l'espace des applications lin{\'e}aires de $V$ dans
$V'$. Pour $\phi$ dans $(V')^{*}$, on note \hbox{$\q_{\phi}=\{s \in \q ~ | ~
\phi(\rho(s)v)=0, \ \forall v \in V \}$} le \guillemotleft
stabilisateur de $\phi$\guillemotright \  et
on note $\q \cdot_{\rho} \phi$ l'image dans $V^{*}$ de 
l'application qui {\`a}
$s$ dans $\q$ associe la forme lin{\'e}aire $v \mapsto -\phi(\rho(s)v)$
d{\'e}finie sur $V$. Lorsqu'il n'y a pas d'ambiguit{\'e}
on omet l'indice $\rho$. On pose, par
analogie avec l'indice, 
$$r(\q,V,V')=\dim V - 
\max\limits_{\phi \in (V')^{*}} (\dim \q \cdot_{\rho} \phi) \cdot$$
Notons que $\rho$ n'est pas un morphisme
d'alg{\`e}bres de Lie en g{\'e}n{\'e}ral; l'espace ${\rm L}(V,V')$ n'est m{\^e}me
pas une alg{\`e}bre de Lie! En revanche, si $V=V'$ et si $\rho$ est une repr{\'e}sentation de $\q$ dans $V$, on retrouve
l'indice ${\rm ind \; }(\q,V)$ du $\q$-module $V$.
On consid{\`e}re, toujours par analogie avec l'indice, la forme bilin{\'e}aire {\`a} valeurs dans $V'$\,:
$${\cal K}(\q,V,V') \; : \; \q \times V \rightarrow V' \; ; \; (s,v) \mapsto
-\rho(s)v \; \cdot$$
En composant cette application avec un {\'e}l{\'e}ment $\phi$ de $(V')^{*}$, on
obtient une forme bilin{\'e}aire {\`a} valeurs dans $\C$\,:
$${\cal K}(\q,V,V')_{\phi} \; : \; \q \times V \rightarrow V'
\stackrel{\phi}{\rightarrow} \C \; \cdot$$  
On v{\'e}rifie facilement l'{\'e}galit{\'e}\,: 
$$r(\q,V,V') =\dim V - \max\limits_{\phi \in (V')^{*}}({\rm rang \; }
{\cal K}(\q,V,V')_{\phi}) \cdot $$
Cette {\'e}criture permet de voir que l'ensemble des formes lin{\'e}aires
$\phi$ de $(V')^*$ telles que 
$\dim \q \cdot_{\rho} \phi= \dim V - r(\q,V,V')$ est un ouvert dense de
$(V')^{*}$. 
En choisissant des bases sur $\q$ et $V$, on 
peut consid{\'e}rer ${\cal
K}(\q,V,V')$ comme une matrice de taille $\dim V \times \dim \q$ {\`a} coefficients dans $V'$, o{\`u} $V'$
est identifi{\'e} {\`a} la composante de degr{\'e} $1$ de l'alg{\`e}bre sym{\'e}trique
${\cal S}(V')$. Ainsi, on obtient l'{\'e}galit{\'e}\,: 
$$r(\q,V,V') =  \dim V -{\rm rang \; } {\cal K}(\q,V,V') \cdot$$

Prouvons maintenant la proposition \ref{rang}\,:\\
\\ 
\dem D'apr{\`e}s la relation (\ref{indice(q,V)}), on a\,:
$${\rm ind \;} (\n(\g^{\xi}),\g^{\xi})=\dim \g^{\xi}-{\rm rang \; } {\cal K}(
\n(\g^{\xi}),\g^{\xi}) \cdot$$
Or la structure de la matrice ${\cal K}(\n(\g^{\xi}),\g^{\xi})$ montre que
l'on a l'in{\'e}galit{\'e}\,: 
$${\rm rang \; } {\cal K}(\n(\g^{\xi}),\g^{\xi}) \leq {\rm rang \;}
\mathfrak{C} + \dim \z(\g^{\xi}) \cdot$$
En utilisant la relation, d{\'e}j{\`a} vue, \hbox{${\rm rang \;} \mathfrak{C}=\dim
\g^{\xi}-{\rm ind \;} \g^{\xi}=\dim \g^{\xi}-{\rm rg \;} \g$}, 
on obtient l'in{\'e}galit{\'e}\,:
$${\rm ind \;} (\n(\g^{\xi}),\g^{\xi}) \geq {\rm rg \;} \g - \dim
\z(\g^{\xi}) \cdot$$

Il y a {\'e}galit{\'e}
dans la relation pr{\'e}c{\'e}dente si, et seulement si, la condition 
\hbox{${\rm rang \; } {\cal K}(\n(\g^{\xi}),\g^{\xi}) 
= {\rm rang \;} \mathfrak{C} + \dim \z(\g^{\xi})$} 
est satisfaite, autrement dit si, et seulement si, le rang de la
matrice 
$\left[
\begin{array}{cc} 
0 & \mathfrak{D} \\
 \mathfrak{C} & \mathfrak{E}
\end{array}
\right]$, 
{\`a} coefficients dans ${\cal S}(\g^{\xi})$, est {\'e}gal {\`a} ${\rm rang \;} \mathfrak{C} + \dim
\z(\g^{\xi})$. Par suite, si l'{\'e}galit{\'e} du th{\'e}or{\`e}me \ref{intro2} est
satisfaite, n{\'e}cessairement la matrice 
$\left[
\begin{array}{c} 
\mathfrak{D} \\
\mathfrak{E}
\end{array}
\right]$ 
est de rang maximal {\'e}gal {\`a} $m=\dim \z(\g^{\xi})$.\\

R{\'e}ciproquement, supposons que le rang de la matrice 
$\left[
\begin{array}{c} 
\mathfrak{D} \\
\mathfrak{E}
\end{array}
\right]$
soit maximal, {\'e}gal {\`a} \hbox{$m=\dim \z(\g^{\xi})$}, et montrons la relation du
th{\'e}or{\`e}me \ref{intro2}. Soit $\rho$ l'application lin{\'e}aire de
$\z(\g^{\xi})$ dans ${\rm L}([\eta,\g^{\xi}],\g^{\xi})$ donn{\'e}e par la relation\,:
$$\rho(s)v=[s,v] \in \g^{\xi},$$
pour $s$ dans $\z(\g^{\xi})$ et $v$ dans $[\eta,\g^{\xi}]$. Puisque
$\z(\g^{\xi})$ et $\g^{\eta}$ ont une intersection nulle, il existe
un suppl{\'e}mentaire $\mathfrak{r}$ de $\g^{\eta} \cap \g^{\xi}$ dans
$\g^{\xi}$ contenant $\z(\g^{\xi})$. Posons \hbox{$r=\dim
  \mathfrak{r}$}. Alors on a $m \leq r \leq n$. On peut supposer
que les {\'e}l{\'e}ments $e_{1},\ldots,e_{r}$ forment une base de
$\mathfrak{r}$ et que les {\'e}l{\'e}ments $e_{r+1},\ldots,e_{n}$ forment une
base de $\g^{\eta} \cap \g^{\xi}$. La famille $\{[\eta,e_{1}],\ldots,[\eta,e_{r}]\}$ est
libre dans $[\eta,\g^{\xi}]$ et de cardinal \hbox{$r=\dim\g^{\xi}-\dim
\g^{\eta} \cap \g^{\xi}=\dim([\eta,\g^{\xi}])$}; c'est donc une base de
$[\eta,\g^{\xi}]$. Dans les bases $\{e_{1},\ldots,e_{m}\}$ de
$\z(\g^{\xi})$ et $\{[\eta,e_{1}],\ldots,[\eta,e_{r}]\}$ de
$[\eta,\g^{\xi}]$, la matrice 
${\cal K}(\z(\g^{\xi}),[\eta,\g^{\xi}],\g^{\xi})$ {\`a} 
coefficients dans ${\cal S}(\g^{\xi})$ est donn{\'e}e par 
$${\cal K}(\z(\g^{\xi}),[\eta,\g^{\xi}],\g^{\xi})=\left[
\begin{array}{ccc}
[e_{1},[\eta,e_{1}]] & \cdots & [e_{m},[\eta,e_{1}]] \\
\vdots &  & \vdots \\
\left[ e_{1},[\eta,e_{m}] \right] & \cdots & [e_{m},[\eta,e_{m}]]\\
\left[ e_{1},[\eta,e_{m+1}] \right] & \cdots & [e_{m},[\eta,e_{m+1}]]\\
\vdots &  & \vdots \\
\left[ e_1,[\eta,e_{r}] \right] & \cdots & [e_m,[\eta,e_{r}]]
\end{array}
\right] \cdot$$
On dispose des relations\,:
$$[[\eta,e_{j}],e_{i}]=[[\eta,e_{i}],e_{j}]=-[e_j,[\eta,e_{i}]],$$
pour $i=1,\ldots,n$ et $j=1,\ldots,m$.
Puisque les vecteurs $e_{r+1},\ldots,e_{n}$ sont dans le
centralisateur $\g^{\eta}$ de $\eta$, les crochets $[[\eta,e_{j}],e_{i}]$ sont nuls, pour
$i=r+1,\ldots,n$ et $j=1,\ldots,m$, et la matrice $\left[
\begin{array}{c} 
\mathfrak{D} \\
\mathfrak{E}
\end{array}
\right]$ de la proposition est donn{\'e}e par\,:
$$\left[
\begin{array}{c} 
\mathfrak{D} \\
\mathfrak{E}
\end{array}
\right]=
\left[
\begin{array}{ccc}
[[\eta,e_{1}],e_{1}] & \cdots & [[\eta,e_{m}],e_{1}] \\
\vdots &  & \vdots \\
\left[ [\eta,e_{1}],e_{m} \right] & \cdots & [[\eta,e_{1}],e_{m}]\\
\left[ [\eta,e_{1}],e_{m+1} \right] & \cdots & [[\eta,e_{m}],e_{m+1}]\\
\vdots &  & \vdots \\
\left[ [\eta,e_{1}],e_{r} \right] & \cdots & [[\eta,e_{m}],e_{r}]\\
0 & \cdots & 0 \\
\vdots &  & \vdots \\
0 & \cdots & 0
\end{array}
\right] \cdot$$ 
Par suite les deux
matrices ${\cal K}(\z(\g^{\xi}),[\eta,\g^{\xi}],\g^{\xi})$ et 
$\left[
\begin{array}{c} 
\mathfrak{D} \\
\mathfrak{E}
\end{array}
\right]$ ont le m{\^e}me rang. D'apr{\`e}s l'hypoth{\`e}se sur le rang de la
matrice $\left[
\begin{array}{c} 
\mathfrak{D} \\
\mathfrak{E}
\end{array}
\right]$, on en d{\'e}duit que l'entier 
$\max\limits_{\phi \in (\g^{\xi})^{*}} \dim \z(\g^{\xi}) \cdot_{\rho} \phi$ 
est {\'e}gal {\`a} $\dim \z(\g^{\xi})$. L'ensemble 
$$\Omega'=\{ \phi \in (\g^{\xi})^{*} \ | \ \dim \z(\g^{\xi}) \cdot_{\rho} \phi
= \dim \z(\g^{\xi}) \}$$
est donc un ouvert dense de $(\g^{\xi})^{*}$. Par ailleurs, de l'{\'e}galit{\'e} ${\rm ind \ } \g^{\xi}=
{\rm rg \ } \g$, il r{\'e}sulte que l'ensemble 
$$\Omega=\{ \phi \in (\g^{\xi})^{*} \ | \ \dim \g^{\xi} \cdot \phi
= \dim \g^{\xi}- {\rm rg \ } \g \}$$ 
est un ouvert dense de $(\g^{\xi})^{*}$, o{\`u} $\g^{\xi} \cdot \phi$ d{\'e}signe l'orbite coadjointe
de $\phi$. L'intersection $\Omega \cap \Omega'$ est donc un ouvert non vide
de $(\g^{\xi})^{*}$. Soit alors $\lambda$ une forme lin{\'e}aire sur
$\g^{\xi}$ appartenant {\`a} cette intersection. On consid{\`e}re l'orbite
$\n(\g^{\xi}) \cdot \lambda$ de $\lambda$ sous l'action naturelle de $\n(\g^{\xi})$ dans
$(\g^{\xi})^{*}$. Autrement dit,
$$\n(\g^{\xi}) \cdot \lambda =\{ v \in \g^{\xi} \mapsto -\lambda([s,v]) \ | \ s \in
\n(\g^{\xi}) \} \cdot$$
Les sous-espaces $\g^{\xi}$ et $[\eta,\g^{\xi}]$ ont une intersection
nulle d'apr{\`e}s la proposition \ref{centralisateur} et on pose 
$E=\g^{\xi} \oplus [\eta,\g^{\xi}]$. En identifiant $(\g^{\xi})^{*}$,
respectivement $([\eta,\g^{\xi}])^*$,  
au sous-espace de $E^{*}$ form{\'e} des formes lin{\'e}aires nulles sur
$[\eta,\g^{\xi}]$, respectivement $\g^{\xi}$, on obtient la
d{\'e}composition suivante\,:
$$E^*=(\g^{\xi})^* \oplus ([\eta,\g^{\xi}])^* \cdot$$
Par suite, les ensembles $\g^{\xi} \cdot \lambda$ et
$\z(\g^{\xi}) \cdot_{\rho} \lambda$ peuvent {\^e}tre vus comme deux sous-espaces 
de $E^*$ dont l'intersection est nulle. En
effet le premier est contenu dans $(\g^{\xi})^*$ et le deuxi{\`e}me dans
$([\eta,\g^{\xi}])^*$. On va montrer qu'il existe un
isomorphisme d'espaces vectoriels entre $\n(\g^{\xi}) \cdot \lambda$
et $\g^{\xi} \cdot \lambda \oplus \z(\g^{\xi}) \cdot_{\rho}
\lambda$. Soit $s$ dans $\n(\g^{\xi})$. Selon la
d{\'e}composition 2) de la proposition \ref{normalisateur}, il s'{\'e}crit de 
mani{\`e}re unique $s=s_{1}+[\eta,s_{2}]$, avec $s_{1}$ dans
$\g^{\xi}$ et $s_{2}$ dans $\z(\g^{\xi})$. On pose 
$$\Phi(s \cdot \lambda)=s_{1} \cdot \lambda +s_{2} \cdot_{\rho} \lambda
\cdot$$
L'application $\Phi$ d{\'e}finit ainsi une application lin{\'e}aire de $\n(\g^{\xi})
\cdot \lambda$ dans $\g^{\xi} \cdot \lambda \oplus \z(\g^{\xi}) \cdot_{\rho}
\lambda$. C'est clairement un isomorphisme d'espaces vectoriels, dont
l'inverse est donn{\'e} par\,:
$$s_{1} \cdot \lambda +s_{2} \cdot_{\rho} \lambda  \longmapsto
(s_{1}+[\eta,s_{2}]) \cdot \lambda \cdot$$ 
De cet isomorphisme, on d{\'e}duit une {\'e}galit{\'e} sur les dimensions,
$$\dim \n(\g^{\xi}) \cdot \lambda=\dim \g^{\xi} \cdot \lambda + \dim
\z(\g^{\xi}) \cdot_{\rho} \lambda,$$
car les deux sous-espaces $ \g^{\xi} \cdot \lambda$ et $\z(\g^{\xi})
\cdot_{\rho} \lambda$ ont une intersection nulle. D'o{\`u}, puisque $\lambda$
appartient {\`a} l'intersection $\Omega \cap \Omega'$,
$$\dim \n(\g^{\xi}) \cdot \lambda=(\dim \g^{\xi}-{\rm rg \ } \g) +
\dim \z(\g^{\xi}) \cdot$$
Par suite, on obtient 
$$\max\limits_{\phi \in (\g^{\xi})^{*}} (\dim \n(\g^{\xi}) \cdot \phi)
\geq \dim \g^{\xi}-{\rm rg \ } \g +
\dim \z(\g^{\xi}),$$
 ce qui donne encore 
$${\rm ind \;} (\n(\g^{\xi}),\g^{\xi}) \leq {\rm rg \;} \g - \dim
\z(\g^{\xi}) \cdot$$
L'in{\'e}galit{\'e} oppos{\'e}e {\'e}tant d{\'e}j{\`a} connue, on en d{\'e}duit le th{\'e}or{\`e}me
\ref{intro2}. L'{\'e}quivalence de la proposition est ainsi d{\'e}montr{\'e}e.
\qed

Notons tout d'abord que la condition de la proposition pr{\'e}c{\'e}dente est en particulier
v{\'e}rifi{\'e}e si la matrice $\mathfrak{D}$ est non singuli{\`e}re. Dire que la
matrice $\mathfrak{D}$ est non singuli{\`e}re signifie que l'indice ${\rm ind
  \;} (\n(\g^{\xi}),\z(\g^{\xi}))$ du $\n(\g^{\xi})$-module
$\z(\g^{\xi})$ est nul, ce qui se traduit encore en disant que le
groupe $N_{\g}(\xi)$ a une orbite ouverte dans $\z(\g^{\xi})^*$. On
retrouve ainsi la condition de D. Panyushev vue en introduction qui lui
permet de conclure dans un certain nombre de cas. Pr{\'e}cis{\'e}ment, il
prouve en \cite{Panyushev} que la matrice $\mathfrak{D}$ est non singuli{\`e}re dans les cas suivants\,: \begin{itemize}
\item[1)] si $\g$ est isomorphe {\`a} l'une des alg{\`e}bres de Lie
  $\mathfrak{sl}_{n},\mathfrak{so}_{2n+1}$ ou $\mathfrak{sp}_{n}$, (Th{\'e}or{\`e}me 4.7 (i)),
\item[2)] pour certains {\'e}l{\'e}ments nilpotents de $\g$, si
  $\g$ est isomorphe {\`a} $\mathfrak{so}_{2n}$, (Th{\'e}or{\`e}me 4.7 (ii)),
\item[3)] si $\dim \z(\g^{\xi}) \leq 2$, (Th{\'e}or{\`e}me
  4.7 (iii)), 
\item[4)] si $\xi$ est un {\'e}l{\'e}ment nilpotent r{\'e}gulier de $\g$, (Corollaire 5.6).
\end{itemize}
Notons que dans les trois premiers cas, D. Panyushev montre une
propri{\'e}t{\'e} plus forte que celle de l'orbite ouverte, {\`a} savoir\,: le
groupe $N_{\g}(\xi)$ a un nombre fini d'orbites dans $\z(\g^{\xi})$. 
Remarquons enfin qu'il existe des cas o{\`u} la matrice $\mathfrak{D}$ est
singuli{\`e}re; par exemple, si $\xi$ est un {\'e}l{\'e}ment nilpotent
sous-r{\'e}gulier de $\mathfrak{so}_{8}$, on peut montrer
(\cite{Panyushev}, partie 4) que le groupe $N_{\g}(\xi)$ n'a pas
d'orbite ouverte dans $\z(\g^{\xi})^*$.\\ 

Soit $\xi$ un {\'e}l{\'e}ment nilpotent non distingu{\'e} de $\g$. En
raisonnant comme dans \cite{Tauvel}, Chapitre XVII, preuve de la
  proposition 5.10, on montre qu'il existe une sous-alg{\`e}bre de Lie semi-simple $\mathfrak{t}$ de $\g$,
contenant $\xi$, et telle que les trois conditions suivantes soient
v{\'e}rifi{\'e}es\,: 
\begin{enumerate}
\item[1)] $\mathfrak{t}$ est l'alg{\`e}bre d{\'e}riv{\'e}e du centralisateur d'un
  {\'e}l{\'e}ment semi-simple de $\g$,
\item[2)] $\z(\g^{\xi}) \subset \z({\mathfrak{t}}^{\xi}) \subset
  {\mathfrak{t}}^{\xi} \subset \g^{\xi}$,
\item[3)] $\xi$ est distingu{\'e} dans $\mathfrak{t}$.
\end{enumerate}

\begin{prop}\label{distingue} On suppose que la relation 
$${\rm ind}(\n(\mathfrak{t}^{\xi}),\mathfrak{t}^{\xi})={\rm rg \ }
\mathfrak{t} - \dim \z(\mathfrak{t}^{\xi})$$
est satisfaite. Alors la relation
$${\rm ind}(\n(\g^{\xi}),\g^{\xi})={\rm rg \ }
\g - \dim \z(\g^{\xi})$$
est satisfaite.
\end{prop}
\dem Soit $\{\xi, \rho',\eta'\}$ un $\mathfrak{sl}_{2}$-triplet dans
$\mathfrak{t}$ contenant $\xi$. On construit une base $e_{1}, \ldots, e_{n}$ de
$\g^{\xi}$ telle que $e_{1},\ldots,e_{m}$ soit une base de
$\z(\g^{\xi})$, $e_{1},\ldots,e_{m'}$ une base de
$\z({\mathfrak{t}}^{\xi})$ et $e_{1},\ldots,e_{n'}$ une base de
${\mathfrak{t}}^{\xi}$, avec $m \leq m' \leq n' \leq n$. Comme la
relation \hbox{${\rm ind}(\n(\mathfrak{t}^{\xi}),\mathfrak{t}^{\xi})={\rm rg \ }
\mathfrak{t} - \dim \z(\mathfrak{t}^{\xi})$} est satisfaite, la
proposition \ref{rang} assure qu'il existe un vecteur $v$ de ${\mathfrak{t}}^{\eta}$
tel que la matrice $(\langle v, [[\eta',e_{i}],e_{j}] \rangle)_{1 \leq
  i \leq n' \atop 1 \leq j \leq m'}$ est de rang maximal, car
${\mathfrak{t}}^{\eta}$ s'identifie au dual de ${\mathfrak{t}}^{\xi}$ via la forme
de Killing. On suppose
par l'absurde que la matrice $(\langle v, [[\eta',e_{i}],e_{j}] \rangle)_{1 \leq
  i \leq n \atop 1 \leq j \leq m}$ n'est pas de rang maximal. Alors il
existe des complexes $a_{1},\ldots,a_{m}$ non tous nuls qui satisfont
l'{\'e}galit{\'e}\,:
$$\langle v, [[\eta',e_{i}], \sum\limits_{j=1}^{m}a_{j}e_{j}]
\rangle=0,$$
pour $i=1, \ldots,n$. On a donc aussi la relation\,:
$$\langle v, [[\eta',e_{i}], \sum\limits_{j=1}^{m}a_{j}e_{j}]
\rangle=0,$$
pour $i=1, \ldots,n'$, car $n' \leq n$. En posant $a_{j}=0$ pour 
$j=m+1,\ldots,m'$, on obtient encore\,:
$$\langle v, [[\eta',e_{i}], \sum\limits_{j=1}^{m'}a_{j}e_{j}]
\rangle=0,$$
pour $i=1, \ldots,n'$. Ceci contredit le fait que la matrice $(\langle v, [[\eta',e_{i}],e_{j}] \rangle)_{1 \leq
  i \leq n' \atop 1 \leq j \leq m'}$ est de rang maximal. Par suite, la
proposition \ref{rang} donne la relation  
${\rm ind}(\n(\g^{\xi}),\g^{\xi})={\rm rg \ }
\g - \dim \z(\g^{\xi})$, car cette derni{\`e}re ne d{\'e}pend pas du choix d'un $\mathfrak{sl}_{2}$-triplet $\{\xi,
\rho,\eta\}$ dans $\g$ contenant $\xi$.
\qed

Effectuons enfin une derni{\`e}re r{\'e}duction. Une alg{\`e}bre de Lie complexe semi-simple est somme directe d'id{\'e}aux
simples et il est clair qu'il suffit de d{\'e}montrer les th{\'e}or{\`e}mes
\ref{intro1} et \ref{intro2} pour chaque composante simple de
$\g$. On peut supposer d{\'e}sormais que $\g$ est une alg{\`e}bre de Lie
complexe simple. 

\section{D{\'e}monstration des th{\'e}or{\`e}mes \ref{intro1} et \ref{intro2} dans le cas classique.}

On suppose que $\g$ est une alg{\`e}bre de Lie simple classique,
c'est-{\`a}-dire on suppose que $\g=\mathfrak{sl}_{n}$,
$\mathfrak{so}_{2n+1}$, $\mathfrak{sp}_{2n}$ ou $\mathfrak{so}_{2n}$. Soit $\xi$ un {\'e}l{\'e}ment
nilpotent de $\g$ et soit $\{\xi,\rho,\eta\}$ un $\mathfrak{sl}_{2}$-triplet contenant $\xi$. On consid{\`e}re le sous-espace $\z'$
de $\g$ engendr{\'e} par les puissances de $\xi$. Autrement dit, on a\,: 
$$\z'=\{ \textrm{polyn{\^o}mes en } \xi \} \cap \g \cdot$$ 
Alors clairement $\z'$ est une sous-alg{\`e}bre de $\z(\g^{\xi})$ stable
par ${\rm ad} \rho$ et on dispose, d'apr{\`e}s \cite{Kurtzke}, preuve de
la proposition 3.2 et proposition 3.3, des r{\'e}sultats suivants\,:
\begin{lemme}\label{etude_z'} Si $\g=\mathfrak{sl}_{n}$, le sous-espace
  $\z'$ est engendr{\'e} par les
  puissances de $\xi$. Si $\g=\mathfrak{so}_{n}$ ou
  $\mathfrak{sp}_{2n}$, le sous-espace $\z'$ est engendr{\'e} 
par les puissances impaires
  de $\xi$. De plus, on a $\z'=\z(\g^{\xi})$ sauf dans le cas o{\`u}
  $\g=\mathfrak{so}_{2n}$ et o{\`u} la partition de $2n$ associ{\'e}e {\`a} $\xi$ ne
  contient que deux parties. Dans ce cas, le sous-espace $\z'$ est de
  codimension 1 dans $\z(\g^{\xi})$. 
\end{lemme}
On reprend les notations de la partie pr{\'e}c{\'e}dente. On va montrer
\guillemotleft directement\guillemotright, {\`a} l'aide de calculs
explicites, que la matrice 
$\left[
\begin{array}{c}
\mathfrak{D}\\
\mathfrak{E}
\end{array}
\right]$ est de rang maximal. Dans le cas o{\`u} \hbox{$\z'=\z(\g^{\xi})$}, la
proposition \ref{indice_nul} donne un r{\'e}sultat plus fort. Commen{\c c}ons par quelques calculs pr{\'e}liminaires. On dispose des relations suivantes\,:
\begin{eqnarray}\label{rho}
[\xi^{k},\eta] = \xi^{k} \eta - \eta \xi^{k} 
=\sum\limits_{\alpha+\beta=k-1} \xi^{\alpha} \rho \xi^{\beta} 
\textrm{ et } \left[ \rho , \xi^{i} \right] =2i \cdot \xi^{i},
\end{eqnarray}
\begin{eqnarray}\label{eta}
\left[ \left[\xi^{k},\eta \right],\xi^{i} \right] =  
\sum\limits_{\alpha+\beta=k-1}\xi^{\alpha} \rho  \xi^{\beta+i}- 
\sum\limits_{\alpha+\beta=k-1}\xi^{\alpha+i} \rho  \xi^{\beta}=
\sum\limits_{\alpha+\beta=k-1}  \xi^{\alpha} [\rho,\xi^{i}]
\xi^{\beta}
=2 ki \cdot \xi^{k+i-1} .
\end{eqnarray} 

L'alg{\`e}bre $\n(\g^{\xi})$ agit sur le sous-espace $\z(\g^{\xi})$ par la
repr{\'e}sentation adjointe et l'on dispose de la proposition suivante\,:
\begin{prop}\label{indice_nul}
On suppose $\z'=\z(\g^{\xi})$. Alors l'indice ${\rm ind \ }(\n(\g^{\xi}),\z(\g^{\xi}))$ est nul. En particulier, la
relation 
$${\rm ind \ }(\n(\g^{\xi}),\g^{\xi})={\rm rg \ }\g-\dim
\z(\g^{\xi})$$
est satisfaite. 
\end{prop} 
\dem Notons $d$ le degr{\'e} du polyn{\^o}me minimal de $\xi$. On suppose tout
d'abord $\g=\mathfrak{sl}_{n}$. Posons $e_{i}=\xi^{i}$, pour
$i=1,\ldots,d-1$, alors les {\'e}l{\'e}ments $e_{1},\ldots,e_{d-1}$ forment une base
de $\z'$. Notons $\mathfrak{D}=(\mathfrak{D}_{ij})_{1 \leq i,j \leq
  d-1}$ la matrice carr{\'e}e d'ordre $d-1$ donn{\'e}e par\,:
$$\mathfrak{D}_{ij}=[[\eta,e_{j}],e_{i}],$$
pour $i$ et $j$ dans $\{1,\ldots,d-1\}$. D'apr{\`e}s la relation
(\ref{eta}), on a 
\begin{equation*}
\mathfrak{D}_{ij}=
\begin{cases}
0 & {\rm \ si \ } i+j > d,\\
-2ij \xi^{i+j-1} & {\rm \ sinon.} 
\end{cases}
\end{equation*}
En particulier on connait $\mathfrak{D}$ sur l'anti-diagonale. On
a $\mathfrak{D}_{ij}=-2ij \xi^{d-1}$ si $i+j=d$. Soit $\phi$ une forme lin{\'e}aire sur
$\z(\g^{\xi})$. De ce qui pr{\'e}c{\`e}de, il r{\'e}sulte que le d{\'e}terminant de la
matrice $\mathfrak{D}$ {\'e}valu{\'e} en $\phi$ vaut, au signe pr{\`e}s\,:
$$2^{d-1} \times (d-1)!^2 \times (\langle \phi,\xi^{d-1} \rangle
)^{d-1} \cdot$$
C'est donc un {\'e}l{\'e}ment non nul de
${\cal S}(\z(\g^{\xi}))$; il suffit en effet de l'{\'e}valuer en une forme lin{\'e}aire
qui ne s'annule pas en $\xi^{d-1}$. Par suite, la matrice
$\mathfrak{D}$ est non singuli{\`e}re.\\
    
On suppose maintenant $\g=\mathfrak{sp}_{2n}$, $\mathfrak{so}_{2n+1}$
ou $\mathfrak{so}_{2n}$. Posons $r=\left[ \frac{d}{2}
\right]$. On proc{\`e}de de la m{\^e}me fa{\c c}on; on pose $e_{i}=\xi^{2i-1}$, pour
$i=1,\ldots,r$. D'apr{\`e}s le lemme \ref{etude_z'}, les {\'e}l{\'e}ments $e_{1},\ldots,e_{r}$ forment une base
de $\z'$. Notons $\mathfrak{D}=(\mathfrak{D}_{ij})_{1 \leq i,j \leq
  r}$ la matrice carr{\'e}e d'ordre $r$ donn{\'e}e par\,:
$$\mathfrak{D}_{ij}=[[\eta,e_{j}],e_{i}],$$
pour $i$ et $j$ dans $\{1,\ldots,r\}$. D'apr{\`e}s la relation
(\ref{eta}), on a 
\begin{equation*}
\mathfrak{D}_{ij}=
\begin{cases}
0 & {\rm \ si \ } i+j > r+1,\\
-2(2i-1)(2j-1) \xi^{2i+2j-3} & {\rm \ sinon.} 
\end{cases}
\end{equation*}
En particulier, on connait $\mathfrak{D}$ sur l'anti-diagonale. On
a $\mathfrak{D}_{ij}=-2(2i-1)(2j-1)\xi^{2r-1}$ si $i+j=r+1$. Soit $\phi$ une forme lin{\'e}aire sur
$\z(\g^{\xi})$. De ce qui pr{\'e}c{\`e}de, il r{\'e}sulte que le d{\'e}terminant de la
matrice $\mathfrak{D}$ {\'e}valu{\'e} en $\phi$ vaut, au signe pr{\`e}s\,:
$$2^r \times ((2s-1) \times (2s-3) \times \cdots \times 1)^2 
\times (\langle \phi,\xi^{2r-1} \rangle )^{r} \cdot$$
C'est donc un {\'e}l{\'e}ment non nul de
${\cal S}(\z(\g^{\xi}))$; il suffit en effet de l'{\'e}valuer en une forme lin{\'e}aire
qui ne s'annule pas en $\xi^{2r-1}$. Par suite, la matrice
$\mathfrak{D}$ est non singuli{\`e}re.\\

Dans tous les cas, on obtient que la matrice
$\mathfrak{D}$ est non singuli{\`e}re. Notons que l'on n'a pas encore utilis{\'e}
l'hypoth{\`e}se $\z'=\z(\g^{\xi})$. On pose maintenant $m=\dim \z'=\dim
\z(\g^{\xi})$, d'o{\`u} $m=d-1$ dans le cas $\mathfrak{sl}_n$, et $m=r$
sinon.  On compl{\`e}te la base $\{e_{1},\ldots,e_{m}\}$ de $\z(\g^{\xi})$ en
une base $\{e_1,\ldots,e_n\}$ de $\g^{\xi}$. Alors, dans les bases
$\{e_1,\ldots,e_n,[\eta,e_1],\ldots,[\eta,e_m]\}$ de $\n(\g^{\xi})$ et
$\{e_1,\ldots,e_m\}$ de $\z(\g^{\xi})$, on a\,:
$${\cal
  K}(\n(\g^{\xi}),\z(\g^{\xi}))=  
\left[
\begin{array}{ccc}
0 & 0 & \mathfrak{D} 
\end{array}
\right] \cdot$$ 
Puisque $\mathfrak{D}$ est non singuli{\`e}re, la rang de la matrice
pr{\'e}c{\'e}dente est {\'e}gal {\`a} $\dim \z(\g^{\xi})$ et l'indice ${\rm ind \
}(\n(\g^{\xi}),\z(\g^{\xi}))$ est nul. L'autre assertion r{\'e}sulte des
remarques qui suivent la proposition \ref{rang}.
\qed
\\
\rem La preuve pr{\'e}c{\'e}dente permet en outre de voir que les formes
lin{\'e}aires r{\'e}guli{\`e}res de $\z(\g^{\xi})^*$ pour l'action naturelle de
$\n(\g^{\xi})$ dans $\z(\g^{\xi})$ sont exactement celles qui ne s'annulent pas sur le
vecteur de plus haut poids pour l'action de ${\rm ad} \rho$.\\ 

On d{\'e}duit de cette proposition le
th{\'e}or{\`e}me \ref{intro2} pour tous les {\'e}l{\'e}ments nilpotents de  
$\mathfrak{sl}_n,\mathfrak{sp}_{2n}$ et $\mathfrak{so}_{2n+1}$ et pour
les {\'e}l{\'e}ments nilpotents de $\mathfrak{so}_{2n}$ dont la partition
associ{\'e}e poss{\`e}de au moins trois parties. En effet dans chacune de ces
situations les sous-espaces $\z'$ et $\z(\g^{\xi})$ co{\"\i}ncident d'apr{\`e}s
le lemme \ref{etude_z'}. Remarquons que jusqu'ici, on a rien obtenu de
nouveau par rapport {\`a} ce qui a d{\'e}j{\`a} {\'e}t{\'e} fait par D. Panyushev dans
\cite{Panyushev}, si ce n'est que l'on a utilis{\'e}
la relation ${\rm ind \ } \g^{\xi}={\rm rg \ } \g$. D'ailleurs, la
preuve pr{\'e}c{\'e}dente est, {\`a} quelques modifications pr{\`e}s, celle de
D. Panyushev; elle m'a {\'e}t{\'e} sugg{\'e}r{\'e}e par M. Ra{\"\i}s.\\

On consacre la fin de cette partie {\`a} la preuve du th{\'e}or{\`e}me
\ref{intro2} pour $\g=\mathfrak{so}_{2n}$.
\begin{prop}\label{so2n}. On suppose que
  $\g=\mathfrak{so}_{2n}$. Alors on a la relation\,:
$${\rm ind \ }(\n(\g^{\xi}),\g^{\xi})={\rm rg \ } \g -\dim
\z(\g^{\xi}) \cdot$$
 \end{prop}
\dem On raisonne par r{\'e}currence sur le rang $n$ de $\g$. Pour $n=1$, $\g$ est
isomorphe {\`a} $\mathfrak{sl}_{2}$ et le r{\'e}sultat est connu. Soit $n \geq
2$. On suppose la proposition d{\'e}montr{\'e}e pour $\g=\mathfrak{so}_{2k}$, avec $k$ dans
$\{1,\ldots,n-1\}$. Il r{\'e}sulte de l'hypoth{\`e}se de r{\'e}currence et de ce
qui pr{\'e}c{\`e}de que la relation 
${\rm ind \ }(\n(\mathfrak{t}^{\xi}),\mathfrak{t}^{\xi})={\rm rg \ }
\mathfrak{t} -\dim \z(\mathfrak{t}^{\xi})$ est satisfaite pour
$\mathfrak{t}$ et $\xi$ satisfaisant les conditions 1), 2) et 3) qui pr{\'e}c{\`e}dent la proposition
\ref{distingue} et pour $\mathfrak{t}$ strictement contenue dans $\g$,
car alors $\mathfrak{t}$ est produit direct d'alg{\`e}bres de Lie simples
classiques de rang strictement inf{\'e}rieur {\`a} celui de $\g$. D'apr{\`e}s la proposition \ref{distingue}, il suffit donc de
prouver le r{\'e}sultat pour les {\'e}l{\'e}ments nilpotents distingu{\'e}s de $\g$. On peut
supposer en outre que $\xi$ n'est pas un {\'e}l{\'e}ment nilpotent r{\'e}gulier de
$\g$, d'apr{\`e}s les remarques qui pr{\'e}c{\`e}dent la proposition \ref{distingue}. On suppose d{\'e}sormais que $\xi$ est un {\'e}l{\'e}ment nilpotent
distingu{\'e} non r{\'e}gulier de $\g$.\\

D'apr{\`e}s \cite{Collinwood}, Th{\'e}or{\`e}me 8.2.14, les orbites nilpotentes distingu{\'e}es de
$\mathfrak{so}_{2n}$ sous l'action du groupe adjoint sont les orbites nilpotentes dont la partition associ{\'e}e n'a que des parties
impaires n'apparaissant qu'une seule fois. D'autre part, toujours
d'apr{\`e}s \cite{Collinwood} (Proposition 5.4.1, (iv)), l'orbite nilpotente r{\'e}guli{\`e}re de
$\mathfrak{so}_{2n}$ est donn{\'e}e par la partition $[2n-1,1]$. Or,
d'apr{\`e}s la proposition \ref{indice_nul} et le lemme \ref{etude_z'}, il
reste seulement {\`a} traiter le cas o{\`u} la partition associ{\'e}e {\`a} $\xi$ ne poss{\`e}de
que deux parties. En r{\'e}sum{\'e}, il suffit de prouver la proposition pour les orbites nilpotentes associ{\'e}es aux partitions de la forme
$[2s+1,2t+1]$, avec $1 \leq t < s \leq n-2$.\\

On fixe un {\'e}l{\'e}ment
nilpotent $\xi$ dans $\mathfrak{so}_{2n}$ associ{\'e} {\`a} une partition de
$2n$ de la forme $[2s+1,2t+1]$, avec $1 \leq t
  < s \leq n-2$. Le polyn{\^o}me minimal de $\xi$ est de degr{\'e} $2s+1$. La
  dimension de $\z'$ est donc {\'e}gale {\`a} $\left[ \frac{2s+1}{2} \right]=s$ et la
  dimension de $\z(\g^{\xi})$ est $s+1$, d'apr{\`e}s le lemme
  \ref{etude_z'}. Soit $b$ une forme bilin{\'e}aire sym{\'e}trique non d{\'e}g{\'e}n{\'e}r{\'e}e sur $\C$. On
r{\'e}alise $\g$ comme l'ensemble des endomorphismes anti-sym{\'e}triques de
$V=\C^{2n}$ relativement {\`a} la forme $b$. Comme $2s+1$ est diff{\'e}rent
de $2t+1$, il r{\'e}sulte de la preuve du lemme
5.1.17 de \cite{Collinwood} qu'il existe deux sous-espaces orthogonaux
$V_{1}$ et $V_{2}$, stables par $\xi$, et de dimensions respectives $2s+1$ et $2t+1$. La restriction $\xi_1$ de
$\xi$ {\`a} $V_1$ est un {\'e}l{\'e}ment nilpotent r{\'e}gulier dans
$\mathfrak{so}_{2s+1}$. De m{\^e}me, la restriction $\xi_2$ de
$\xi$ {\`a} $V_2$ est un {\'e}l{\'e}ment nilpotent r{\'e}gulier dans
$\mathfrak{so}_{2t+1}$. Il existe alors deux $\mathfrak{sl}_{2}$-triplets
$\{\xi_1,\rho_1,\eta_1\}$ et $\{\xi_2,\rho_2,\eta_2\}$ dans
$\mathfrak{so}_{2s+1}$ et $\mathfrak{so}_{2t+1}$ contenant $\xi_1$ et
$\xi_2$ respectivement. On note encore, de mani{\`e}re
abusive, $\xi_1$ l'{\'e}l{\'e}ment de $\g$ qui co{\"\i}ncide avec $\xi_{1}$ sur
$V_1$ et qui est nul sur $V_2$. On adopte la m{\^e}me convention pour $\rho_1$, $\eta_1$, $\xi_2$,
$\rho_2$ et $\eta_2$. On
pose $\rho=\rho_1 \oplus \rho_2$ et $\eta=\eta_1 \oplus \eta_2$. Le
triplet $\{ \xi,\rho,\eta \}$ forme ainsi un
$\mathfrak{sl}_{2}$-triplet dans $\g$. Soit ${\cal B}_1$ et ${\cal B}_2$ des bases
orthonorm{\'e}es de $V_1$ et $V_2$ et soit ${\cal
  B}$ la base orthonorm{\'e}e de $V$ obtenue en concat{\'e}nant les bases
${\cal B}_1$ et ${\cal B}_2$. On note $X$ (respectivement $X_1$,
$H_1$, $Y_1$ et $X_2$, $H_2$, $Y_2$) la matrice de $\xi$
(respectivement $\xi_1$, $\rho_1$, $\eta_1$ et $\xi_2$, $\rho_2$, $\eta_2$) dans la base
${\cal B}$ (respectivement ${\cal B}_1$ et ${\cal B}_2$). Ainsi,
$$X=\left[
\begin{array}{cc}
X_1 & 0 \\
0 & X_2
\end{array}
\right], \hspace{1cm}
H=\left[
\begin{array}{cc}
H_1 & 0 \\
0 & H_2
\end{array}
\right], \hspace{1cm}
Y=\left[
\begin{array}{cc}
Y_1 & 0 \\
0 & Y_2
\end{array}
\right].$$

Soit $A$ une matrice non nulle de taille $(2s+1) \times (2t+1)$ v{\'e}rifiant les
relations $X_1A=0$ et $AX_2=0$. On v{\'e}rifie ais{\'e}ment qu'une telle
matrice existe. C'est nec{\'e}ssairment une matrice de rang
1. En effet, commme $X_1$ et $X_2$ sont des matrices nilpotentes
r{\'e}guli{\`e}res, leur noyau est de dimension $1$. On note $w$ l'{\'e}l{\'e}ment de $\g$
dont la matrice dans la base ${\cal B}$ est 
$$W=\left[
\begin{array}{cc}
0 & A \\
-A^t & 0
\end{array}
\right].$$     
On a le lemme suivant\,:
\begin{lemme}\label{centre_w}
L'{\'e}l{\'e}ment $w$ appartient au centre $\z(\g^{\xi})$ et c'est un
vecteur propre pour ${\rm ad} \rho$.  
\end{lemme}
\dem Prouvons d'abord la premi{\`e}re assertion. On montre sans difficult{\'e} que les {\'e}l{\'e}ments de $\g$ qui commutent
avec $\xi$ sont les endomorphismes de $\g$ dont la matrice dans la
base ${\cal B}$ est de la forme\,:
$$\left[
\begin{array}{cc}
X_{1}^{2k+1} & R \\
-R^t & X_{2}^{2l+1}
\end{array}
\right],$$
avec $k$ et $l$ deux entiers et $R$ une matrice de taille
$(2s+1)\times(2t+1)$ v{\'e}rifiant la relation $X_1 R-RX_2=0$. Il est
imm{\'e}diat que la matrice $W$ commute avec les puissances de $X_1$ et 
$X_2$. Il reste {\`a} montrer que $W$ commute avec les
matrices de la forme\,:
$$\left[
\begin{array}{cc}
0 & R \\
-R^t & 0
\end{array}
\right],$$
o{\`u} $R$ est une matrice de taille
$(2s+1)\times(2t+1)$ v{\'e}rifiant la relation $X_1 R-RX_2=0$. Soit $Z$ 
une telle matrice. Prouvons que le crochet $[W,Z]=WZ-ZW$ est nul. Il
s'agit de montrer que les
deux relations suivantes sont satisfaites\,:
$$\left\lbrace
\begin{array}{l}
AR^t=RA^t \\
A^t R=R^t A
\end{array}
\right.$$
Montrons d'abord la premi{\`e}re {\'e}galit{\'e}. Ceci revient {\`a} montrer que la
matrice carr{\'e}e $AR^t$ d'ordre $2s+1$ est sym{\'e}trique. Si cette
matrice est nulle, le r{\'e}sultat est clair. On suppose que la matrice $AR^t$ n'est
pas nulle, alors il en est de m{\^e}me de sa transpos{\'e}e $RA^t$. De la
relation $AX_1=0$, on tire la relation\,: $X_1(AR^t)=0$ et des
relations $X_1 R=RX_2$ et $X_2 A^t=0$, on tire la relation\,: $X_1
(RA^t)=RX_2 A^t=0$. Les images des deux matrices $AR^t$ et $RA^t$ sont
donc contenues dans le noyau de la matrice $X_1$, qui est de dimension
1. Par suite, les deux matrices $AR^t$ et $RA^t$ sont de rang $1$ et ont
la m{\^e}me image. On note $u$ l'endomorphisme de $\mathfrak{gl}_{2s+1}$ dont la matrice dans la base ${\cal B}_1$
est $AR^t$. Comme $RA^t$ est la transpos{\'e}e de $AR^t$, la matrice
$RA^t$ repr{\'e}sente l'endomorphisme adjoint $u^*$ de $u$ dans la base ${\cal B}_1$. De ce qui pr{\'e}c{\`e}de, il resulte
que les endomorphismes $u$ et $u^*$ s'{\'e}crivent sous la forme 
$u=\phi v$ et $u^*=\psi v$, avec $v$ un vecteur non nul de
$V_1$ et $\phi$ et $\psi$ deux formes lin{\'e}aires non nulles sur
$V_1$. Pour tout $x$ dans $V_1$, on a l'{\'e}galit{\'e}\,: 
$b(u(x),x)=b(\phi(x)v,x)=\phi(x)b(v,x)$. Or, par d{\'e}finition de l'adjoint, on a la relation\,:
$b(u(x),x)=b(x,u^*(x))=\psi(x)b(v,x)$.
Ainsi, pour tout $x$ non orthogonal {\`a} $v$, on a
$\phi(x)=\psi(x)$. Comme l'ensemble des vecteurs $x$ de $V_1$ v{\'e}rifiant la relation $b(x,v)
\not= 0$ est un ouvert non vide de $V_1$, les formes lin{\'e}aires $\phi$
et $\psi$ sont {\'e}gales. On a obtenu l'{\'e}galit{\'e} souhait{\'e}e\,:
$AR^t=RA^t$. Un raisonnement similaire permet d'obtenir la deuxi{\`e}me relation.\\

On en d{\'e}duit la deuxi{\`e}me assertion. L'{\'e}l{\'e}ment $w$ appartient au centre $\z(\g^{\xi})$ de
$\g^{\xi}$ d'apr{\`e}s ce qui pr{\'e}c{\`e}de et il est clair que ce n'est pas un
polyn{\^o}me en $\xi$. Par suite, il r{\'e}sulte du lemme \ref{etude_z'} que les {\'e}l{\'e}ments $\xi,\xi^{3},\ldots,\xi^{2s-1},w$ forment une base de
$\z(\g^{\xi})$. Montrons que le crochet $[H,W]$ est colin{\'e}aire {\`a} $W$. On a\,:
$$[H,W]=
\left[
\begin{array}{cc}
0 & H_1 A -A H_2 \\
-H_2 A^t+A^tH_1
\end{array}
\right].$$
Le crochet $[H,W]$ est un {\'e}l{\'e}ment du centre et ses
coordonn{\'e}es en les puissances de $X$ sont nulles. Il est donc 
colin{\'e}aire {\`a} $W$ et le lemme s'ensuit. 
\qed

On note $\lambda$ la valeur propre
associ{\'e}e {\`a} $w$.
\begin{lemme}\label{crochet} \begin{itemize}
\item[(i)] On a\,:\begin{equation*}
[[\eta,w],\xi^{2i-1}]=
\begin{cases}
-\lambda w & {\rm \ si \ } i=1,\\
0 & {\rm \ si \ } i=2,\ldots,s. 
\end{cases}
\end{equation*}
\item[(ii)] L'{\'e}l{\'e}ment $[[\eta,w],\xi_2]$ est non nul dans
  $\g^{\xi}$. On note $x$ cet {\'e}l{\'e}ment. 
\item[(iii)] Le crochet $[[\eta,\xi^{2s-1}],\xi_2]$ est nul. 
\end{itemize}
\end{lemme}
\dem i) Si $i=1$, on a
$[[\eta,w],\xi^{2i-1}]=[[\eta,\xi],w]=-[\rho,w]=-\lambda w$. On
suppose $i>1$. En utilisant la relation (\ref{rho}), on montre sans
peine que le crochet $[[\eta,w],\xi^{2i-1}]$ est donn{\'e} par la matrice\,:
$$\left[
\begin{array}{cc}
0 & U\\
-U^t & 0
\end{array}
\right],$$
o{\`u} $U$ est la matrice $-\sum\limits_{\alpha+\beta=2i-2} X_{1}^{\alpha} H_1 X_{1}^{\beta} A
+ A \sum\limits_{\alpha+\beta=2i-2} X_{2}^{\alpha} H_2 X_{2}^{\beta}$. 
Dans la somme $\sum\limits_{\alpha+\beta=2i-2} X_{1}^{\alpha} H_1
X_{1}^{\beta} A$, les termes pour lesquel $\beta$ est strictement
positif sont nuls,
puisqu'on a la relation $X_1 A=0$. Il reste le terme $X_{1}^{2i-2}H_1
A$. Comme $w$ est un vecteur propre de ${\rm ad} \rho$ relativement {\`a}
la valeur propre $\lambda$, on a la relation $H_1 A=A
H_2 + \lambda A$. D'o{\`u}\,:
$$X_{1}^{2i-2} H_1 A=X_{1}^{2i-2} A H_2 + \lambda X_{1}^{2i-2} A=0,$$
car $X_1 A=0$ et car $i>1$. Par suite, on a\,: 
$\sum\limits_{\alpha+\beta=2i-2} X_{1}^{\alpha} H_1 X_{1}^{\beta}
A=0$. De m{\^e}me, la relation $AX_2=0$ entraine la
relation 
$A \sum\limits_{\alpha+\beta=2i-2} X_{2}^{\alpha} H_2
X_{2}^{\beta}=0$. Par suite la matrice $U$ est nulle et 
le crochet $[[\eta,w],\xi^{2i-1}]$ est nul pour $i>1$, d'o{\`u} (i).\\

ii) On a la relation
$x=[[\eta,w],\xi_2]=[[\eta,\xi_2],w]=-[\rho_2,w]$. On suppose par
l'absurde que $[\rho_2,w]$ est nul. Il s'agit d'aboutir {\`a} une
contradiction. De la relation $AX_2=0$, on tire la relation
$[\xi_2,w]=0$. Comme le triplet $\{\xi_2,\rho_2,\eta_2\}$ forme un
$\mathfrak{sl}_{2}$-triplet dans $\g$, on d{\'e}duit des {\'e}galit{\'e}s
$[\rho_2,w]=0$ et $[\xi_2,w]=0$ la relation $[\eta_2,w]=0$, d'o{\`u} il
vient $A Y_2=0$. Ainsi, on a\,: $X_2 A^t=0$ et $Y_2 A^t=0$, donc l'image
de la matrice $A^t$ est contenue dans l'intersection des noyaux de $X_2$ et
$Y_2$. Comme $X_2$ est une matrice nilpotente r{\'e}guliere, cette
intersection est nulle d'apr{\`e}s la th{\'e}orie des repr{\'e}sentations de
$\mathfrak{sl}_2$. D'o{\`u} la contradiction car $A$ n'est pas nulle.\\

iii) Comme $\xi_{2}^{2s-1}$ est nul, le crochet $[\eta,\xi^{2s-1}]$ est donn{\'e} par la matrice
$$\left[
\begin{array}{cc}
-\sum\limits_{\alpha+\beta=2s-2} X_{1}^{\alpha} H_{1} X_{1}^{\beta} & 0 \\
0 & 0
\end{array}
\right].$$
Par suite le crochet $[[\eta,\xi^{2s-1}],\xi_{2}]$ est nul. 
\qed

On pose $e_{i}=\xi^{2i-1}$, pour $i=1,\ldots,s$, $e_{s+1}=w$ et
$e_{s+2}=\xi_2$. Comme $\xi_{2}$ n'appartient pas au centre
$\z(\g^{\xi})$, la famille $\{e_{1},\ldots,e_{s+2}\}$ est libre dans
$\g^{\xi}$. On la compl{\`e}te en une base $\{e_{1},\ldots,e_{n}\}$ de
$\g^{\xi}$ et on note 
$\mathfrak{M}=(\mathfrak{M}_{ij})_{1 \leq i \leq n \atop 1 \leq j
  \leq s+1}$ la matrice de taille $n \times (s+1)$, {\`a} coefficients
dans ${\cal S}(\g^{\xi})$, donn{\'e}e par\,:
$$\mathfrak{M}_{ij}=[[\eta,e_{j}],e_{i}],$$
pour $i=1,\ldots,n$ et $j=1,\ldots,s+1$. On consid{\`e}re la sous-matrice
carr{\'e}e $\mathfrak{M}'$ de $\mathfrak{M}$ d'ordre $s+1$ dont les
lignes correspondent aux {\'e}l{\'e}ments $e_{1},\ldots,e_{s}$ et $e_{s+2}$ et
dont les colonnes correspondent aux {\'e}l{\'e}ments $e_{1},\ldots,e_{s+1}$ de
$\z(\g^{\xi})$.

Gr{\^a}ce au lemme \ref{crochet} et aux calculs effectu{\'e}s dans la preuve
de la proposition \ref{indice_nul}, on obtient la structure de la
matrice $\mathfrak{M}'$\,:
\newsavebox{\matriceM}
\savebox{\matriceM}[4cm]{
$\mathfrak{M}'=
\left[
\begin{array}{ccccc}
\ast   & \cdots & \ast & \mu_{1} \xi^{2s-1} & -\lambda w \\
\ast   & \cdots & \mu_{2}\xi^{2s-1} & 0 & 0 \\
\vdots &        &                                     & \vdots &
\vdots \\
\mu_{s} \xi^{2s-1} & 0 & \cdots & 0 & 0 \\
\ast & \cdots & \ast & 0 & x
\end{array}
\right],$
}
\begin{center}
\psset{unit=1cm}
\begin{pspicture}(-3,-1.5)(3,1.5)
\rput(0,0){\usebox{\matriceM}}
\psline[linestyle=dotted,linewidth=0.43mm](-1.2,-0.15)(-0.75,0.1)
\end{pspicture}
\end{center}
o{\`u} $\mu_{i}=-2 \times (2s-(2i-1))\times (2i-1)$, pour $i=1,\ldots,s$. 
Soit $\phi$ une forme lin{\'e}aire sur $\g^{\xi}$. D'apr{\`e}s la structure de la
matrice $\mathfrak{M}'$, l'{\'e}valuation en $\phi$ du d{\'e}terminant de
$\mathfrak{M}'$ vaut, au signe pr{\`e}s\,:
$$ 2^{s} \times ((2s-1) \times (2s-3) \times \cdots \times 1)^2 \times 
\langle \phi,x \rangle \times (\langle \phi,\xi^{2s-1}
\rangle)^{s} \cdot$$
Si $\phi$ ne s'annule ni en $x$, ni en $\xi^{2s-1}$, ce d{\'e}terminant
est non nul. Donc, la matrice
$\mathfrak{M}'$ est non singuli{\`e}re. Avec les
notations de la proposition \ref{rang}, la matrice $\mathfrak{M}$ 
est la matrice 
$\left[
\begin{array}{c}
\mathfrak{D} \\
\mathfrak{E}
\end{array}
\right]$. Cette derni{\`e}re est donc de rang maximal {\'e}gal {\`a} $s+1=\dim
\z(\g^{\xi})$, car la sous-matrice carr{\'e}e $\mathfrak{M}'$ d'ordre 
$s+1$ est non
singuli{\`e}re. La proposition r{\'e}sulte alors de la proposition
\ref{rang}.
\qed
\\
\rem Les calculs de la preuve pr{\'e}c{\'e}dente permettent de donner la
structure de la matrice $\mathfrak{D}$ dans le cas o{\`u} $\xi$ est
associ{\'e} {\`a} une partition de la forme $[2s+1,2t+1]$, avec $1 \leq t < 
s \leq n-2$. On a\,:\\
\newsavebox{\matrice}
\savebox{\matrice}[4cm]{
$\mathfrak{D} \ =
\left[
\begin{array}{ccccc}
\ast   & \cdots & \ast & \mu_{1} \xi^{2s-1} & -\lambda w \\
\ast   & \cdots & \mu_{2}\xi^{2s-1} & 0 & 0 \\
\vdots &        &                                     & \vdots &
\vdots \\
\mu_{s} \xi^{2s-1} & 0 & \cdots & 0 & 0 \\
-\lambda w & 0 & \cdots & 0 & 0
\end{array}
\right].$
}
\begin{center}
\psset{unit=1cm}
\begin{pspicture}(-3,-1.5)(3,1.5)
\rput(0,0){\usebox{\matrice}}
\psline[linestyle=dotted,linewidth=0.43mm](-1.2,-0.15)(-0.75,0.1)
\end{pspicture}
\end{center}
La matrice $\mathfrak{D}$ est donc singuli{\`e}re dans ce cas et la
relation ${\rm ind}(\n(\g^{\xi}),\z(\g^{\xi}))=0$ n'a pas lieu.\\  

\section{D{\'e}monstration du th{\'e}or{\`e}me \ref{intro2} sous l'hypoth{\`e}se que
  $\xi$ v{\'e}rifie la propri{\'e}t{\'e} $(P)$.}

On suppose que $\g$  est une alg{\`e}bre de Lie complexe simple et 
que $\xi$ est un {\'e}l{\'e}ment nilpotent non
r{\'e}gulier de $\g$. On fixe
un $\mathfrak{sl}_{2}$-triplet $\{\xi,
\rho,\eta\}$ dans $\g$ contenant $\xi$. On note $W$ le sous-espace
propre de la restriction de ${\rm ad} \rho$ {\`a} $\z(\g^{\xi})$ associ{\'e} {\`a} la plus grande
valeur propre et on suppose que $\xi$ v{\'e}rifie la propri{\'e}t{\'e} $(P)$ de la
d{\'e}finition \ref{prop_P}. Le but de cette partie est d'obtenir le
corollaire \ref{rang_z'} qui servira au cas exceptionnel dans la
partie suivante.

On va utiliser un certain nombre de notations et de
r{\'e}sultats de \cite{Charbonnel}. On note $\pi_{\xi}$ l'application $(g,x) \mapsto
g(x)$ de $G \times (\xi + \g^{\eta})$ dans $\g$. D'apr{\`e}s
\cite{Charbonnel}, corollaire 5.7, il
existe un voisinage ouvert $W$ de $\xi$ dans $\xi + \g^{\eta}$ tel que
la restriction de $\pi_{\xi}$ {\`a} $G \times W$ soit un morphisme lisse
de $G \times W$ sur un ouvert $G$-invariant de $\g$ qui contient
$\xi$. On note $X$ l'{\'e}clatement en $\xi$ de $W$ et $\sigma$ le morphisme
d'{\'e}clatement. Ce qui pr{\'e}c{\`e}de montre
que l'ouvert non vide $\pi_{\xi}(G \times W)$ rencontre l'ouvert dense de $\g$ des {\'e}l{\'e}ments r{\'e}guliers de
$\g$. Par suite, le sous-ensemble de $X$ des points $x$
de $X$ pour lesquels $\sigma(x)$ est un {\'e}l{\'e}ment r{\'e}gulier de $\g$ est
un ouvert non vide de $X$. On note $X_{r}$ cet ouvert. En particulier, $\sigma^{-1}(\{\xi\})$ est une hypersurface de
$X$, contenue dans $X \setminus X_{r}$, car l'{\'e}l{\'e}ment $\xi$ n'est pas un {\'e}l{\'e}ment
r{\'e}gulier de $\g$. On note enfin $X_{*}$ le plus grand ouvert de $X$ auquel l'application $x \mapsto \g^{\sigma(x)}$ de $X_{r}$ dans
${\rm Gr}_{{\rm rg \; } \g}(\g)$ a un
prolongement r{\'e}gulier not{\'e} $\alpha$. D'apr{\`e}s \cite{Charbonnel}, Lemme 2.3, il existe un ouvert affine $Y$ de $X$ qui rencontre
$\sigma^{-1}(\{\xi\})$ et qui satisfait les conditions suivantes\,:
\begin{itemize}
\item[ 1) ] $Y$ est contenu dans $X_{*}$,
\item[ 2) ] $Y \setminus X_{r}$ est une hypersurface lisse,
  irr{\'e}ductible, contenue dans $\sigma^{-1}(\{\xi\})$ et dont l'id{\'e}al de
  d{\'e}finition dans $\C[Y]$ est engendr{\'e} par un {\'e}l{\'e}ment $q$,
\item[ 3) ] il existe un sous-espace $\m$ de $\g$ qui est un
  suppl{\'e}mentaire de $\alpha(x)$ dans $\g$ pour tout $x$ dans $Y$
  et qui contient un suppl{\'e}mentaire $\p$ de $\g^{\xi}$ dans $\g$.\\  
\end{itemize}

Par d{\'e}finition de $\alpha$, on a l'inclusion $\alpha(x)
\subset \g^{\sigma(x)}$, pour tout $x$ dans $Y$. En particulier, on a
l'inclusion $\alpha(x)  \subset \g^{\xi}$, pour tout $x$ dans $ Y \setminus
X_{r}$, car $Y \setminus X_{r}$ est contenu dans
$\sigma^{-1}(\{\xi\})$. On dispose alors du lemme suivant\,:
\begin{lemme}\label{supplementaire}
Soit $x_{0}$ un point de $Y \setminus X_{r}$ et soit $\{v_{1},\ldots,v_{s}\}$ une base d'un suppl{\'e}mentaire $\n$ de
$\alpha(x_{0})$ dans $\g^{\xi}$. Alors l'{\'e}l{\'e}ment $\det ([v_{i},v_{j}])_{1
  \leq i,j \leq s}$ est un {\'e}l{\'e}ment non  nul de ${\cal
  S}(\g^{\xi})$.
\end{lemme}
\dem D'apr{\`e}s la proposition \ref{centralisateur}, $[\eta,
\g]$ est un suppl{\'e}mentaire de $\g^{\xi}$ dans $\g$. Par suite, on a les 
d{\'e}compositions suivantes\,:
$$\alpha(x_{0}) \oplus \n=\g^{\xi} {\rm \ \  et \ \ } \alpha(x_{0}) \oplus \n
\oplus [\eta, \g]=\g \cdot$$
On en d{\'e}duit qu'il existe un ouvert non vide $Y'$ de $Y$ contenant
$x_{0}$ tel que, pour tout $x'$ dans $Y'$, on ait la d{\'e}composition\,:
$$ \alpha(x') \oplus \n \oplus [\eta, \g]=\g \cdot$$
Puisque $\alpha(x')$ est contenu dans $\g^{\xi}$ pour tout $x'$ dans
$Y' \setminus X_{r}$, on obtient de plus
$$ \alpha(x') \oplus \n=\g^{\xi},$$
pour tout $x'$ dans $Y' \setminus X_{r}$.\\

Soit $f_{1},\ldots,f_{t}$ une base de $[\eta,\g]$. Il suffit de
montrer que l'{\'e}l{\'e}ment 
$$\Delta= \det 
\left[
\begin{array}{cc}
\begin{array}{ccc} 
[v_{1},v_{1}]  & \cdots & [v_{1},v_{s}]\\
\vdots & & \vdots\\ 
\left[v_{s},v_{1}\right] & \cdots & [v_{s},v_{s}] 
\end{array} & 0  \\

\begin{array}{ccc}
[f_{1},v_{1} ] & \cdots & [ f_{1},v_{s} ]  \\
\vdots        &        & \vdots\\
\left[f_{t},v_{1}\right]  & \cdots & [ f_{t},v_{s} ] 
\end{array} & 
\begin{array}{ccc}
[ f_{1},f_{1} ]  & \cdots & [ f_{1},f_{t} ] \\
\vdots        &        & \vdots\\
\left[f_{t},f_{1}\right]  & \cdots & [ f_{t},f_{t} ] 
\end{array}
\end{array}
\right]
$$
de ${\cal S}(\g)$ est non nul pour obtenir le lemme.\\

Soit $i$ dans $\{1,\ldots,s\}$. L'application qui {\`a} $x$ dans $Y$
associe $[\sigma(x),v_{i}]$ est nulle en tout point de $Y \setminus X_{r}$ car $Y \setminus
X_{r}$ est inclus dans $\sigma^{-1}(\{\xi\})$ et $\g^{\xi}$ contient
$v_{i}$. Puisque $q$ engendre l'id{\'e}al de d{\'e}finition de $Y \setminus
X_{r}$, on en d{\'e}duit qu'il existe une application r{\'e}guli{\`e}re  $\mu_{i}$
de $Y$ dans $\g$ non identiquement nulle sur $Y \setminus X_{r}$ et un
entier $m_{i} \geq 1$ qui satisfont l'{\'e}galit{\'e}
$[\sigma(x),v_{i}]=q(x)^{m_{i}} \mu_{i}(x)$, pour tout $x$ de $Y$.

Montrons que l'entier $m_{i}$ est {\'e}gal {\`a} $1$ pour tout $i$ dans
$\{1,\ldots,s\}$. On suppose par l'absurde que $m_{i} > 1$ pour un
certain $i$ dans $\{1,\ldots,s\}$. Il s'agit d'aboutir {\`a} une contradiction.
 
Soit $T$ une ind{\'e}termin{\'e}e, $\tau $ et  $\tau _{0}$ les images
respectives de $T$ par les applications canoniques de l'anneau de
polyn{\^o}mes  ${\mathbb C}[T]$ dans les quotients de ${\mathbb C}[T]$ par les
id{\'e}aux $T^{m_{i}+1}{\mathbb C}[T]$ et $T^{2}{\mathbb C}[T]$. Pour $\nu
$ un ${\mathbb C}[\tau ]$-point de $X$, on note $\gamma _{\nu }$ l'{\'e}valuation en
$\nu $ de l'anneau $\an Xx$ o{\`u} $x$ est l'image de $\nu $ par la
projection canonique de l'ensemble des ${\mathbb C}[\tau ]$-points de $X$
sur l'ensemble des ${\mathbb C}$-points de $X$. Alors $\gamma _{\nu }$ est
un morphisme de l'anneau $\an Xx$ dans l'anneau ${\mathbb C}[\tau ]$. Soit
$\gamma '_{\nu }$ le morphisme
$$ \tk {{\mathbb C}}{\an Xx}{\mathfrak g}\rightarrow 
 \tk {{\mathbb C}}{{\mathbb C}[\tau ]}{\mathfrak g} \mbox{ , }
 \varphi \tens u \mapsto \gamma _{\nu }(\varphi )\tens u \mbox{ ,}$$
o{\`u} $\varphi $ est dans $\an Xx$ et o{\`u} $u$ est dans ${\mathfrak g}$. Puisque $\sigma$ est le morphisme d'{\'e}clatement de $W$ centr{\'e} en $\xi$,
les fibres de $\sigma$ sont irr{\'e}ductibles et l'espace tangent {\`a} $W$ en
un point $z$ de $W$ est la r{\'e}union des images des applications
lin{\'e}aires tangentes {\`a} $\sigma$ en les points de la fibre de $\sigma$
en $z$. En particulier, puisque $Y \setminus X_{r}$ est un ouvert de
la fibre $\sigma^{-1}(\xi)$ de $\sigma$ en $\xi$, la r{\'e}union des images des
applications lin{\'e}aires tangentes {\`a} $\sigma$ en les points de $Y
\setminus X_{r}$ contient un ouvert non vide de $\g^{\eta}$.

Soit $v$ un vecteur de $\g^{\eta}$ appartenant {\`a} cet ouvert. Soit
alors $x$ dans $Y \setminus X_{r}$ et $v'$ un
vecteur tangent {\`a} $Y$ en $x$ tels que l'image de $v'$ par l'application lin{\'e}aire tangente
$\sigma'(x)$ {\`a} $\sigma$ en $x$ soit {\'e}gale {\`a} $v$. Soit enfin $\nu$ un ${\mathbb C}[\tau ]$-point de $Y$ au dessus du ${\mathbb C}[\tau _{0}]$-point de $X$ d{\'e}fini par $v'$. On note $q'(x)$ la diff{\'e}rentielle de $q$ en $x$. Puisque
l'image de $q^{m_{i}} \mu_{i}$ par $\gamma' _{\nu }$ est {\'e}gale {\`a}
$\tau ^{m_{i}}q'(x)(v') \mu_{i}(x)$,  de l'{\'e}galit{\'e}:
$$ [\sigma (x),v_{i}] = q(x)^{m_{i}}\mu _{i}(x) \mbox{ ,}$$
on tire l'{\'e}galit{\'e}:
$$ \tau [v,v_{i}] + \cdots = \tau ^{m_{i}} q'(x)(v')\mu_{i}(x) \mbox{ ,}$$
car le terme de degr{\'e} $1$ en $\tau $ de $\gamma' _{\nu }(\sigma )$ est 
$\tau v$. Puisque $m_{i}>1$, on obtient que $[v,v_{i}]$ est
nul. Il en r{\'e}sulte que pour tout $v$ dans un ouvert non vide de $\g^{\eta}$, $[v,v_{i}]$ est
nul. Par suite, $[v,v_{i}]$ est nul pour tout $v$ dans ${\mathfrak
  g}^{\eta }$. On en d{\'e}duit que $v_{i}$ appartient au centre $\z(\g^{\eta})$ de $\g^{\eta}$. Ceci est
impossible car $\z(\g^{\eta})$ et $\g^{\xi}$ ont une intersection nulle. On a finalement obtenu la relation
$[\sigma(x),v_{i}]=q(x)\mu_{i}(x)$, pour tout $x$ dans $Y$.
\\

Pour $x$ dans $Y$, la valeur de $\Delta$ en la forme lin{\'e}aire 
$v \mapsto \langle \sigma(x), v \rangle$ est:
\begin{eqnarray*}
& & \det 
\left[
\begin{array}{cc}
\begin{array}{ccc} 
\langle \sigma(x), [v_{1},v_{1}] \rangle & \cdots & \langle \sigma(x),
[v_{1},v_{s}] \rangle\\
\vdots & & \vdots\\ 
\langle \sigma(x), \left[v_{s},v_{1}\right] \rangle & \cdots & \langle \sigma(x),  [v_{s},v_{s}] \rangle
\end{array} & 0  \\

\begin{array}{ccc}
\langle \sigma(x), [f_{1},v_{1} ] \rangle & \cdots & \langle \sigma(x), [ f_{1},v_{s} ] \rangle \\
\vdots        &        & \vdots\\
\langle \sigma(x), \left[f_{t},v_{1}\right] \rangle & \cdots & \langle \sigma(x), [ f_{t},v_{s} ] \rangle
\end{array} & 
\begin{array}{ccc}
\langle \sigma(x), [ f_{1},f_{1} ]  \rangle & \cdots & \langle \sigma(x), [ f_{1},f_{t} ] \rangle\\
\vdots        &        & \vdots\\
\langle \sigma(x), \left[f_{t},f_{1}\right] \rangle & \cdots & \langle
\sigma(x), [ f_{t},f_{t} ] \rangle 
\end{array}
\end{array}
\right]\\
& = & 
(-q(x))^{s} \det 
\left[
\begin{array}{cc}
\begin{array}{ccc} 
\langle \mu_{1}(x), v_{1} \rangle & \cdots & \langle \mu_{s}(x), v_{1} \rangle\\
\vdots & & \vdots\\ 
\langle \mu_{1}(x), v_{s} \rangle & \cdots & \langle \mu_{s}(x), v_{s} \rangle
\end{array} & 0  \\

\begin{array}{ccc}
\langle \mu_{1}(x), f_{1} \rangle & \cdots & \langle \mu_{s}(x), f_{1} \rangle \\
\vdots        &        & \vdots\\
\langle \mu_{1}(x), f_{t} \rangle & \cdots & \langle
\mu_{s}(x),f_{t} \rangle
\end{array} & 
\begin{array}{ccc}
\langle \sigma(x), [ f_{1},f_{1} ]  \rangle & \cdots & \langle \sigma(x), [ f_{1},f_{t} ] \rangle\\
\vdots        &        & \vdots\\
\langle \sigma(x), \left[f_{t},f_{1}\right] \rangle & \cdots & \langle
\sigma(x), [ f_{t},f_{t} ] \rangle 
\end{array}
\end{array}
\right] \cdot
\end{eqnarray*}

Soit $\zeta$ la fonction r{\'e}guli{\`e}re d{\'e}finie sur $Y$ qui {\`a} $x$ dans $Y$
associe $(-q(x))^{-s}\Delta(\langle
\sigma(x), \cdot \rangle)$. On suppose par l'absurde
que $\zeta$ est identiquement nulle sur $Y \setminus X_{r}$. Alors il existe des fonctions r{\'e}guli{\`e}res non
toutes identiquement nulles, $a_{1},\ldots,a_{s},b_{1},\ldots,b_{t}$
sur $Y \setminus X_{r}$, qui satisfont aux {\'e}galit{\'e}s\,:
$$\left\{ \begin{array}{ll}
\langle \sum\limits_{j=i}^{s}a_{j}(x)\mu_{j}(x),v_{i} \rangle
=0, \; \forall i=1,\ldots,s \\
\langle \sum\limits_{j=1}^{s}a_{j}(x)\mu_{j}(x) + 
\sum\limits_{k=1}^{t}b_{k}(x) [\sigma(x),f_{k}],f_{i}\rangle
=0, \; \forall i=1,\ldots,t \; ,
\end{array}
\right. $$
pour tout $x \in Y \setminus X_{r}$. En utilisant les inclusions 
$Y \setminus X_{r} \subset \sigma^{-1}(\{\xi\})$ et 
$\n \subset \g^{\xi}$, on trouve la relation 
$\langle [\sigma(x),f_{k}],w \rangle=0$,
pour tout $w$, $x$ et $k$ dans $\n$, $Y \setminus X_{r}$ et 
$\{1,\ldots,t\}$ respectivement, d'o{\`u} l'{\'e}galit{\'e}\,:
\begin{eqnarray}\label{relation}
& \langle \sum\limits_{j=1}^{s}a_{j}(x)\mu_{j}(x) + 
\sum\limits_{k=1}^{t}b_{k}(x) [\sigma(x),f_{k}], w \rangle
=0,
\end{eqnarray}
pour tout $w$ dans $\n \oplus [\eta,\g]$ et tout $x$ dans $Y \setminus X_{r}$.\\

On fixe une base $u_{1},\ldots,u_{n}$ de $\g^{\eta}$ et on note
$x_{1},\ldots,x_{n}$ la base duale. Pour $i=1,\ldots,n$, on d{\'e}signe
aussi par $x_{i}$ la forme affine sur $\xi + \g^{\eta}$ dont la valeur
au point $\xi +x$ de $\xi + \g^{\eta}$ est la valeur de $x_{i}$ en
$x$. Soit $x$ un point de $Y \setminus X_{r}$. Quitte {\`a} changer de base
$u_{1},\ldots,u_{n}$, on peut supposer que l'ensemble des fonctions 
$$y_{1}=x_{1}, \; y_{2}=\frac{x_{2}}{x_{1}},
\ldots,y_{n}=\frac{x_{n}}{x_{1}}$$
est un syst{\`e}me de coordonn{\'e}es de l'anneau local ${\cal O}_{X,x}$ de
$X$ en $x$. Puisque $x_{1}$ engendre l'id{\'e}al de d{\'e}finition de
$\sigma^{-1}(\{\xi\})$ dans ${\cal O}_{X,x}$, la fonction $q/x_{1}$
est un {\'e}l{\'e}ment inversible de ${\cal O}_{X,x}$. L'{\'e}galit{\'e} (\ref{relation})
revient {\`a} dire qu'il existe des polyn{\^o}mes 
$p_{1},\ldots,p_{s},q_{1},\ldots,q_{t}$ en $n-1$ ind{\'e}termin{\'e}es, non tous nuls, qui satisfont l'{\'e}galit{\'e}\,:
\begin{eqnarray*}
& \langle
\sum\limits_{j=1}^{s}p_{j}(y_{2},\ldots,y_{n})
[\xi+u_{1}+y_{2}u_{2}+\cdots+y_{n}u_{n},v_{j}] \\
& \hspace{3cm} +\sum\limits_{k=1}^{t}q_{k}(y_{2},\ldots,y_{n})
[\xi+u_{1}+y_{2}u_{2}+\cdots+y_{n}u_{n},f_{k}] , w \rangle=0,
\end{eqnarray*}
pour tout $w$ dans $\n \oplus [\eta,\g]$. L'{\'e}galit{\'e}
pr{\'e}c{\'e}dente reste valable pour tout $w$ qui centralise
$\xi+u_{1}+y_{2}u_{2}+\cdots+y_{n}u_{n}$. Or pour tout $v$ dans un
ouvert non vide de $\g^{\eta}$, on a la d{\'e}composition 
$\g=\g^{\xi+v}\oplus \n \oplus
[\eta,\g]$, car $\alpha(x)=\g^{\sigma(x)}$, pour tout $x$ dans
l'intersection $Y' \cap X_{r}$. On en d{\'e}duit que l'{\'e}l{\'e}ment
\begin{eqnarray*}
& \sum\limits_{j=1}^{s}p_{j}(y_{2},\ldots,y_{n})[\xi+u_{1}+y_{2}u_{2}+\cdots+y_{n}u_{n},v_{j}] \\
& \hspace{3.5cm} +\sum\limits_{k=1}^{t}q_{k}(y_{2},\ldots,y_{n})
[\xi+u_{1}+y_{2}u_{2}+\cdots+y_{n}u_{n},f_{k}]
\end{eqnarray*}
est orthogonal {\`a} $\g$, donc est nul. D{\'e}signons par $d$ le plus grand des degr{\'e}s des polyn{\^o}mes
$p_{1},\ldots,p_{s},q_{1},\ldots,q_{t}$ et par $\chi$ la fonction
polynomiale sur $\g^{\eta}$\,:
$$x_{1}^{d}\sum\limits_{j=1}^{s}
p_{j}\left( \frac{x_{2}}{x_{1}},\ldots,\frac{x_{n}}{x_{1}} \right) v_{j}+
x_{1}^{d}\sum\limits_{k=1}^{t}
q_{k}\left( \frac{x_{2}}{x_{1}},\ldots,\frac{x_{n}}{x_{1}} \right) f_{k} \; \cdot$$  
Alors $\chi(z)$ centralise $\xi+x_{1}^{-1}z$ pour tout
$z=x_{1}u_{1}+ \cdots + x_{n} u_{n}$ dans $\g^{\eta}$, avec $x_{1}$
non nul. Or
pour tout $v$ dans un ouvert non vide de $\g^{\eta}$, les sous-espaces
$\g^{\xi+v}$ et $\n \oplus [\eta,\g]$ de $\g$ ont une intersection
nulle. Par suite $\chi$ est nulle. Ceci est absurde car les
polyn{\^o}mes $p_{1},\ldots,p_{s},q_{1},\ldots,q_{t}$ ne sont pas tous
nuls.\\

L'application $\zeta$ n'est donc pas identiquement nulle sur $Y
\setminus X_{r}$ et par cons{\'e}quent $\Delta$ n'est pas un {\'e}l{\'e}ment nul de
${\cal S}(\g)$. Ceci termine la preuve du lemme.
\qed
\\
On en d{\'e}duit le r{\'e}sultat suivant\,:
\begin{prop}\label{inclusion}
On a les inclusions\,:
$$\z(\g^{\xi}) \subset \alpha(x) \subset \g^{\xi}, {\rm \; pour \;
  tout \; } x {\rm \ dans \ } Y \setminus X_{r} \cdot$$
\end{prop}
\dem La deuxi{\`e}me inclusion est connue. On s'int{\'e}resse {\`a} la premi{\`e}re. On fixe un point $x_{0}$ dans $Y \setminus X_{r}$ et il
s'agit de montrer l'inclusion\,:
$$\z(\g^{\xi}) \subset \alpha(x_{0}) \cdot$$
On suppose par l'absurde que $\z(\g^{\xi})$ n'est pas contenu dans $\alpha(x_{0})$. Il existe alors
un vecteur $v'$ de  $\z(\g^{\xi})$ qui n'appartient par {\`a}
$\alpha(x_{0})$ et on construit un suppl{\'e}mentaire $\n'$ de
$\alpha(x_{0})$ dans $\g^{\xi}$ admettant une base
$\{v_{1}',\ldots,v_{s}'\}$ telle que $v_{1}'=v'$. La matrice
$([v_{i}',v_{j}'])_{1 \leq i,j \leq s}$ a une premi{\`e}re colonne nulle
car $v_{1}'=v'$ appartient au centre $\z(\g^{\xi})$ de
$\g^{\xi}$. Ceci contredit le lemme \ref{supplementaire} appliqu{\'e} au
suppl{\'e}mentaire $\n=\n'$ de $\alpha(x_{0})$ dans $\g^{\xi}$.
\qed

Soit $e_{1},\ldots,e_{m}$ une base de $\z(\g^{\xi})$. D'apr{\`e}s la
proposition pr{\'e}c{\'e}dente, on a l'inclusion $\z(\g^{\xi}) \subset \alpha(x)$, pour tout
$x$ dans $Y \setminus X_{r}$. Ceci permet de construire une base
$\epsilon_{1}(x),\ldots,\epsilon_{{\rm rg \; } \g}(x)$ de $\alpha(x)$
pour $x$ dans $Y$ qui v{\'e}rifie
$$\epsilon_{i}(x)=e_{i}, {\rm \ pour \ tout \ } x {\rm \ dans \ } Y \setminus
X_{r} {\rm \ et \ } i=1,\ldots,m \cdot$$
Il existe donc des applications r{\'e}guli{\`e}res
$\widetilde{\epsilon_{1}},\ldots,\widetilde{\epsilon_{m}}$ sur $Y$, non identiquement nulles
sur $Y \setminus X_{r}$, qui v{\'e}rifient les relations\,:
$$\epsilon_{i}(x)=e_{i}+q(x)^{m_{i}} \widetilde{\epsilon_{i}}(x), {\rm
  \ pour \ tout \ } x  {\rm \ dans \ } Y  {\rm \ et \ } i=1,\ldots,m,$$
o{\`u} $m_{i}$ est un entier strictement positif.
\begin{lemme}\label{mi}
Pour tout $i$ de $\{1,\ldots,m\}$, l'entier $m_{i}$ est {\'e}gal {\`a} 1.
\end{lemme} 
\dem On suppose par l'absurde $m_{i} >1$ pour un $i$ dans $\{1,\ldots,m\}$. Il s'agit d'aboutir {\`a} une
contradiction. On raisonne comme dans la d{\'e}monstration du lemme
\ref{supplementaire} et on reprend les m{\^e}mes notations. Puisque l'image de $\epsilon_{i}=e_{i}+q^{m_{i}} \widetilde{\epsilon_{i}}$ par
$\gamma'_{\nu}$ est {\'e}gale {\`a} $e_{i}+\tau^{m_{i}} q'(x)(v')
\widetilde{\epsilon_{i}}(x)$, de l'{\'e}galit{\'e}
$$[\sigma(x),\epsilon_{i}(x)]=0, $$
on tire l'{\'e}galit{\'e}
$$[\xi+ \tau v + \cdots, e_{i}+\tau^{m_{i}} q'(x)(v')
\widetilde{\epsilon_{i}}(x)]=0,$$
car l'image de $\sigma$ par $\gamma'_{\nu}$ est {\'e}gale {\`a} $\xi + \tau v
{\rm \; mod \;} \tau^{2}$.  Puisque $m_{i}>1$, le terme de degr{\'e} $1$ en $\tau$ du membre de gauche de cette {\'e}galit{\'e}
est nul; donc $[v,e_{i}]$ est nul. Il en r{\'e}sulte que pour tout $v$ dans un ouvert non vide de
$\g^{\eta}$, $[v,e_{i}]$ est nul. Par suite, $[v,e_{i}]$ est nul pour
tout $v$ dans $\g^{\eta}$. Ceci est absurde car $e_{i}$ est un {\'e}l{\'e}ment
de $\z(\g^{\xi})$ et $\z(\g^{\eta})$ et
$\z(\g^{\xi})$ ont une intersection nulle. 
\qed

Du lemme pr{\'e}c{\'e}dent, on d{\'e}duit la relation $\epsilon_{i}(x)=e_{i}+q(x)
\widetilde{\epsilon_{i}}(x)$, pour tout $x$ dans $Y$ et tout $i$ dans
$\{1,\ldots,m\}$. D'apr{\`e}s la proposition \ref{centralisateur}, on a la
d{\'e}composition\,: $\g={\g^{\xi}}^{\perp} \oplus \g^{\eta}$. On
note, pour $x$ dans $Y$ et $i$ dans $\{1,\ldots,m\}$,  $\widetilde{\epsilon_{i,1}}(x)$ et  $\widetilde{\epsilon_{i,2}}(x)$
les composantes de  $\widetilde{\epsilon_{i}}(x)$ sur ${\g^{\xi}}^{\perp}$
et $\g^{\eta}$ respectivement.

\begin{prop}\label{libre}
Pour tout $x$ dans un ouvert non vide de $Y \setminus X_{r}$,
l'ensemble des {\'e}l{\'e}ments $\widetilde{\epsilon_{1,2}}(x),\ldots,\widetilde{\epsilon_{m,2}}(x)$
est une partie libre de $\g$.
\end{prop}
\dem Supposons l'assertion fausse. Il s'agit d'aboutir {\`a} une
contradiction. Il existe des fonctions r{\'e}guli{\`e}res $a_{1},\ldots,a_{m}$
sur $Y$, non toutes
identiquement nulles sur $Y \setminus X_{r}$, qui
satisfont l'{\'e}galit{\'e}\,:
$$a_{1}(x) \widetilde{\epsilon_{1,2}}(x) + \cdots + a_{m}(x) \widetilde{\epsilon_{m,2}}(x)=0
\; ,$$
pour tout $x$ dans $Y \setminus X_{r}$. Soit $\mu$ l'application
r{\'e}guli{\`e}re sur $Y$ {\`a} valeurs dans $\z(\g^{\xi})$ qui {\`a}
$x$ dans $Y$ associe $a_{1}(x) e_{1} + \cdots + a_{m}(x)
e_{m}$.\\

On reprend les notations de la
d{\'e}monstration du lemme \ref{inclusion} avec $m_{i}=1$, autrement dit
$\tau=\tau_{0}$. Pour $v$ dans $\g^{\eta}$ non nul, on note $[v]$ l'image de $v$ dans
$\mathbb{P}(\g^{\eta})$. L'ensemble des points $v$ de $\g^{\eta}$ pour
lesquels le point $\xi \times [v]$ de $(\xi + \g^{\eta}) \times
\mathbb{P}(\g^{\eta})$ appartient {\`a} $Y \setminus X_{r}$ est un ouvert
non vide de $\g^{\eta}$. Soit $v$ dans $\g^{\eta}$ appartenant {\`a} cet
ouvert et posons $x=\xi \times [v] \in Y \setminus X_{r}$. Alors il existe un
vecteur $v'$ tangent {\`a} $Y$ en $x$ tel que l'image de $v'$ par l'application
lin{\'e}aire tangente {\`a} $\sigma$ en $x$ soit {\'e}gale {\`a} $v$. En effet, on
consid{\`e}re la courbe $\gamma$ qui au complexe $t$ associe le point $(\xi
+tv) \times [v]$ de $(\xi + \g^{\eta}) \times
\mathbb{P}(\g^{\eta})$. Pour $t$ dans un voisinage ouvert de
$0$, la courbe $\gamma$ est {\`a} valeurs dans $Y$ et l'on a $\gamma(0)=x$;
par suite la d{\'e}riv{\'e}e $\gamma'(0)$ de $\gamma$ en
$0$ est un vecteur de l'espace tangent {\`a} $Y$ en $x$. Or on v{\'e}rifie
ais{\'e}ment que $\sigma'(x)(\gamma'(0))=v$, d'o{\`u} l'existence du vecteur
$v'$.

Puisque l'image de \hbox{$\sum\limits_{j=1}^{m} a_{j} \epsilon_{j}=
\sum\limits_{j=1}^{m} a_{j} e_{j}+q \sum\limits_{j=1}^{m} a_{j}
\widetilde{\epsilon_{j}}$} par $\gamma'_{\nu}$ est {\'e}gale {\`a}
\hbox{$\sum\limits_{j=1}^{m} a_{j}(x) e_{j} +\tau \sum\limits_{j=1}^{m} a_{j}'(x)(v') e_{j}+\tau q'(x)(v')
\sum\limits_{j=1}^{m} a_{j}(x) \widetilde{\epsilon_{j}}(x)$}, de
l'{\'e}galit{\'e} 
$$[\sigma(x),\sum\limits_{j=1}^{m} a_{j}(x) \epsilon_{j}(x)]=0,$$
on tire l'{\'e}galit{\'e}
$$[\xi+\tau v,\sum\limits_{j=1}^{m} a_{j}(x) e_{j} +\tau \sum\limits_{j=1}^{m} a_{j}'(x)(v') e_{j}+\tau q'(x)(v')
\sum\limits_{j=1}^{m} a_{j}(x) \widetilde{\epsilon_{j,1}}(x)]=0 \cdot$$
Cela provient des deux {\'e}galit{\'e}s suivantes\,: $\gamma'_{\nu}(\sigma)=\xi+\tau v$ et 
$\sum\limits_{j=1}^{p} a_{j}(x) \widetilde{\epsilon_{j,2}}(x)=0$. Le terme de degr{\'e} $1$ en $\tau$ du
membre de gauche de l'{\'e}galit{\'e} pr{\'e}c{\'e}dente est nul, d'o{\`u} l'{\'e}galit{\'e}
$$[\xi, q'(x)(v')\sum\limits_{j=1}^{m}
a_{j}(x)\widetilde{\epsilon_{i,1}}(x)]=-[v,\sum\limits_{j=1}^{m}
a_{j}(x) e_{j}] \cdot$$

Les applications $\widetilde{\epsilon_{i,1}}$ sont {\`a} valeurs dans
${\g^{\xi}}^{\perp}=[\xi,\g]$ et ce qui pr{\'e}c{\`e}de donne la relation\,:
$$[v,\sum\limits_{j=1}^{m} a_{j}(x) e_{j}] \in [\xi,[\xi,\g]]
=(\g^{\xi} \oplus [\eta,\g^{\xi}])^{\perp},$$
d'apr{\`e}s la proposition \ref{perp}, ce qui signifie encore, par d{\'e}finition de l'application $\mu$,
\begin{eqnarray*}
\langle v, [[\eta,\g^{\xi}],\mu(x)] \rangle = \{0\} \cdot
\end{eqnarray*}
En r{\'e}sum{\'e}, on a montr{\'e}\,: pour tout $\xi \times [v]$ dans $Y
\setminus X_{r}$, on a la relation
\begin{eqnarray}\label{relation_nu}
v \in {[[\eta,\g^{\xi}],\mu(\xi \times [v])]}^{\perp}.
\end{eqnarray}
On note
$$Y'=\{ x \in Y \setminus X_{r} \ | \ \mu(x) \not= 0 \} \cdot$$
L'ensemble $Y'$ est un ouvert non vide de $Y \setminus X_{r}$ et comme
$\xi$ v{\'e}rifie la propri{\'e}t{\'e} $(P)$, on a l'inclusion\,:
$${[[\eta,\g^{\xi}],\mu(x)]}^{\perp} \subset {W}^{\perp},$$
pour tout $x$ dans $Y'$. Comme le sous-espace $W$ est
contenu dans $\g^{\xi}$, son orthogonal $W^{\perp}$ ne contient pas
$\g^{\eta}$; sinon, ce dernier contiendrait la somme $\g^{\eta} \oplus
{\g^{\xi}}^{\perp}=\g$ et $W$ serait nul. L'ensemble
$$X^{\circ}=\{ \xi \times [v] \in X \setminus X_{r} \ | \ v \not\in
{W}^{\perp} \}$$
est donc un ouvert non vide de $X \setminus X_{r}$ et l'intersection $Y'
\cap X^{\circ}$ est non vide. Il en r{\'e}sulte que pour tout $x=\xi \times
[v]$ dans l'intersection $Y' \cap X^{\circ}$, le sous-espace 
${[[\eta,\g^{\xi}],\mu(x)]}^{\perp}$ n'est pas contenu dans
${W}^{\perp}$, d'apr{\`e}s la relation (\ref{relation_nu}). Ceci contredit
le fait que $x$ appartient {\`a} $Y'$. 
\qed
\\
On en d{\'e}duit le corollaire suivant\,:
\begin{cor}\label{rang_z'}
On note $e_{1},\ldots,e_{n}$ une base de $\g^{\xi}$ telle que
$e_{1},\ldots,e_{m}$ est une base de $\z(\g^{\xi})$. On suppose que $\xi$ v{\'e}rifie la propri{\'e}t{\'e} $(P)$. Alors la matrice 
$([[\eta,e_{i}],e_{j}])_{1 \leq i \leq n \atop  1 \leq j \leq m}$ 
{\`a} coefficients dans ${\cal S}(\g^{\xi})$ est de rang maximal {\'e}gal {\`a}
$\dim \z(\g^{\xi})$.
\end{cor}
\dem On
{\'e}value la matrice $([[\eta,e_{i}],e_{j}])_{1 \leq i \leq n \atop  1
  \leq j \leq m}$ en la forme lin{\'e}aire $\langle \sigma(x), \cdot
\rangle$ sur $\g^{\xi}$, pour $x$ dans $Y$\,:
\begin{eqnarray*}
(\langle \sigma(x),[[\eta,e_{i}],e_{j}] \rangle)_{1 \leq i \leq n \atop  1
  \leq j \leq m} & = & 
(\langle [\sigma(x),e_{j}],[\eta,e_{i}] \rangle)_{1 \leq i \leq n \atop  1
  \leq j \leq m}\\  
& = & -q(x)(\langle [\sigma(x),\widetilde{\epsilon_{j}}(x)],[\eta,e_{i}] \rangle)_{1 \leq i \leq n \atop  1
  \leq j \leq m}.
\end{eqnarray*}
On suppose par l'absurde que, pour tout $x$ dans $Y \setminus X_{r}$, la matrice $(\langle [\sigma(x),\widetilde{\epsilon_{j}}(x)],[\eta,e_{i}] \rangle)_{1 \leq i \leq n \atop  1
  \leq j \leq m}$ n'est pas de rang maximal; autrement dit on suppose
  que, pour tout $x$ dans $Y \setminus X_{r}$, les
  vecteurs {\`a} coefficients complexes
$$\left[
\begin{array}{c}
\langle [\sigma(x),\widetilde{\epsilon_{1}}(x)],[\eta,e_{1}] \rangle \\
\vdots \\
\langle [\sigma(x),\widetilde{\epsilon_{1}}(x)],[\eta,e_{n}] \rangle
\end{array}
\right], \ldots, 
\left[
\begin{array}{c}
\langle [\sigma(x),\widetilde{\epsilon_{m}}(x)],[\eta,e_{1}] \rangle  \\
\vdots \\
\langle [\sigma(x),\widetilde{\epsilon_{m}}(x)],[\eta,e_{n}] \rangle 
\end{array}
\right],$$
sont li{\'e}s. Il existe donc des fonctions r{\'e}guli{\`e}res $a_{1},\ldots,a_{m}$ sur $Y \setminus X_{r}$, non toutes identiquement nulles, qui satisfont l'{\'e}galit{\'e}\,:
$$\langle [\sigma(x),\sum\limits_{j=1}^{m} a_{j}(x)
\widetilde{\epsilon_{j}}(x)], [\eta,e_{i}] \rangle=0,$$
pour $i=1,\ldots,n$ et $x$ dans $Y \setminus X_{r}$. Notons
$\chi$ l'application $\chi(x)=\sum\limits_{j=1}^{m} a_{j}(x)
\widetilde{\epsilon_{j}}(x)$. L'inclusion $Y \setminus X_{r} \subset
\sigma^{-1}(\{\xi\})$ donne la relation
$$\langle \chi(x), [\xi,[\eta,\g^{\xi}]] \rangle =\{0\},$$
pour tout $x$ dans $Y \setminus X_{r}$. Puisque
$[\xi,[\eta,\g^{\xi}]]=\g^{\xi}$, on en d{\'e}duit que $\chi(x)$ appartient
{\`a} ${\g^{\xi}}^{\perp}$, pour tout $x$ dans $Y \setminus X_{r}$. Comme $\chi(x)$ s'{\'e}crit $\chi(x)=\sum\limits_{i=1}^{m} a_{i}(x)
\widetilde{\epsilon_{i,1}}(x) +\sum\limits_{i=1}^{m} a_{i}(x)
\widetilde{\epsilon_{i,2}}(x)$, avec $\sum\limits_{i=1}^{m} a_{i}(x)
\widetilde{\epsilon_{i,1}}(x)$ dans ${\g^{\xi}}^{\perp}$, l'{\'e}l{\'e}ment $\sum\limits_{i=1}^{} a_{i}(x)
\widetilde{\epsilon_{i,2}}(x)$ appartient {\`a} l'intersection
${\g^{\xi}}^{\perp} \cap \g^{\eta}$, donc est nul. D'apr{\`e}s la proposition \ref{libre}, les {\'e}l{\'e}ments
$\widetilde{\epsilon_{i,2}}(x)$ sont lin{\'e}airement ind{\'e}pendants pour tout $x$ dans
un ouvert non vide de $Y \setminus X_{r}$. On en d{\'e}duit que les
fonctions $a_{i}$ sont toutes identiquement nulles sur $Y \setminus
X_{r}$, ce qui contredit les hypoth{\`e}ses. \\

Par cons{\'e}quent la matrice $(\langle [\sigma(x),\widetilde{\epsilon_{j}}],[\eta,e_{i}] \rangle)_{1 \leq i \leq n \atop  1
  \leq j \leq m}$ {\`a} coefficients dans ${\cal S}(\g^{\xi})$ est de rang
  maximal, pour tout $x$ dans un ouvert non vide de $Y \setminus X_{r}$. 
\qed 

On en d{\'e}duit le th{\'e}or{\`e}me \ref{intro2} lorsque $\xi$ v{\'e}rifie la propri{\'e}t{\'e}
$(P)$ de la d{\'e}finition \ref{prop_P}. En effet, avec les notations de
la proposition \ref{rang}, la matrice du corollaire pr{\'e}c{\'e}dent est la matrice 
$\left[
\begin{array}{c}
\mathfrak{D} \\
\mathfrak{E}
\end{array}
\right]$. Cela r{\'e}sulte des relations, d{\'e}j{\`a} vues, 
$[[\eta,e_{j}],e_{i}]=[[\eta,e_{i}],e_{j}]$, pour $i=1,\ldots,n$ et
$j=1,\ldots,m$. En outre le th{\'e}or{\`e}me \ref{intro2} est connu pour les
{\'e}l{\'e}ments nilpotents r{\'e}guliers de $\g$.\\
\\
\rem Sous les hypoth{\`e}ses et avec les notations de la partie 3, il
n'est pas difficile de montrer que les {\'e}l{\'e}ments nilpotents de $\g$
v{\'e}rifient la propri{\'e}t{\'e} $(P)$ lorsque $\z'=\z(\g^{\xi})$. En revanche,
les {\'e}l{\'e}ments nilpotents de $\mathfrak{so}_{2n}$ ne v{\'e}rifient pas la
propri{\'e}t{\'e} $(P)$ en g{\'e}n{\'e}ral.

\section{{\'E}tude de la propri{\'e}t{\'e} $(P)$ et d{\'e}monstration des th{\'e}or{\`e}mes
  \ref{intro1} et \ref{intro2} dans le cas exceptionnel.}

Pour terminer la preuve des th{\'e}or{\`e}mes \ref{intro1} et \ref{intro2}, il
reste essentiellement {\`a} {\'e}tudier la propri{\'e}t{\'e} $(P)$ pour les alg{\`e}bres de Lie simples
exceptionnelles. On suppose que $\g$ est une alg{\`e}bre de Lie
simple exceptionnelle  et
que $\xi$ est un
{\'e}l{\'e}ment nilpotent distingu{\'e} non r{\'e}gulier de $\g$. Il s'agit de montrer
que l'{\'e}l{\'e}ment $\xi$ v{\'e}rifie la propri{\'e}t{\'e} $(P)$. On note
$m_{1},\ldots, m_{r}$ les valeurs propres de la restriction de ${\rm
  ad} \rho$ au sous-espace $\z(\g^{\xi})$. Les entiers
$m_{1},\ldots,m_{r}$ sont pairs et on a 
$$2=m_{1} < m_{2} < \cdots < m_{r} \cdot$$
On note $V(m_{l})$ le sous-espace propre correspondant {\`a} la valeur
propre $m_{l}$, pour $l=1,\ldots,r$. Avec les notations de la
d{\'e}finition \ref{prop_P}, on a\,: 
$$W=V(m_{r}) \cdot$$ 
On choisit une base 
$${\cal B}=
\xi_{m_{1}}^{1},\ldots,\xi_{m_{1}}^{d_{1}},\xi_{m_{2}}^{1},\ldots
,\xi_{m_{2}}^{d_{2}}, \ldots,\xi_{m_{r}}^{1},\ldots,\xi_{m_{r}}^{d_{r}}$$
de $\g^{\xi}$ de vecteurs propres telle que
$\xi_{m_{l}}^{1},\ldots,\xi_{m_{l}}^{d_{l}}$ forme une base de
$V(m_{l})$, pour $l=1,\ldots,r$ et telle qu'il existe une base de
$\z(\g^{\xi})$ form{\'e}e de vecteurs de ${\cal B}$. On peut supposer que
$\xi_{m_{1}}^{1}=\xi_{2}^{1}=\xi$. On fixe un {\'e}l{\'e}ment non nul $v$ de
$\z(\g^{\xi})$ et on souhaite montrer que le sous-epsace $V(m_r)$ est
contenu dans le sous-espace $[[\eta,\g^{\xi}],v]$. Commen{\c c}ons par prouver la proposition suivante\,:
\begin{lemme}\label{composante_xi} Si la coordonn{\'e}e de $v$ en
  $\xi_{m_{1}}^{1}=\xi$ est non nulle, alors le sous-epsace $V(m_r)$ est
contenu dans le sous-espace $[[\eta,\g^{\xi}],v]$.  
\end{lemme}
\dem Soit $v$ un {\'e}l{\'e}ment non nul de $\z(\g^{\xi})$ dont la
  coordonn{\'e}e $\lambda$ de $v$ en $\xi_{m_{1}}^{1}=\xi$ est non
  nulle. Soit $p$ dans $\{1,\ldots,d_{r}\}$. L'{\'e}l{\'e}ment
  $\xi_{m_{r}}^{p}$ de $V(m_{r})$ s'{\'e}crit
$$\xi_{m_{r}}^{p}=[[\eta,-\displaystyle{\frac{1}{m_{r}}}\xi_{m_{r}}^{p}],\xi]
\cdot$$
D'apr{\`e}s \cite{Panyushev}, (Th{\'e}or{\`e}me 2.3), $\xi$ est le seul {\'e}l{\'e}ment de
$\z(\g^{\xi})$ de poids 2. Comme $m_{r}$ est le plus haut poids, on en
d{\'e}duit que l'{\'e}l{\'e}ment $\xi_{m_{r}}^{p}$ s'{\'e}crit sous la forme\,:
$$\xi_{m_{r}}^{p}=[[\eta,-\displaystyle{\frac{1}{m_{r} \times
    \lambda}} \xi_{m_{r}}^{p}],v] \cdot$$
Par suite, on a l'inclusion\,: 
$$V(m_{r}) \subseteq [[\eta,\g^{\xi}],v],$$
d'o{\`u} le lemme.
\qed 

Soit $i_{1} < \cdots < i_{p}$ dans
$\{1,\ldots,r\}$ et
$k_{(1,1)},\ldots,k_{(1,\delta_{1})},\ldots,k_{(p,1)},\ldots,k_{(p,\delta_{p})}$
des indices tels que les {\'e}l{\'e}ments 
$$\xi_{m_{i_{1}}}^{k_{(1,1)}},\ldots
,\xi_{m_{i_{1}}}^{k_{(1,\delta_{1})}},\ldots,
\xi_{m_{i_{p}}}^{k_{(p,1)}},\ldots
,\xi_{m_{i_{p}}}^{k_{(p,\delta_{p})}}$$  
appartiennent {\`a} $\z(\g^{\xi})$ et tels que $v$ s'{\'e}crive sous la
forme  
$$v=\alpha_{1}^{1}\xi_{m_{i_{1}}}^{k_{(1,1)}}+\cdots
+\alpha_{1}^{\delta_{1}}\xi_{m_{i_{1}}}^{k_{(1,\delta_{1})}}+\cdots+
\alpha_{p}^{1}\xi_{m_{i_{p}}}^{k_{(p,1)}}+\cdots
+\alpha_{p}^{\delta_{p}}\xi_{m_{i_{p}}}^{k_{(p,\delta_{p})}},$$
avec
$(\alpha_{1}^{1},\ldots,\alpha_{1}^{\delta_{1}},\ldots,\alpha_{p}^{1},\ldots,\alpha_{p}^{\delta_{p}})$
dans $(\C^{*})^{K}$ o{\`u} $K=\delta_{1}+\cdots+\delta_{p}$. D'apr{\`e}s la
proposition pr{\'e}c{\'e}dente, on peut supposer que $i_{1}$ est strictement
plus grand que $1$.\\

Pour $i,j$ et $k$ dans $\{1,\ldots,r\}$ et $s,p$ et $q$ dans
$\{1,\ldots,d_{i}\}$, $\{1,\ldots,d_{j}\}$ et $\{1,\ldots,d_{k}\}$
respectivement, on note $\lambda_{(m_{k},q),(m_{i},s),(m_{j},p)}$ la
coordonn{\'e}e de l'{\'e}l{\'e}ment $[[\eta,\xi_{m_{k}}^{q}],\xi_{m_{i}}^{s}]$ en
$\xi_{m_{j}}^{p}$. Notons que si
$\lambda_{(m_{k},q),(m_{i},s),(m_{j},p)}$ est non nul, on a la
relation $m_{k}=m_{j}-m_{i}+2$.\\

On suppose qu'il existe $k_{1}$ dans
$\{1,\ldots,r\}$ tel que $m_{k_{1}}=m_{r}-m_{i_{1}}+2$. Soit
$w=\sum\limits_{p=1}^{d_{r}} b_{p}\xi_{m_{r}}^{p}$ un {\'e}l{\'e}ment de $V({m_{r}})$. On cherche $u$
dans $\g^{\xi}$ sous la forme $u=\sum\limits_{q=1}^{d_{k_{1}}}
a_{q}\xi_{m_{k_{1}}}^{q}$ tel que 
$$[[\eta,u],v]=w \cdot$$
Ce probl{\`e}me revient {\`a} r{\'e}soudre un syst{\`e}me lin{\'e}aire d'inconnue le
vecteur 
$\left[
\begin{array}{c}
a_{1} \\
a_{2} \\
\vdots \\
a_{d_{k_{1}}}
\end{array}
\right]
$, de second membre 
$\left[ 
\begin{array}{c}
b_{1} \\
b_{2} \\
\vdots \\
b_{d_{r}}
\end{array}
\right]
$ 
et de matrice associ{\'e}e la matrice $M(v)$ de taille $d_{r} \times d_{k_{1}}$
dont les coefficients $M(v)_{p,q}$ sont donn{\'e}s par\,:
$$M(v)_{p,q}=\sum_{s=1}^{\delta_{i_{1}}} \alpha_{1}^{s}
\lambda_{(m_{k_{1}},q),(m_{i_{1}},s),(m_{r},p)},$$
pour $p$ dans $\{1,\ldots,d_{r}\}$ et $q$ dans
$\{1,\ldots,d_{k_{1}}\}$. La matrice $M(v)$ s'{\'e}crit aussi de mani{\`e}re
plus agr{\'e}able comme une somme de matrices\,:
$$M(v)=\sum_{s=1}^{\delta_{i_{1}}} \alpha_{1}^{s}
M(\xi_{m_{i_{1}}}^{k_{(1,s)}}) \cdot$$
Il est clair que si la matrice $M(v)$ est surjective, le sous-espace
$V(m_r)$ est contenu dans le sous-epsace $[[\eta,\g^{\xi}],v]$. On s'int{\'e}resse d{\'e}sormais {\`a} cette
condition. Remarquons que cette condition ne d{\'e}pend que des termes de
plus bas poids $m_{i_{1}}$ intervenant dans l'{\'e}criture de $v$.\\

Pour chaque
alg{\`e}bre de Lie exceptionnelle, on trouve dans
\cite{Triplets} une liste de $\mathfrak{sl}_{2}$-triplets
correspondant aux orbites nilpotentes. Le logiciel GAP4 permet en outre d'effectuer
des calculs dans les alg{\`e}bres de Lie. Il permet notamment de calculer
le centralisateur d'un {\'e}l{\'e}ment, le centre d'une sous-alg{\`e}bre,
etc. L'{\'e}tude pr{\'e}c{\'e}dente donne les conditions {\`a} v{\'e}rifier
pour que tout {\'e}l{\'e}ment non nul de $\z(\g^{\xi})$ v{\'e}rifie $(P)$. Gr{\^a}ce
aux donn{\'e}es de \cite{Triplets} et au logiciel GAP4, on exhibe 
une base ${\cal B}$ v{\'e}rifiant les conditions pr{\'e}c{\'e}dentes. Si $i_1$ est un
entier de $\{1,\ldots,r\}$ tel que $V(m_{i_1})$ et $\z(\g^{\xi})$
ont une intersection non nulle et si les {\'e}l{\'e}ments 
$\xi_{m_{i_1}}^{l_1},\ldots,\xi_{m_{i_1}}^{l_{\delta}}$ 
forment une base de l'intersection $V(m_{i_1}) \cap
\z(\g^{\xi})$, on v{\'e}rifie qu'il existe un entier $k_{1}$ tel que $m_{k_1}=m_{j_1}-m_{i_1}+2$, on calcule {\`a} l'aide de GAP4 les matrices 
$M(\xi_{m_{i_1}}^{l_1}),\ldots
,M(\xi_{m_{i_1}}^{l_{\delta}})$ et on v{\'e}rifie la surjectivit{\'e} de
la matrice 
$$M(v)=\sum_{s=1}^{\delta} \alpha_{l_s}
M(\xi_{m_{i}}^{l_s}),$$
pour tout $\delta$-uplet $(\alpha_{l_1},\ldots,\alpha_{l_{\delta}})$
non nul. Le lemme \ref{composante_xi} permet de ne traiter que
le cas o{\`u} $m_{i_1}$ est strictement plus grand que $2$. De plus, si
$m_{i_1}=m_{r}$ et $\delta=1$, alors on a \,: 
$1=\delta=\dim V(m_{i_1}) \cap \z(\g^{\xi})=\dim V(m_r)$ et on a clairement l'inclusion $V(m_r) \subseteq
[[\eta,\g^{\xi}],v]$, pour $v=\sum_{s=1}^{\delta} \alpha_{l_s}
\xi_{m_{i_1}}^{l_s}$. 

On donne en annexe les calculs permettant de v{\'e}rifier ces conditions
et par suite la propri{\'e}t{\'e} $(P)$. On s'aper{\c c}oit que dans la plupart des cas la somme
pr{\'e}c{\'e}dente n'a qu'un seul terme. L'{\'e}tude de la surjectivit{\'e} ne d{\'e}pend
alors d'aucun param{\`e}tre, ce qui facilite le travail.\\ 

Ces calculs montrent que les alg{\`e}bres $E_{6}$, $E_{7}$, $E_{8}$ et
$F_{4}$ v{\'e}rifient la propri{\'e}t{\'e} $(P)$. Reste le cas de $G_{2}$\,; dans $G_{2}$, il n'y a qu'une
seule orbite nilpotente distingu{\'e}e non r{\'e}guli{\`e}re et pour cette orbite
le centre est de dimension $2$. Le lemme 
\ref{composante_xi} assure que cette orbite v{\'e}rifie $(P)$. Par suite $G_{2}$
v{\'e}rifie $(P)$ aussi. Des propositions \ref{rang} et \ref{distingue} et
du corollaire \ref{rang_z'}, il r{\'e}sulte que
la relation du th{\'e}or{\`e}me \ref{intro2} est satisfaite pour ces cinq 
alg{\`e}bres de Lie. Notons ici l'importance d'avoir d{\'e}montr{\'e} le th{\'e}or{\`e}me
\ref{intro2} dans le cas classique avant le cas exceptionnel pour
appliquer la proposition \ref{distingue}. En
effet, les sous-alg{\`e}bres de Lie semi-simples d'une alg{\`e}bre de Lie exceptionnelle
sont isomorphes {\`a} des
produits finis d'id{\'e}aux simples et ces composantes simples peuvent
{\^e}tre isomorphes {\`a} des alg{\`e}bres de Lie de type $A$, $B$, $C$ ou $D$. La d{\'e}monstration des th{\'e}or{\`e}mes \ref{intro1} et \ref{intro2}
est ainsi achev{\'e}e.
 
\section*{Annexe} On pr{\'e}sente dans cette annexe les calculs faits {\`a}
l'aide du logiciel GAP4 permettant
de v{\'e}rifier la propri{\'e}t{\'e} $(P)$ pour $E_{6}$, $E_{7}$, $E_{8}$ et
$F_{4}$. On reprend les notations de la partie pr{\'e}c{\'e}dente et la d{\'e}marche g{\'e}n{\'e}rale est la suivante\,: on d{\'e}finit l'alg{\`e}bre
de Lie {\tt L} dans laquelle on veut travailler gr{\^a}ce {\`a} la commande {\tt
  SimpleLieAlgebra}, on d{\'e}finit un syst{\`e}me de racines ({\tt
  RootSystem}), un syst{\`e}me de racines positives correspondant ({\tt
  PositiveRoots}) puis des syst{\`e}mes de vecteurs \guillemotleft
positifs\guillemotright \ et  \guillemotleft n{\'e}gatifs\guillemotright \ 
associ{\'e}s ({\tt PositiveRootVectors} et {\tt
  NegatitiveRootVectors}). La commande {\tt CanonicalGenerators} donne
une base de la sous-alg{\`e}bre de Cartan. On peut d{\'e}sormais faire des calculs dans
l'alg{\`e}bre de Lie {\tt L}. Il s'agit ensuite d'{\'e}tudier les orbites nilpotentes distingu{\'e}es non
r{\'e}guli{\`e}res. Pour chacune d'entre elles, on d{\'e}finit un
$\mathfrak{sl}_{2}$-triplet {\tt $\{$e,h,f$\}$} gr{\^a}ce aux donn{\'e}es de
\cite{Triplets}. On calcule ensuite le centralisateur {\tt g} de l'{\'e}l{\'e}ment positif {\tt e} avec la
commande {\tt LieCentralizer} puis le centre {\tt z} du centralisateur avec
{\tt LieCentre}. Pour chaque cas, on pr{\'e}cise la valeur du plus haut
poids $m_{r}$ et on donne le nombre de matrices {\`a} {\'e}tudier. Ensuite,
pour chacune d'entre elles on donne les valeurs de $m_{i_{1}}$ et
$m_{k_{1}}$ et on fait les calculs nec{\'e}ssaires. Les calculs de la premi{\`e}re orbite de la premi{\`e}re alg{\`e}bre (il s'agit de
l'orbite sous-r{\'e}guli{\`e}re de $E_{6}$) sont d{\'e}taill{\'e}s; les autres le sont
un peu moins.

\subsection*{Cas de $E_{6}$.}
\noindent On commence par d{\'e}finir {\tt L} et les g{\'e}n{\'e}rateurs de {\tt L}\,:
\begin{verbatim}
> L:=SimpleLieAlgebra("E",6,Rationals);
<Lie algebra of dimension 78 over Rationals>
> R:=RootSystem(L);
<root system of rank 6>
> P:=PositiveRoots(R);;
> x:=PositiveRootVectors(R);
[ v.1, v.2, v.3, v.4, v.5, v.6, v.7, v.8, v.9, v.10, v.11, v.12, v.13, v.14, 
  v.15, v.16, v.17, v.18, v.19, v.20, v.21, v.22, v.23, v.24, v.25, v.26, 
  v.27, v.28, v.29, v.30, v.31, v.32, v.33, v.34, v.35, v.36 ]
> y:=NegativetiveRootVectors(R);
[ v.37, v.38, v.39, v.40, v.41, v.42, v.43, v.44, v.45, v.46, v.47, v.48, 
  v.49, v.50, v.51, v.52, v.53, v.54, v.55, v.56, v.57, v.58, v.59, v.60, 
  v.61, v.62, v.63, v.64, v.65, v.66, v.67, v.68, v.69, v.70, v.71, v.72 ]
> CanonicalGenerators(R)[3]
[ v.73, v.74, v.75, v.76, v.77, v.78 ]
\end{verbatim}
Dans $E_{6}$, il y a deux orbites nilpotentes distingu{\'e}es non
r{\'e}guli{\`e}res\,: \\

\begin{enumerate} 
\item Caract{\'e}ristique\,: 
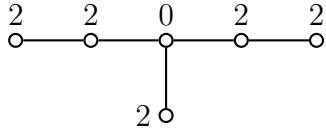
\begin{center}
\begin{pspicture}(-5,-1.2)(10,1)

\pscircle(-2,0){1mm}
\pscircle(-1,0){1mm}
\pscircle(0,0){1mm}
\pscircle(1,0){1mm}
\pscircle(2,0){1mm}
\pscircle(0,-1){1mm}

\psline(-1.9,0)(-1.1,0)
\psline(-0.1,0)(-0.9,0)
\psline(0.1,0)(0.9,0)
\psline(1.1,0)(1.9,0)
\psline(0,-0.1)(0,-0.9)

\rput[b](-2,0.2){$2$}
\rput[b](-1,0.2){$2$}
\rput[b](0,0.2){$0$}
\rput[b](1,0.2){$2$}
\rput[b](2,0.2){$2$}
\rput[r](-0.2,-1){$2$}
\end{pspicture}
\end{center}
On d{\'e}finit les {\'e}l{\'e}ments {\tt e} et {\tt f} du
$\mathfrak{sl}_{2}$-triplet correspondant dans les donn{\'e}es de \cite{Triplets}\,:
\begin{verbatim}
> e:=x[1]+x[2]+x[5]+x[6]+x[8]+x[9];
v.1+v.2+v.5+v.6+v.8+v.9
> f:=(12)*y[1]+(8)*y[2]+(-8)*y[3]+(22)*y[5]+(12)*y[6]+(8)*y[8]+
(22)*y[9]+(8)*y[10];;
\end{verbatim}
On v{\'e}rifie que le crochet {\tt e*f} est {\'e}gal {\`a} l'{\'e}l{\'e}ment neutre de
la caract{\'e}ristique et on pose {\tt h:=e*f}\,:
\begin{verbatim}
> e*f;
(12)*v.73+(16)*v.74+(22)*v.75+(30)*v.76+(22)*v.77+(12)*v.78
> h:=e*f;;
\end{verbatim}
On calcule le centralisateur {\tt g} de {\tt e} et on en donne une base
{\tt Bg}. On calcule ensuite le centre {\tt z} et on donne une base {\tt Bz} de {\tt z}\,:
\begin{verbatim}
> g:=LieCentralizer(L,Subspace(L,[e]));
<Lie algebra of dimension 8 over Rationals>
\end{verbatim}
Le centralisateur est de dimension 8 (ce que l'on savait d{\'e}j{\`a}).
\begin{verbatim}
> Bg:=BasisVectors(Basis(g));;
> z:=LieCentre(g);
<two-sided ideal in <Lie algebra of dimension 8 over Rationals>, 
(dimension 5)>
\end{verbatim}
Le centre est un id{\'e}al de dimension 5 dans {\tt g}.
\begin{verbatim}
> Bz:=BasisVectors(Basis(z));
[ v.1+v.2+v.5+v.6+v.8+v.9, v.23+(-1)*v.25+v.26, 
  v.27+(-1)*v.29+(-1)*v.30+(-1)*v.31, v.34+v.35, v.36 ]
\end{verbatim}
On calcule les \guillemotleft poids\guillemotright \  de {\tt z} en {\'e}valuant {\tt h*Bz[i]} pour
$i=1,\ldots,5$. On sait d{\'e}j{\`a} que {\tt h*Bz[1]=(2)*Bz[1]} car {\tt Bz[1]=e}.
\begin{verbatim}
> h*Bz[2];
(8)*v.23+(-8)*v.25+(8)*v.26
> h*Bz[3];
(10)*v.27+(-10)*v.29+(-10)*v.30+(-10)*v.31
> h*Bz[4];
(14)*v.34+(14)*v.35
> h*Bz[5];
(16)*v.36 
\end{verbatim}
On obtient que les poids sont $2,8,10,14,16$, d'o{\`u} $m_{r}=16$. Il y a
donc trois matrices {\`a} {\'e}tudier. 

\begin{enumerate}
\item $m_{i_{1}}=8$, $m_{k_{1}}=10$. On cherche les {\'e}l{\'e}ments de $V(10)$
parmi {\tt Bg} en calculant {\tt h*Bg[i]}, pour $i=1,\ldots,8$ et on effectue le calcul correspondant\,:
\begin{verbatim} 
> h*Bg[5];
(10)*v.27+(-10)*v.29
> h*Bg[6];
(10)*v.30+(10)*v.31
\end{verbatim}
Le sous-espace $V(10)$ 
est engendr{\'e} par les vecteurs {\tt Bg[5]} et {\tt Bg[6]}.
\begin{verbatim}
> ((f*Bg[5])*Bz[2]);
(-20)*v.36
> ((f*Bg[6])*Bz[2]);
(20)*v.36
\end{verbatim}
Notons que {\tt v.36} est bien l'{\'e}l{\'e}ment {\tt Bz[5]} correspondant au
plus haut poids 16. La matrice qu'il faut consid{\'e}rer est {\tt
  M(Bz[2])}; elle est donn{\'e}e par 
$\left[
\begin{array}{cc}
-20 & 20
\end{array}
\right]$. C'est une matrice de rang 1.\\
\\
\rem Lorsque $d_{r}=1$, il suffit de trouver \underline{un}
{\'e}l{\'e}ment de $V(m_{k_{1}})$ qui donne un crochet non nul; dans la suite
on donnera seulement le calcul correspondant {\`a} cet {\'e}l{\'e}ment.

\item $m_{i_{1}}=10$, $m_{k_{1}}=8$.
\begin{verbatim}
> ((f*Bz[2])*Bz[3]);
(-40)*v.36
\end{verbatim}

\item $m_{i_{1}}=14$, $m_{k_{1}}=4$.
\begin{verbatim}
> h*Bg[2];
(4)*v.7+(2)*v.11+(2)*v.12+(2)*v.13+(-2)*v.14+(2)*v.15+(-4)*v.16

> ((f*Bg[2])*Bz[4]);
(14)*v.36
\end{verbatim}
\end{enumerate}
Ces trois calculs montrent que les conditions de la partie 5 
sont v{\'e}rifi{\'e}es pour tout {\'e}l{\'e}ment non  nul de {\tt z}. Ainsi cette
orbite v{\'e}rifie $(P)$.\\

\item Caract{\'e}ristique\,: 
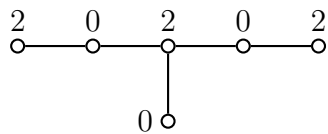
\begin{center}
\begin{pspicture}(-5,-1.2)(10,1)

\pscircle(-2,0){1mm}
\pscircle(-1,0){1mm}
\pscircle(0,0){1mm}
\pscircle(1,0){1mm}
\pscircle(2,0){1mm}
\pscircle(0,-1){1mm}

\psline(-1.9,0)(-1.1,0)
\psline(-0.1,0)(-0.9,0)
\psline(0.1,0)(0.9,0)
\psline(1.1,0)(1.9,0)
\psline(0,-0.1)(0,-0.9)

\rput[b](-2,0.2){$2$}
\rput[b](-1,0.2){$0$}
\rput[b](0,0.2){$2$}
\rput[b](1,0.2){$0$}
\rput[b](2,0.2){$2$}
\rput[r](-0.2,-1){$0$}
\end{pspicture}
\end{center}
D{\'e}finition du $\mathfrak{sl}_{2}$-triplet\,:
\begin{verbatim}
>e:=x[7]+x[8]+x[9]+x[10]=x[11]+x[19];
v.7+v.8+v.9+v.10+v.11+v.19
> f:=(8)*y[7]+(9)*y[8]+(5)*y[9]+(5)*y[10]+(8)*y[11]+y[19];
(8)*v.43+(9)*v.44+(5)*v.45+(5)*v.46+(8)*v.47+v.55
> e*f;                                                    
(8)*v.73+(10)*v.74+(14)*v.75+(20)*v.76+(14)*v.77+(8)*v.78
> h:=e*f;;
\end{verbatim}
Calcul de {\tt g}, {\tt Bg}, {\tt z} et {\tt Bz}\,:
\begin{verbatim}
> g:=LieCentralizer(L,Subspace(L,[e]));
<Lie algebra of dimension 12 over Rationals>
> z:=LieCentre(g);
<two-sided ideal in <Lie algebra of dimension 12 over Rationals>, 
  (dimension 4)>
> Bg:=BasisVectors(Basis(g));;
> Bz:=BasisVectors(Basis(z));
[ v.7+v.8+v.9+v.10+v.11+v.19, v.32+(-1)*v.33, v.35, v.36 ]
> h*Bz[2];
(8)*v.32+(-8)*v.33
> h*Bz[3];
(10)*v.35
> h*Bz[4];
(10)*v.36
\end{verbatim}  
Les poids de {\tt z} sont 2,8,10,10; d'o{\`u} $m_ {r}=10$. Il y a deux matrices {\`a} {\'e}tudier.

\begin{enumerate}
\item $m_{i_{1}}=8$,  $m_{k_{1}}=4$
\begin{verbatim}
> h*Bg[4];
(4)*v.17+(-4)*v.18+(4)*v.20+(-4)*v.21
> h*Bg[5];
(4)*v.12+(4)*v.16+(-8)*v.22+(4)*v.24
> h*Bg[6];
(4)*v.22+(-4)*v.24+(-4)*v.25

> ((f*Bg[4])*Bz[2]);        
(-16)*v.35
\end{verbatim}

\item $m_{i_{1}}=10$,  $m_{k_{1}}=2$
\begin{verbatim}
> h*Bg[2];
(2)*v.1+(2)*v.4+(2)*v.6+(2)*v.13+(2)*v.14+(-6)*v.15
> h*Bg[3];
(2)*v.19

> ((f*Bg[1])*Bz[3]); ((f*Bg[2])*Bz[3]); ((f*Bg[3])*Bz[3]);
(-10)*v.35
(10)*v.36
0*v.1
> ((f*Bg[1])*Bz[4]); ((f*Bg[2])*Bz[4]); ((f*Bg[3])*Bz[4]);
(-9)*v.36
(6)*v.35
(-1)*v.36
\end{verbatim}
La matrice {\`a} {\'e}tudier est
$$\alpha 
\left[
\begin{array}{ccc}
-10 & 0 & 0 \\
0 & 10 & 0 
\end{array}
\right] + \beta
\left[
\begin{array}{ccc}
0 & 6 & 0 \\
-9 & 0 & -1 
\end{array}
\right]=
\left[
\begin{array}{ccc}
-10\alpha & -\beta & 0 \\
-9 \beta & 10\alpha & -\beta 
\end{array}
\right]
 \cdot$$
On v{\'e}rifie que c'est une matrice de rang 2 pour tout
couple $(\alpha,\beta)$ non nul. 
\end{enumerate}
La propri{\'e}t{\'e} $(P)$ est v{\'e}rifi{\'e}e pour cette orbite.\\
\end{enumerate}
{\bf Conclusion\,:} $E_{6}$ v{\'e}rifie $(P)$.

\subsection*{Cas de $E_{7}$.}
\noindent D{\'e}finition de {\tt L}\,:
\begin{verbatim}
> L:=SimpleLieAlgebra("E",7,Rationals);
<Lie algebra of dimension 133 over Rationals>
> R:=RootSystem(L);
<root system of rank 7>
> P:=PositiveRoots(R);; 
> x:=PositiveRootVectors(R);
[ v.1, v.2, v.3, v.4, v.5, v.6, v.7, v.8, v.9, v.10, v.11, v.12, v.13, v.14, 
  v.15, v.16, v.17, v.18, v.19, v.20, v.21, v.22, v.23, v.24, v.25, v.26, 
  v.27, v.28, v.29, v.30, v.31, v.32, v.33, v.34, v.35, v.36, v.37, v.38, 
  v.39, v.40, v.41, v.42, v.43, v.44, v.45, v.46, v.47, v.48, v.49, v.50, 
  v.51, v.52, v.53, v.54, v.55, v.56, v.57, v.58, v.59, v.60, v.61, v.62, 
  v.63 ]
> y:=NegativeRootVectors(R);
[ v.64, v.65, v.66, v.67, v.68, v.69, v.70, v.71, v.72, v.73, v.74, v.75, 
  v.76, v.77, v.78, v.79, v.80, v.81, v.82, v.83, v.84, v.85, v.86, v.87, 
  v.88, v.89, v.90, v.91, v.92, v.93, v.94, v.95, v.96, v.97, v.98, v.99, 
  v.100, v.101, v.102, v.103, v.104, v.105, v.106, v.107, v.108, v.109, 
  v.110, v.111, v.112, v.113, v.114, v.115, v.116, v.117, v.118, v.119, 
  v.120, v.121, v.122, v.123, v.124, v.125, v.126 ]
> CanonicalGenerators(R)[3];
[ v.127, v.128, v.129, v.130, v.131, v.132, v.133 ]
\end{verbatim}
Dans $E_{7}$, il y a cinq orbites nilpotentes distingu{\'e}es non
r{\'e}guli{\`e}res\,:\\

\begin{enumerate}
\item Caract{\'e}ristique\,:
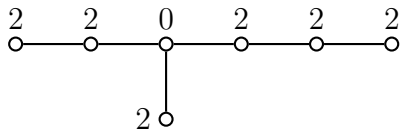
\begin{center}
\begin{pspicture}(-5,-1.2)(10,1)

\pscircle(-2,0){1mm}
\pscircle(-1,0){1mm}
\pscircle(0,0){1mm}
\pscircle(1,0){1mm}
\pscircle(2,0){1mm}
\pscircle(3,0){1mm}

\pscircle(0,-1){1mm}

\psline(-1.9,0)(-1.1,0)
\psline(-0.1,0)(-0.9,0)
\psline(0.1,0)(0.9,0)
\psline(1.1,0)(1.9,0)
\psline(2.1,0)(2.9,0)
\psline(0,-0.1)(0,-0.9)

\rput[b](-2,0.2){$2$}
\rput[b](-1,0.2){$2$}
\rput[b](0,0.2){$0$}
\rput[b](1,0.2){$2$}
\rput[b](2,0.2){$2$}
\rput[b](3,0.2){$2$}
\rput[r](-0.2,-1){$2$}
\end{pspicture}
\end{center} 
D{\'e}finition du $\mathfrak{sl}_{2}$-triplet\,:
\begin{verbatim}
> e:=x[1]+x[2]+x[3]+x[6]+x[7]+x[9]+x[11];
v.1+v.2+v.3+v.6+v.7+v.9+v.11
>f:=(26)*y[1]+(22)*y[2]+(50)*y[3]+(22)*y[5]+(40)*y[6]
(21)*y[7]+(15)*y[9]+(-15)*y[10]+(57)*y[11];
(26)*v.64+(22)*v.65+(50)*v.66+(22)*v.68+(40)*v.69
+(21)*v.70+(15)*v.72+(-15)*v.73+(57)*v.74
> e*f;
(26)*v.127+(37)*v.128+(50)*v.129+(72)*v.130+(57)*v.131+(40)*v.132+(21)*v.133
> h:=e*f;;
\end{verbatim}
Calcul de {\tt g}, {\tt Bg}, {\tt z} et {\tt Bz}\,:
\begin{verbatim}
> g:=LieCentralizer(L,Subspace(L,[e]));
<Lie algebra of dimension 9 over Rationals>
> z:=LieCentre(g);
<two-sided ideal in <Lie algebra of dimension 9 over Rationals>, 
(dimension 6)>
> Bg:=BasisVectors(Basis(g));;
> Bz:=BasisVectors(Basis(z));
[ v.1+v.2+v.3+v.6+v.7+v.9+v.11, v.33+(2)*v.34+(-1)*v.36+(-1)*v.37
+(-1)*v.38+v.40+(-3)*v.41, v.46+v.47+v.48+v.49+v.50, v.56+(-1)*v.57, v.60, 
  v.63 ]
> h*Bz[2];
(10)*v.33+(20)*v.34+(-10)*v.36+(-10)*v.37+(-10)*v.38+(10)*v.40+(-30)*v.41
> h*Bz[3];
(14)*v.46+(14)*v.47+(14)*v.48+(14)*v.49+(14)*v.50
> h*Bz[4];
(18)*v.56+(-18)*v.57
> h*Bz[5];
(22)*v.60
> h*Bz[6];
(26)*v.63
\end{verbatim}  
Les poids de {\tt z} sont 2,10,14,18,22,26; d'o{\`u} $m_{r}=26$. Il y a quatre matrices {\`a} {\'e}tudier.

\begin{enumerate}
\item $m_{i_{1}}=10$,  $m_{k_{1}}=18$
\begin{verbatim}
> ((f*Bz[4])*Bz[2]);
(90)*v.63
\end{verbatim}

\item $m_{i_{1}}=14$,  $m_{k_{1}}=14$
\begin{verbatim}
> ((f*Bz[3])*Bz[3]);
(-98)*v.63
\end{verbatim}

\item $m_{i_{1}}=10$,  $m_{k_{1}}=18$
\begin{verbatim}
> ((f*Bz[4])*Bz[2]);
(90)*v.63
\end{verbatim}

\item $m_{i_{1}}=18$,  $m_{k_{1}}=10$
\begin{verbatim}
> h*Bg[3];
(10)*v.33+(20)*v.34+(-10)*v.36+(20)*v.37+(20)*v.38+(10)*v.40
> h*Bg[4];
(10)*v.37+(10)*v.38+(10)*v.41

> ((f*Bg[4])*Bz[4]);
(-30)*v.63
\end{verbatim}
4) $m_{i_{1}}=22$,  $m_{k_{1}}=6$
\begin{verbatim}
> h*Bg[2];
(6)*v.19+(4)*v.20+(-2)*v.21+(2)*v.22+(4)*v.23+(-2)*v.24+(-2)*v.25+(-2)*v.28

> ((f*Bg[2])*Bz[5]);
(-22)*v.63
\end{verbatim}
\end{enumerate}
La propri{\'e}t{\'e} $(P)$ est v{\'e}rifi{\'e}e pour cette orbite.\\

\item Caract{\'e}ristique\,:
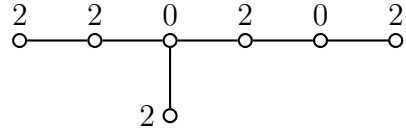
\begin{center}
\begin{pspicture}(-5,-1.2)(10,1)

\pscircle(-2,0){1mm}
\pscircle(-1,0){1mm}
\pscircle(0,0){1mm}
\pscircle(1,0){1mm}
\pscircle(2,0){1mm}
\pscircle(3,0){1mm}

\pscircle(0,-1){1mm}

\psline(-1.9,0)(-1.1,0)
\psline(-0.1,0)(-0.9,0)
\psline(0.1,0)(0.9,0)
\psline(1.1,0)(1.9,0)
\psline(2.1,0)(2.9,0)
\psline(0,-0.1)(0,-0.9)

\rput[b](-2,0.2){$2$}
\rput[b](-1,0.2){$2$}
\rput[b](0,0.2){$0$}
\rput[b](1,0.2){$2$}
\rput[b](2,0.2){$0$}
\rput[b](3,0.2){$2$}
\rput[r](-0.2,-1){$2$}
\end{pspicture}
\end{center} 
D{\'e}finition du $\mathfrak{sl}_{2}$-triplet\,:
\begin{verbatim}
> e:=x[1]+x[2]+x[3]+x[5]+x[7]+x[9]+x[18];                               
v.1+v.2+v.3+v.5+v.7+v.9+v.18
>f:=(22)*y[1]+(3)*y[2]+(42)*y[3]+(15)*y[5]+(17)*y[7]
+(28)*y[9]+(-28)*y[10]+(3)*y[12]+(3)*y[13]+(32)*y[18];
(22)*v.64+(3)*v.65+(42)*v.66+(15)*v.68+(17)*v.70
+(28)*v.72+(-28)*v.73+(3)*v.75+(3)*v.76+(32)*v.81
> e*f;
(22)*v.127+(31)*v.128+(42)*v.129+(60)*v.130+(47)*v.131+(32)*v.132+(17)*v.133
> h:=e*f;;
\end{verbatim}
Calcul de {\tt g}, {\tt Bg}, {\tt z} et {\tt Bz}\,:
\begin{verbatim}
> g:=LieCentralizer(L,Subspace(L,[e]));
<Lie algebra of dimension 11 over Rationals>
> z:=LieCentre(g);                     
<two-sided ideal in <Lie algebra of dimension 11 over Rationals>, 
  (dimension 5)>
> Bg:=BasisVectors(Basis(g));;
> Bz:=BasisVectors(Basis(z));
[ v.1+v.2+v.3+v.5+v.7+v.9+v.18, v.39+(-1)*v.42+(-1)*v.43+(-1)*v.44
+(-2)*v.45+v.49, v.51+v.53+v.55+v.57, v.60, v.63 ]
> h*Bz[2];
(10)*v.39+(-10)*v.42+(-10)*v.43+(-10)*v.44+(-20)*v.45+(10)*v.49
> h*Bz[3];
(14)*v.51+(14)*v.53+(14)*v.55+(14)*v.57
> h*Bz[4];
(18)*v.60
> h*Bz[5];
(22)*v.63
\end{verbatim}  
Les poids de {\tt z} sont 2,10,14,18,22; d'o{\`u} $m_{r}=22$. Il y a trois matrices {\`a} {\'e}tudier.

\begin{enumerate}
\item $m_{i_{1}}=10$,  $m_{k_{1}}=14$
\begin{verbatim}
> h*Bg[7];
(14)*v.53+(-14)*v.54
> h*Bg[8];
(14)*v.51+(14)*v.54+(14)*v.55+(14)*v.57

> ((f*Bg[7])*Bz[2]);
(-70)*v.63
\end{verbatim}

\item $m_{i_{1}}=14$,  $m_{k_{1}}=10$
\begin{verbatim}
> h*Bg[5];
(10)*v.42+(10)*v.43+(10)*v.44+(20)*v.45
> h*Bg[6];
(10)*v.39+(10)*v.49

> ((f*Bg[5])*Bz[3]);
(70)*v.63
\end{verbatim}

\item $m_{i_{1}}=18$,  $m_{k_{1}}=4$
\begin{verbatim}
> h*Bg[3];
(6)*v.20+(9)*v.21+(-3)*v.22+(-3)*v.27+(-9)*v.28+
(3)*v.29+(6)*v.30+(-3)*v.31+(-3)*v.35
> ((f*Bg[3])*Bz[4]);
(18)*v.63
\end{verbatim}
\end{enumerate}
La propri{\'e}t{\'e} $(P)$ est v{\'e}rifi{\'e}e pour cette orbite.\\

\item Caract{\'e}ristique\,:
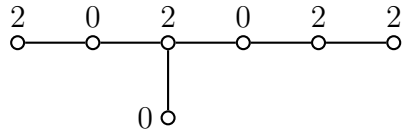
\begin{center}
\begin{pspicture}(-5,-1.2)(10,1)

\pscircle(-2,0){1mm}
\pscircle(-1,0){1mm}
\pscircle(0,0){1mm}
\pscircle(1,0){1mm}
\pscircle(2,0){1mm}
\pscircle(3,0){1mm}

\pscircle(0,-1){1mm}

\psline(-1.9,0)(-1.1,0)
\psline(-0.1,0)(-0.9,0)
\psline(0.1,0)(0.9,0)
\psline(1.1,0)(1.9,0)
\psline(2.1,0)(2.9,0)
\psline(0,-0.1)(0,-0.9)

\rput[b](-2,0.2){$2$}
\rput[b](-1,0.2){$0$}
\rput[b](0,0.2){$2$}
\rput[b](1,0.2){$0$}
\rput[b](2,0.2){$2$}
\rput[b](3,0.2){$2$}
\rput[r](-0.2,-1){$0$}
\end{pspicture}
\end{center} 
D{\'e}finition du $\mathfrak{sl}_{2}$-triplet\,:
\begin{verbatim}
> e:=x[7]+x[8]+x[9]+x[10]+x[11]+x[12]+x[22];                            
v.7+v.8+v.9+v.10+v.11+v.12+v.22
> f:=(15)*y[7]+(18)*y[8]+(24)*y[9]+(15)*y[10]+(10)*y[11]+
(28)*y[12]+y[22];
(15)*v.70+(18)*v.71+(24)*v.72+(15)*v.73+(10)*v.74+(28)*v.75+v.85
> e*f;
(18)*v.127+(25)*v.128+(34)*v.129+(50)*v.130+(39)*v.131+(28)*v.132+(15)*v.133
> h:=e*f;;
\end{verbatim}
Calcul de {\tt g}, {\tt Bg}, {\tt z} et {\tt Bz}\,:
\begin{verbatim}
> g:=LieCentralizer(L,Subspace(L,[e]));
<Lie algebra of dimension 13 over Rationals>
> z:=LieCentre(g);
<two-sided ideal in <Lie algebra of dimension 13 over Rationals>, 
  (dimension 5)>
> Bz:=BasisVectors(Basis(z));
[ v.7+v.8+v.9+v.10+v.11+v.12+v.22, v.47+(-3)*v.48+(-1)*v.49+(-2)*v.50, 
  v.58+v.59, v.62, v.63 ]
> h*Bz[2];
(10)*v.47+(-30)*v.48+(-10)*v.49+(-20)*v.50
> h*Bz[3];
(14)*v.58+(14)*v.59
> h*Bz[4];
(16)*v.62
> h*Bz[5];
(18)*v.63
\end{verbatim}  
Les poids de {\tt z} sont 2,10,14,16,18; d'o{\`u} $m_{r}=18$. Il y a trois matrices {\`a}
{\'e}tudier.

\begin{enumerate}
\item $m_{i_{1}}=10$,  $m_{k_{1}}=10$
\begin{verbatim}
> ((f*Bz[2])*Bz[2]);
(-150)*v.63
\end{verbatim}

\item $m_{i_{1}}=14$,  $m_{k_{1}}=6$
\begin{verbatim}
> h*Bg[4];
(6)*v.30+(-12)*v.31+(6)*v.32+(6)*v.33+(-6)*v.35
> h*Bg[5];          
(6)*v.36+(-6)*v.37+(-6)*v.40

> ((f*Bg[4])*Bz[3]); 
(42)*v.63
\end{verbatim}

\item $m_{i_{1}}=16$,  $m_{k_{1}}=4$
\begin{verbatim}
> h*Bg[3];
(4)*v.13+(4)*v.14+(4)*v.18+(-12)*v.26+(8)*v.28+(4)*v.29

> ((f*Bg[3])*Bz[4]);
(12)*v.63
\end{verbatim}
\end{enumerate}
La propri{\'e}t{\'e} $(P)$ est v{\'e}rifi{\'e}e pour cette orbite.\\

\item Caract{\'e}ristique\,:
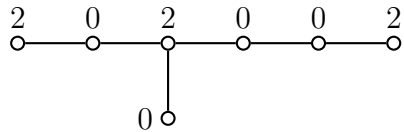
\begin{center}
\begin{pspicture}(-5,-1.2)(10,1)

\pscircle(-2,0){1mm}
\pscircle(-1,0){1mm}
\pscircle(0,0){1mm}
\pscircle(1,0){1mm}
\pscircle(2,0){1mm}
\pscircle(3,0){1mm}

\pscircle(0,-1){1mm}

\psline(-1.9,0)(-1.1,0)
\psline(-0.1,0)(-0.9,0)
\psline(0.1,0)(0.9,0)
\psline(1.1,0)(1.9,0)
\psline(2.1,0)(2.9,0)
\psline(0,-0.1)(0,-0.9)

\rput[b](-2,0.2){$2$}
\rput[b](-1,0.2){$0$}
\rput[b](0,0.2){$2$}
\rput[b](1,0.2){$0$}
\rput[b](2,0.2){$0$}
\rput[b](3,0.2){$2$}
\rput[r](-0.2,-1){$0$}
\end{pspicture}
\end{center} 
D{\'e}finition du $\mathfrak{sl}_{2}$-triplet\,:
\begin{verbatim}
> e:=x[8]+x[9]+x[13]+x[16]+x[17]+x[18]+x[29];                           
v.8+v.9+v.13+v.16+v.17+v.18+v.29
> f:=(14)*y[8]+(9)*y[9]+(-9)*y[10]+(11)*y[13]+(9)*y[16]
+(11)*y[17]+(8)*y[18]+(9)*y[19]+y[29];
(14)*v.71+(9)*v.72+(-9)*v.73+(11)*v.76+(9)*v.79
+(11)*v.80+(8)*v.81+(9)*v.82+v.92
> e*f;
(14)*v.127+(19)*v.128+(26)*v.129+(38)*v.130+(29)*v.131+(20)*v.132+(11)*v.133
> h:=e*f;;
\end{verbatim}
Calcul de {\tt g}, {\tt Bg}, {\tt z} et {\tt Bz}\,:
\begin{verbatim}
> g:=LieCentralizer(L,Subspace(L,[e]));
<Lie algebra of dimension 17 over Rationals>
> z:=LieCentre(g);
<two-sided ideal in <Lie algebra of dimension 17 over Rationals>, 
  (dimension 3)>
> Bz:=BasisVectors(Basis(z)); 
[ v.8+v.9+v.13+v.16+v.17+v.18+v.29, v.56+(3)*v.57+(2)*v.59, v.63 ]
> h*Bz[2];
(10)*v.56+(30)*v.57+(20)*v.59
> h*Bz[3];
(14)*v.63
\end{verbatim}  
Les poids de {\tt z} sont 2,10,14; d'o{\`u} $m_{r}=14$. Il n'y a qu'une matrice {\`a} {\'e}tudier.

\begin{enumerate}
\item $m_{i_{1}}=10$,  $m_{k_{1}}=6$
\begin{verbatim}
> h*Bg[7];
(6)*v.37+(-6)*v.38+(-6)*v.41
> h*Bg[8];
(6)*v.39+(6)*v.43+(-6)*v.45
> h*Bg[9];
(6)*v.42+(6)*v.46+(6)*v.49

> ((f*Bg[7])*Bz[2]);
(30)*v.63
\end{verbatim}
\end{enumerate}
La propri{\'e}t{\'e} $(P)$ est v{\'e}rifi{\'e}e pour cette orbite.\\

\item Caract{\'e}ristique\,:
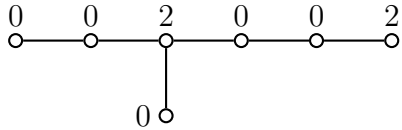
\begin{center}
\begin{pspicture}(-5,-1.2)(10,1)

\pscircle(-2,0){1mm}
\pscircle(-1,0){1mm}
\pscircle(0,0){1mm}
\pscircle(1,0){1mm}
\pscircle(2,0){1mm}
\pscircle(3,0){1mm}

\pscircle(0,-1){1mm}

\psline(-1.9,0)(-1.1,0)
\psline(-0.1,0)(-0.9,0)
\psline(0.1,0)(0.9,0)
\psline(1.1,0)(1.9,0)
\psline(2.1,0)(2.9,0)
\psline(0,-0.1)(0,-0.9)

\rput[b](-2,0.2){$0$}
\rput[b](-1,0.2){$0$}
\rput[b](0,0.2){$2$}
\rput[b](1,0.2){$0$}
\rput[b](2,0.2){$0$}
\rput[b](3,0.2){$2$}
\rput[r](-0.2,-1){$0$}
\end{pspicture}
\end{center} 
D{\'e}finition du $\mathfrak{sl}_{2}$-triplet\,:
\begin{verbatim}
> e:=x[13]+x[14]+x[15]+x[16]+x[17]+x[18]+x[33];                         
v.13+v.14+v.15+v.16+v.17+v.18+v.33
>f:=(9)*y[13]+(5)*y[14]+(2)*y[15]+(8)*y[16]+(8)*y[17]
+(2)*y[18]+(5)*y[33]; 
(9)*v.76+(5)*v.77+(2)*v.78+(8)*v.79+(8)*v.80+(2)*v.81+(5)*v.96
> e*f;
(10)*v.127+(15)*v.128+(20)*v.129+(30)*v.130+(23)*v.131+(16)*v.132+(9)*v.133
> h:=e*f;;
\end{verbatim}
Calcul de {\tt g}, {\tt Bg}, {\tt z} et {\tt Bz}\,:
\begin{verbatim}
> g:=LieCentralizer(L,Subspace(L,[e]));
<Lie algebra of dimension 21 over Rationals>
> Bg:=BasisVectors(Basis(g));;
> z:=LieCentre(g);
<two-sided ideal in <Lie algebra of dimension 21 over Rationals>, 
  (dimension 4)>
> Bz:=BasisVectors(Basis(z));
[ v.13+v.14+v.15+v.16+v.17+v.18+v.33, v.61, v.62, v.63 ]
> h*Bz[2];
(10)*v.61
> h*Bz[3];
(10)*v.62
> h*Bz[4];
(10)*v.63
\end{verbatim}  
Les poids de {\tt z} sont 2,10,10,10; d'o{\`u} $m_{r}=10$. Il n'y a qu'une matrice {\`a}
{\'e}tudier.

\begin{enumerate}
\item $m_{i_{1}}=10$,  $m_{k_{1}}=2$
\begin{verbatim}
> h*Bg[1];
(2)*v.15+(2)*v.18
> h*Bg[2]; 
(2)*v.4+(2/3)*v.19+(-4/3)*v.20+(2/3)*v.21+(2/3)*v.22+(-2/3)*v.23
> h*Bg[3];
(2)*v.19+(-4)*v.20+(2)*v.21+(2)*v.22+(4)*v.23+(-6)*v.24
> h*Bg[4];
(2)*v.7+(-6)*v.9+(4)*v.10+(2)*v.11+(2)*v.26+(-4)*v.27
> h*Bg[5];
(2)*v.9+(-2)*v.10+(2)*v.29
> h*Bg[6];
(2)*v.13+(2)*v.14+(2)*v.16+(2)*v.17+(2)*v.33

> ((f*Bg[1])*Bz[2]);((f*Bg[2])*Bz[2]);((f*Bg[3])*Bz[2]);
((f*Bg[4])*Bz[2]);((f*Bg[5])*Bz[2]);((f*Bg[6])*Bz[2]);
(-2)*v.61
(8/3)*v.62
(2)*v.62
(-4)*v.63
0*v.1
(-8)*v.61
>  ((f*Bg[1])*Bz[3]); ((f*Bg[2])*Bz[3]); ((f*Bg[3])*Bz[3]); 
((f*Bg[4])*Bz[3]);((f*Bg[5])*Bz[3]);((f*Bg[6])*Bz[3]);
(-2)*v.62
(4/3)*v.63
(4)*v.63
(2)*v.61
(2)*v.61
(-8)*v.62
> ((f*Bg[1])*Bz[4]); ((f*Bg[2])*Bz[4]); ((f*Bg[3])*Bz[4]); 
((f*Bg[4])*Bz[4]);((f*Bg[5])*Bz[4]);((f*Bg[6])*Bz[4]);
0*v.1
(-10/3)*v.61
(-10)*v.61
(10)*v.62
0*v.1
(-10)*v.63
\end{verbatim}
La matrice {\`a} {\'e}tudier est
$$\left[
\begin{array}{cccccc}
-2\alpha & -10/3 \gamma & -10 \gamma & 2 \beta & 2 \beta &-8 \alpha\\
-2\beta & -8/3 \alpha & 2\alpha & 10\gamma & 0 & -8\beta\\
0 & 4/3 \beta & 4\beta & -4\alpha & 0 & -10\gamma
\end{array}
\right]
 \cdot$$
Une {\'e}tude {\'e}l{\'e}mentaire montre que cette matrice est de rang 3 pour tout
 3-uplet $(\alpha,\beta,\gamma)$ non nul. 
\end{enumerate}
La propri{\'e}t{\'e} $(P)$ est v{\'e}rifi{\'e}e pour cette orbite.\\

\end{enumerate}
{\bf Conclusion\,:} $E_{7}$ v{\'e}rifie $(P)$.

\subsection*{Cas de $E_{8}$.}
\noindent D{\'e}finition de {\tt L}\,:
\begin{verbatim}
> L:=SimpleLieAlgebra("E",8,Rationals);
<Lie algebra of dimension 248 over Rationals>
> R:=RootSystem(L);
<root system of rank 7>
> P:=PositiveRoots(R);; 
> x:=PositiveRootVectors(R);
[ v.1, v.2, v.3, v.4, v.5, v.6, v.7, v.8, v.9, v.10, v.11, v.12, v.13, v.14, 
  v.15, v.16, v.17, v.18, v.19, v.20, v.21, v.22, v.23, v.24, v.25, v.26, 
  v.27, v.28, v.29, v.30, v.31, v.32, v.33, v.34, v.35, v.36, v.37, v.38, 
  v.39, v.40, v.41, v.42, v.43, v.44, v.45, v.46, v.47, v.48, v.49, v.50, 
  v.51, v.52, v.53, v.54, v.55, v.56, v.57, v.58, v.59, v.60, v.61, v.62, 
  v.63, v.64, v.65, v.66, v.67, v.68, v.69, v.70, v.71, v.72, v.73, v.74, 
  v.75, v.76, v.77, v.78, v.79, v.80, v.81, v.82, v.83, v.84, v.85, v.86, 
  v.87, v.88, v.89, v.90, v.91, v.92, v.93, v.94, v.95, v.96, v.97, v.98, 
  v.99, v.100, v.101, v.102, v.103, v.104, v.105, v.106, v.107, v.108, 
  v.109, v.110, v.111, v.112, v.113, v.114, v.115, v.116, v.117, v.118, 
  v.119, v.120 ]
> y:=NegativeRootVectors(R);
[ v.121, v.122, v.123, v.124, v.125, v.126, v.127, v.128, v.129, v.130, 
  v.131, v.132, v.133, v.134, v.135, v.136, v.137, v.138, v.139, v.140, 
  v.141, v.142, v.143, v.144, v.145, v.146, v.147, v.148, v.149, v.150, 
  v.151, v.152, v.153, v.154, v.155, v.156, v.157, v.158, v.159, v.160, 
  v.161, v.162, v.163, v.164, v.165, v.166, v.167, v.168, v.169, v.170, 
  v.171, v.172, v.173, v.174, v.175, v.176, v.177, v.178, v.179, v.180, 
  v.181, v.182, v.183, v.184, v.185, v.186, v.187, v.188, v.189, v.190, 
  v.191, v.192, v.193, v.194, v.195, v.196, v.197, v.198, v.199, v.200, 
  v.201, v.202, v.203, v.204, v.205, v.206, v.207, v.208, v.209, v.210, 
  v.211, v.212, v.213, v.214, v.215, v.216, v.217, v.218, v.219, v.220, 
  v.221, v.222, v.223, v.224, v.225, v.226, v.227, v.228, v.229, v.230, 
  v.231, v.232, v.233, v.234, v.235, v.236, v.237, v.238, v.239, v.240 ]
> CanonicalGenerators(R)[3];
[ v.241, v.242, v.243, v.244, v.245, v.246, v.247, v.248 ]
\end{verbatim}
Dans $E_{8}$, il y a dix orbites nilpotentes distingu{\'e}es non
r{\'e}guli{\`e}res\,:\\

\begin{enumerate}
\item Caract{\'e}ristique\,:
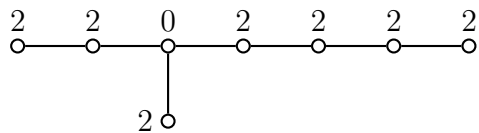
\begin{center}
\begin{pspicture}(-5,-1.2)(10,1)

\pscircle(-2,0){1mm}
\pscircle(-1,0){1mm}
\pscircle(0,0){1mm}
\pscircle(1,0){1mm}
\pscircle(2,0){1mm}
\pscircle(3,0){1mm}
\pscircle(4,0){1mm}

\pscircle(0,-1){1mm}

\psline(-1.9,0)(-1.1,0)
\psline(-0.1,0)(-0.9,0)
\psline(0.1,0)(0.9,0)
\psline(1.1,0)(1.9,0)
\psline(2.1,0)(2.9,0)
\psline(3.1,0)(3.9,0)
\psline(0,-0.1)(0,-0.9)

\rput[b](-2,0.2){$2$}
\rput[b](-1,0.2){$2$}
\rput[b](0,0.2){$0$}
\rput[b](1,0.2){$2$}
\rput[b](2,0.2){$2$}
\rput[b](3,0.2){$2$}
\rput[b](4,0.2){$2$}
\rput[r](-0.2,-1){$2$}
\end{pspicture}
\end{center}
D{\'e}finition du $\mathfrak{sl}_{2}$-triplet\,:
\begin{verbatim}
> e:=x[1]+x[2]+x[3]+x[6]+x[7]+x[8]+x[10]+x[12]; 
v.1+v.2+v.3+v.6+v.7+v.8+v.10+v.12
> f:=(72)*y[1]+(60)*y[2]+(142)*y[3]+(68)*y[5]+(132)*y[6]
+(90)*y[7]+(46)*y[8]+(38)*y[10]+(-38)*y[11]+(172)*y[12];                       
(72)*v.121+(68)*v.122+(142)*v.123+(68)*v.125+(132)*v.126
+(90)*v.127+(46)*v.128+(38)*v.130+(-38)*v.131+(172)*v.132
> e*f;
(72)*v.241+(106)*v.242+(142)*v.243+(210)*v.244+(172)*v.245+(132)*v.246
+(90)*v.247+(46)*v.248
> h:=e*f;;
\end{verbatim}
Calcul de {\tt g}, {\tt Bg}, {\tt z} et {\tt Bz}\,:
\begin{verbatim}
> g:=LieCentralizer(L,Subspace(L,[e]));
<Lie algebra of dimension 10 over Rationals>
> Bg:=BasisVectors(Basis(g));;
> z:=LieCentre(g);
<two-sided ideal in <Lie algebra of dimension 10 over Rationals>, 
  (dimension 7)>
> Bz:=BasisVectors(Basis(z));
[ v.1+v.2+v.3+v.6+v.7+v.8+v.10+v.12, v.54+(-1/2)*v.57+(-1/2)*v.58
+(-1/2)*v.59+(-1/2)*v.60+(-1/2)*v.61+(1/2)*v.62+(-1/2)*v.63, 
  v.84+v.85+(-1)*v.86+(-1)*v.87+(2)*v.88, v.95+(-1)*v.96+v.97+(-1)*v.98+v.99,
  v.109+(-1)*v.112, v.116, v.120 ]
> h*Bz[2];
(14)*v.54+(-7)*v.57+(-7)*v.58+(-7)*v.59+(-7)*v.60+(-7)*v.61
+(7)*v.62+(-7)*v.63
> h*Bz[3];
(22)*v.84+(22)*v.85+(-22)*v.86+(-22)*v.87+(44)*v.88
> h*Bz[4];
(26)*v.95+(-26)*v.96+(26)*v.97+(-26)*v.98+(26)*v.99
> h*Bz[5];
(34)*v.109+(-34)*v.112
> h*Bz[6];
(38)*v.116
> h*Bz[7];
(46)*v.120
\end{verbatim}  
Les poids de {\tt z} sont 2,14,22,26,34,38,46; d'o{\`u} $m_{r}=46$. Il y a cinq matrices {\`a}
{\'e}tudier.

\begin{enumerate}
\item $m_{i_{1}}=14$,  $m_{k_{1}}=34$
\begin{verbatim}
> ((f*Bz[5])*Bz[2]);
(-119)*v.120
\end{verbatim}

\item $m_{i_{1}}=22$,  $m_{k_{1}}=26$
\begin{verbatim}
> ((f*Bz[4])*Bz[3]);
(286)*v.120
\end{verbatim}

\item $m_{i_{1}}=26$,  $m_{k_{1}}=22$
\begin{verbatim}
> ((f*Bz[3])*Bz[4]);
(286)*v.120
\end{verbatim}

\item $m_{i_{1}}=34$,  $m_{k_{1}}=14$
\begin{verbatim}
> ((f*Bz[2])*Bz[5]);
(-119)*v.120
\end{verbatim}

\item $m_{i_{1}}=38$,  $m_{k_{1}}=10$
\begin{verbatim}
> h*Bg[2];
(10)*v.38+(20)*v.39+(-10)*v.41+(-30)*v.42+(20)*v.43+(40)*v.44+(40)*v.45
+(10)*v.48+(20)*v.49

> ((f*Bg[2])*Bz[6]);               
(190)*v.120
\end{verbatim}
\end{enumerate}
La propri{\'e}t{\'e} $(P)$ est v{\'e}rifi{\'e}e pour cette orbite.\\

\item Caract{\'e}ristique\,:
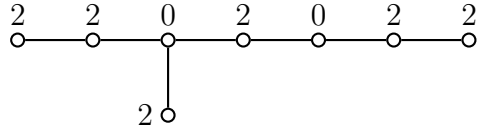
\begin{center}
\begin{pspicture}(-5,-1.2)(10,1)

\pscircle(-2,0){1mm}
\pscircle(-1,0){1mm}
\pscircle(0,0){1mm}
\pscircle(1,0){1mm}
\pscircle(2,0){1mm}
\pscircle(3,0){1mm}
\pscircle(4,0){1mm}

\pscircle(0,-1){1mm}

\psline(-1.9,0)(-1.1,0)
\psline(-0.1,0)(-0.9,0)
\psline(0.1,0)(0.9,0)
\psline(1.1,0)(1.9,0)
\psline(2.1,0)(2.9,0)
\psline(3.1,0)(3.9,0)
\psline(0,-0.1)(0,-0.9)

\rput[b](-2,0.2){$2$}
\rput[b](-1,0.2){$2$}
\rput[b](0,0.2){$0$}
\rput[b](1,0.2){$2$}
\rput[b](2,0.2){$0$}
\rput[b](3,0.2){$2$}
\rput[b](4,0.2){$2$}
\rput[r](-0.2,-1){$2$}
\end{pspicture}
\end{center}
D{\'e}finition du $\mathfrak{sl}_{2}$-triplet\,:
\begin{verbatim}
> e:=x[1]+x[2]+x[3]+x[5]+x[7]+x[8]+x[10]+x[20];                         
v.1+v.2+v.3+v.5+v.7+v.8+v.10+v.20
> f:=(61)*y[1]+(22)*y[2]+(118)*y[3]+(34)*y[5]+(74)*y[7]+(38)*y[8]
+(66)*y[10]+(-66)*y[11]+(22)*y[13]+(22)*y[14]+(108)*y[20];            
(60)*v.121+(22)*v.122+(118)*v.123+(34)*v.125+(74)*v.127+
(38)*v.128+(66)*v.130+(-66)*v.131+(22)*v.133+(22)*v.134+(108)*v.140
> e*f;
(60)*v.241+(88)*v.242+(118)*v.243+(174)*v.244+(142)*v.245+(108)*v.246
+(74)*v.247+(38)*v.248
> h:=e*f;;
\end{verbatim}
Calcul de {\tt g}, {\tt Bg}, {\tt z} et {\tt Bz}\,:
\begin{verbatim}
> g:=LieCentralizer(L,Subspace(L,[e]));
<Lie algebra of dimension 12 over Rationals>
> Bg:=BasisVectors(Basis(g));;
> z:=LieCentre(g);
<two-sided ideal in <Lie algebra of dimension 12 over Rationals>, 
  (dimension 6)>
> Bz:=BasisVectors(Basis(z));
[ v.1+v.2+v.3+v.5+v.7+v.8+v.10+v.20, 
  v.64+v.65+(2)*v.67+v.69+v.71+v.73+v.74+v.76, v.97+(-1)*v.98+v.99+v.100, 
  v.104+v.107+(-1)*v.108+(-1)*v.110, v.117, v.120 ]
> h*Bz[2];
(14)*v.64+(14)*v.65+(28)*v.67+(14)*v.69+(14)*v.71+
(14)*v.73+(14)*v.74+(14)*v.76
> h*Bz[3];
(22)*v.97+(-22)*v.98+(22)*v.99+(22)*v.100 
> h*Bz[4];
(26)*v.104+(26)*v.107+(-26)*v.108+(-26)*v.110
> h*Bz[5];
(34)*v.117
> h*Bz[6];
(38)*v.120
\end{verbatim}  
Les poids de {\tt z} sont 2,14,22,26,34,38; d'o{\`u} $m_{r}=38$. Il y a quatre matrices {\`a}
{\'e}tudier.

\begin{enumerate}
\item $m_{i_{1}}=14$,  $m_{k_{1}}=26$
\begin{verbatim}
> ((f*Bz[4])*Bz[2]);
(-182)*v.120
\end{verbatim}

\item $m_{i_{1}}=22$,  $m_{k_{1}}=18$
\begin{verbatim}
> h*Bg[6];
(18)*v.81+(36)*v.85+(-18)*v.86+(-18)*v.87+(18)*v.88

> ((f*Bg[6])*Bz[3]);                                                    
(-198)*v.120
\end{verbatim}

\item $m_{i_{1}}=26$,  $m_{k_{1}}=14$
\begin{verbatim}
> ((f*Bz[2])*Bz[4]);
(-182)*v.120
\end{verbatim}

\item $m_{i_{1}}=34$,  $m_{k_{1}}=6$
\begin{verbatim}
> h*Bg[2];
(6)*v.23+(9)*v.24+(-3)*v.25+(9)*v.29+(-3)*v.31+(-9)*v.32+(3)*v.33+(6)*v.34
+(-3)*v.35+(-3)*v.36+(-3)*v.40

> ((f*Bg[2])*Bz[5]);
(51)*v.120
\end{verbatim}
\end{enumerate}
La propri{\'e}t{\'e} $(P)$ est v{\'e}rifi{\'e}e pour cette orbite.\\

\item Caract{\'e}ristique\,:
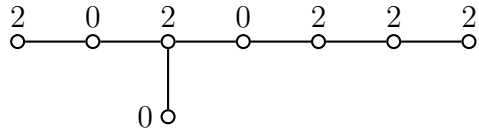
\begin{center}
\begin{pspicture}(-5,-1.2)(10,1)

\pscircle(-2,0){1mm}
\pscircle(-1,0){1mm}
\pscircle(0,0){1mm}
\pscircle(1,0){1mm}
\pscircle(2,0){1mm}
\pscircle(3,0){1mm}
\pscircle(4,0){1mm}

\pscircle(0,-1){1mm}

\psline(-1.9,0)(-1.1,0)
\psline(-0.1,0)(-0.9,0)
\psline(0.1,0)(0.9,0)
\psline(1.1,0)(1.9,0)
\psline(2.1,0)(2.9,0)
\psline(3.1,0)(3.9,0)
\psline(0,-0.1)(0,-0.9)

\rput[b](-2,0.2){$2$}
\rput[b](-1,0.2){$0$}
\rput[b](0,0.2){$2$}
\rput[b](1,0.2){$0$}
\rput[b](2,0.2){$2$}
\rput[b](3,0.2){$2$}
\rput[b](4,0.2){$2$}
\rput[r](-0.2,-1){$0$}
\end{pspicture}
\end{center}
D{\'e}finition du $\mathfrak{sl}_{2}$-triplet\,:
\begin{verbatim}
> e:=x[7]+x[8]+x[9]+x[10]+x[11]+x[12]+x[13]+x[25];                      
v.7+v.8+v.9+v.10+v.11+v.12+v.13+v.25
> f:=(66)*y[7]+(34)*y[8]+(52)*y[9]+(75)*y[10]+(49)*y[11]+(27)*y[12]
+(96)*y[13]+y[25];                                                  
(66)*v.127+(34)*v.128+(52)*v.129+(75)*v.130+(49)*v.131+(27)*v.132
+(96)*v.133+v.145
> e*f;
(52)*v.241+(76)*v.242+(102)*v.243+(152)*v.244+(124)*v.245+(96)*v.246
+(66)*v.247+(34)*v.248
> h:=e*f;;
\end{verbatim}
Calcul de {\tt g}, {\tt Bg}, {\tt z} et {\tt Bz}\,:
\begin{verbatim}
> g:=LieCentralizer(L,Subspace(L,[e]));
<Lie algebra of dimension 14 over Rationals>
> Bg:=BasisVectors(Basis(g));;
> z:=LieCentre(g);
<two-sided ideal in <Lie algebra of dimension 14 over Rationals>, 
  (dimension 6)>
> Bz:=BasisVectors(Basis(z));
[ v.7+v.8+v.9+v.10+v.11+v.12+v.13+v.25,
v.74+(-2)*v.77+(-1)*v.78+(-1)*v.80
+(-1)*v.82, v.104+(-1)*v.105+(-1)*v.106, v.113+v.114, v.117, v.120 ]
> h*Bz[2];
(14)*v.74+(-28)*v.77+(-14)*v.78+(-14)*v.80+(-14)*v.82
> h*Bz[3];
(22)*v.104+(-22)*v.105+(-22)*v.106
> h*Bz[4];
(26)*v.113+(26)*v.114
> h*Bz[5];
(28)*v.117
> h*Bz[6];
(34)*v.120
\end{verbatim}  
Les poids de {\tt z} sont 2,14,22,26,28,34; d'o{\`u} $m_{r}=34$. Il y a quatre matrices {\`a}
{\'e}tudier.

\begin{enumerate}
\item $m_{i_{1}}=14$,  $m_{k_{1}}=22$
\begin{verbatim}
> ((f*Bz[3])*Bz[2]);
(-154)*v.120
\end{verbatim}

\item $m_{i_{1}}=22$,  $m_{k_{1}}=14$
\begin{verbatim}
> ((f*Bz[2])*Bz[3]);
(-154)*v.120
\end{verbatim}

\item $m_{i_{1}}=26$,  $m_{k_{1}}=10$
\begin{verbatim}
> h*Bg[4];
(10)*v.54+(-5)*v.56+(-5)*v.58+(-10)*v.59+(5)*v.61+(-15)*v.63
> h*Bg[5];
(10)*v.62+(-10)*v.64+(-10)*v.69

> ((f*Bg[4])*Bz[4]);
(-65)*v.120
\end{verbatim}

\item $m_{i_{1}}=28$,  $m_{k_{1}}=8$
\begin{verbatim}
> h*Bg[3];
(8)*v.36+(8)*v.39+(8)*v.45+(-24)*v.50+(16)*v.55+(-8)*v.57

> ((f*Bg[3])*Bz[5]);
(24)*v.120
\end{verbatim}
\end{enumerate}
La propri{\'e}t{\'e} $(P)$ est v{\'e}rifi{\'e}e pour cette orbite.\\

\item Caract{\'e}ristique\,:
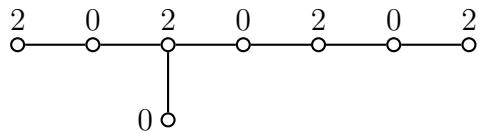
\begin{center}
\begin{pspicture}(-5,-1.2)(10,1)

\pscircle(-2,0){1mm}
\pscircle(-1,0){1mm}
\pscircle(0,0){1mm}
\pscircle(1,0){1mm}
\pscircle(2,0){1mm}
\pscircle(3,0){1mm}
\pscircle(4,0){1mm}

\pscircle(0,-1){1mm}

\psline(-1.9,0)(-1.1,0)
\psline(-0.1,0)(-0.9,0)
\psline(0.1,0)(0.9,0)
\psline(1.1,0)(1.9,0)
\psline(2.1,0)(2.9,0)
\psline(3.1,0)(3.9,0)
\psline(0,-0.1)(0,-0.9)

\rput[b](-2,0.2){$2$}
\rput[b](-1,0.2){$0$}
\rput[b](0,0.2){$2$}
\rput[b](1,0.2){$0$}
\rput[b](2,0.2){$2$}
\rput[b](3,0.2){$0$}
\rput[b](4,0.2){$2$}
\rput[r](-0.2,-1){$0$}
\end{pspicture}
\end{center}
D{\'e}finition du $\mathfrak{sl}_{2}$-triplet\,:
\begin{verbatim}
> e:=x[9]+x[10]+x[11]+x[12]+x[13]+x[14]+x[15]+x[25];                    
v.9+v.10+v.11+v.12+v.13+v.14+v.15+v.25
> f:=(44)*y[9]+(50)*y[10]+(28)*y[11]+(36)*y[12]+(54)*y[13]+(26)*y[14]
+(28)*y[15]+(14)*y[25];                                           
(44)*v.129+(50)*v.130+(28)*v.131+(36)*v.132+(54)*v.133+(26)*v.134
+(28)*v.135+(14)*v.145
> e*f;
(44)*v.241+(64)*v.242+(86)*v.243+(128)*v.244+(104)*v.245+(80)*v.246
+(54)*v.247+(28)*v.248
> h:=e*f;;
\end{verbatim}
Calcul de {\tt g}, {\tt Bg}, {\tt z} et {\tt Bz}\,:
\begin{verbatim}
> g:=LieCentralizer(L,Subspace(L,[e]));
<Lie algebra of dimension 16 over Rationals>
> Bg:=BasisVectors(Basis(g));;
> z:=LieCentre(g);
<two-sided ideal in <Lie algebra of dimension 16 over Rationals>, 
  (dimension 5)>
> Bz:=BasisVectors(Basis(z));
[ v.9+v.10+v.11+v.12+v.13+v.14+v.15+v.25, v.86+v.87+v.88+(4)*v.93+(-3)*v.95, 
  v.111+v.115, v.119, v.120 ]
> h*Bz[2];
(14)*v.86+(14)*v.87+(14)*v.88+(56)*v.93+(-42)*v.95
> h*Bz[3];
(22)*v.111+(22)*v.115
> h*Bz[4];
(26)*v.119
> h*Bz[5];
(28)*v.120
\end{verbatim}  
Les poids de {\tt z} sont 2,14,22,26,28; d'o{\`u} $m_{r}=28$. Il y a trois matrices {\`a}
{\'e}tudier.

\begin{enumerate}
\item $m_{i_{1}}=14$,  $m_{k_{1}}=16$
\begin{verbatim}
> h*Bg[10];
(16)*v.97+(16)*v.98+(-16)*v.99+(-16)*v.100

> ((f*Bg[10])*Bz[2]);
(-112)*v.120
\end{verbatim}

\item $m_{i_{1}}=22$,  $m_{k_{1}}=8$
\begin{verbatim}
> h*Bg[4];
(8)*v.45+(8)*v.47+(-8)*v.57+(-24)*v.58+(24)*v.59+(16)*v.61+(16)*v.62
> h*Bg[5];

> ((f*Bg[4])*Bz[3]); 
(88)*v.120
\end{verbatim}

\item $m_{i_{1}}=26$,  $m_{k_{1}}=4$
\begin{verbatim}
> h*Bg[2];
(4)*v.16+(4)*v.20+(4)*v.22+(-12)*v.30+(8)*v.32+(4)*v.33+(-12)*v.34+(20)*v.35

> ((f*Bg[2])*Bz[4]);
(26)*v.120
\end{verbatim}
\end{enumerate}
La propri{\'e}t{\'e} $(P)$ est v{\'e}rifi{\'e}e pour cette orbite.\\

\item Caract{\'e}ristique\,:
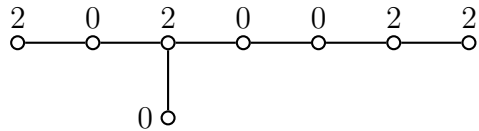
\begin{center}
\begin{pspicture}(-5,-1.2)(10,1)

\pscircle(-2,0){1mm}
\pscircle(-1,0){1mm}
\pscircle(0,0){1mm}
\pscircle(1,0){1mm}
\pscircle(2,0){1mm}
\pscircle(3,0){1mm}
\pscircle(4,0){1mm}

\pscircle(0,-1){1mm}

\psline(-1.9,0)(-1.1,0)
\psline(-0.1,0)(-0.9,0)
\psline(0.1,0)(0.9,0)
\psline(1.1,0)(1.9,0)
\psline(2.1,0)(2.9,0)
\psline(3.1,0)(3.9,0)
\psline(0,-0.1)(0,-0.9)

\rput[b](-2,0.2){$2$}
\rput[b](-1,0.2){$0$}
\rput[b](0,0.2){$2$}
\rput[b](1,0.2){$0$}
\rput[b](2,0.2){$0$}
\rput[b](3,0.2){$2$}
\rput[b](4,0.2){$2$}
\rput[r](-0.2,-1){$0$}
\end{pspicture}
\end{center}
D{\'e}finition du $\mathfrak{sl}_{2}$-triplet\,:
\begin{verbatim}
> e:=x[8]+x[9]+x[10]+x[14]+x[18]+x[19]+x[20]+x[33];                     
v.8+v.9+v.10+v.14+v.18+v.19+v.20+v.33
> f:=(8)*y[8]+(40)*y[9]+(22)*y[10]+(-22)*y[11]+(50)*y[14]+(35)*y[18]
+(37)*y[19]+(21)*y[20]+(35)*y[21]+y[33];                          
(26)*v.128+(40)*v.129+(22)*v.130+(-22)*v.131+(50)*v.134+(37)*v.138
+(37)*v.139+(21)*v.140+(35)*v.141+v.153
> e*f;                                                                  
(40)*v.241+(58)*v.242+(78)*v.243+(116)*v.244+(94)*v.245+(72)*v.246
+(50)*v.247+(26)*v.248
> h:=e*f;;
\end{verbatim}
Calcul de {\tt g}, {\tt Bg}, {\tt z} et {\tt Bz}\,:
\begin{verbatim}
> g:=LieCentralizer(L,Subspace(L,[e]));
<Lie algebra of dimension 18 over Rationals>
> Bg:=BasisVectors(Basis(g));;
> z:=LieCentre(g);
<two-sided ideal in <Lie algebra of dimension 18 over Rationals>, 
  (dimension 4)>
> Bz:=BasisVectors(Basis(z));
[ v.8+v.9+v.10+v.14+v.18+v.19+v.20+v.33, v.94+v.95+v.96+v.97, v.116, v.120 ]
> h*Bz[2];
(14)*v.94+(14)*v.95+(14)*v.96+(14)*v.97
> h*Bz[3];
(22)*v.116
> h*Bz[4];
(26)*v.120
\end{verbatim}  
Les poids de {\tt z} sont 2,14,22,26; d'o{\`u} $m_{r}=26$. Il y a deux matrices {\`a} {\'e}tudier.

\begin{enumerate}
\item $m_{i_{1}}=14$,  $m_{k_{1}}=14$
\begin{verbatim}
> ((f*Bz[2])*Bz[2]);
(-98)*v.120
\end{verbatim}

\item $m_{i_{1}}=22$,  $m_{k_{1}}=6$
\begin{verbatim}
> h*Bg[4];
(6)*v.42+(-12)*v.43+(12)*v.44+(-12)*v.45+(6)*v.46+(-12)*v.49+(6)*v.52
+(-6)*v.55
> h*Bg[5];
(6)*v.50+(-6)*v.51+(-6)*v.57+(-6)*v.61

> ((f*Bg[4])*Bz[3]);                
(-66)*v.120
\end{verbatim}
\end{enumerate}
La propri{\'e}t{\'e} $(P)$ est v{\'e}rifi{\'e}e pour cette orbite.\\

\item Caract{\'e}ristique\,:
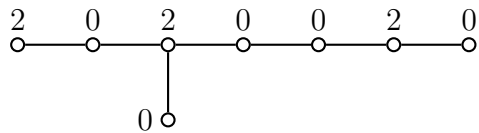
\begin{center}
\begin{pspicture}(-5,-1.2)(10,1)

\pscircle(-2,0){1mm}
\pscircle(-1,0){1mm}
\pscircle(0,0){1mm}
\pscircle(1,0){1mm}
\pscircle(2,0){1mm}
\pscircle(3,0){1mm}
\pscircle(4,0){1mm}

\pscircle(0,-1){1mm}

\psline(-1.9,0)(-1.1,0)
\psline(-0.1,0)(-0.9,0)
\psline(0.1,0)(0.9,0)
\psline(1.1,0)(1.9,0)
\psline(2.1,0)(2.9,0)
\psline(3.1,0)(3.9,0)
\psline(0,-0.1)(0,-0.9)

\rput[b](-2,0.2){$2$}
\rput[b](-1,0.2){$0$}
\rput[b](0,0.2){$2$}
\rput[b](1,0.2){$0$}
\rput[b](2,0.2){$0$}
\rput[b](3,0.2){$2$}
\rput[b](4,0.2){$0$}
\rput[r](-0.2,-1){$0$}
\end{pspicture}
\end{center}
D{\'e}finition du $\mathfrak{sl}_{2}$-triplet\,:
\begin{verbatim}
> e:=x[9]+x[10]+x[14]+x[15]+x[18]+x[19]+x[20]+x[33];                    
v.9+v.10+v.14+v.15+v.18+v.19+v.20+v.33
> f:=(36)*y[9]+(20)*y[10]+(-20)*y[11]+(22)*y[14]+(22)*y[15]+(20)*y[18]
+(22)*y[19]+(30)*y[20]+(20)*y[21]+(12)*y[33];                    
(36)*v.129+(20)*v.130+(-20)*v.131+(22)*v.134+(22)*v.135+(20)*v.138
+(22)*v.139+(30)*v.140+(20)*v.141+(12)*v.153
> e*f;
(36)*v.241+(52)*v.242+(70)*v.243+(104)*v.244+(84)*v.245+(64)*v.246
+(44)*v.247+(22)*v.248
> h:=e*f;;
\end{verbatim}
Calcul de {\tt g}, {\tt Bg}, {\tt z} et {\tt Bz}\,:
\begin{verbatim}
> g:=LieCentralizer(L,Subspace(L,[e]));
<Lie algebra of dimension 20 over Rationals>
> Bg:=BasisVectors(Basis(g));;
> z:=LieCentre(g);
<two-sided ideal in <Lie algebra of dimension 20 over Rationals>, 
  (dimension 4)>
> Bz:=BasisVectors(Basis(z));
[ v.9+v.10+v.14+v.15+v.18+v.19+v.20+v.33, v.97+(-1)*v.99+(3)*v.105+(4)*v.108,
  v.119, v.120 ]
> h*Bz[2];                   
(14)*v.97+(-14)*v.99+(42)*v.105+(56)*v.108
> h*Bz[3];
(22)*v.119
> h*Bz[4];
(22)*v.120
\end{verbatim}  
Les poids de {\tt z} sont 2,14,22,22; d'o{\`u} $m_{r}=22$. Il y a deux matrices {\`a} {\'e}tudier.

\begin{enumerate}
\item $m_{i_{1}}=14$,  $m_{k_{1}}=10$
\begin{verbatim}
> h*Bg[7];
(10)*v.70+(10)*v.80+(10)*v.81+(10)*v.83
> h*Bg[8];
(10)*v.75+(10)*v.76+(10)*v.77+(10)*v.86
> h*Bg[9];
(10)*v.81+(20)*v.83+(10)*v.84+(-10)*v.85+(10)*v.87
> h*Bg[10];
(10)*v.76+(10)*v.77+(-10)*v.79+(10)*v.82+(10)*v.86+(-10)*v.90

> ((f*Bg[7])*Bz[2]);
0*v.1
> ((f*Bg[8])*Bz[2]);
0*v.1
> ((f*Bg[9])*Bz[2]);
(-70)*v.119
> ((f*Bg[10])*Bz[2]); 
(70)*v.120
\end{verbatim}
La matrice correspondante est 
$$\left[
\begin{array}{cccc}
0 & 0 & -70 & 0 \\
0 & 0 & 0 & 70
\end{array}
\right] \cdot$$
C'est clairement une matrice de rang 2.

\item $m_{i_{1}}=22$,  $m_{k_{1}}=2$
\begin{verbatim}
> h*Bg[3];
(2)*v.9+(2)*v.10+(2)*v.14+(2)*v.15+(2)*v.18+(2)*v.19+(2)*v.20+(2)*v.33

> ((f*Bg[1])*Bz[3]); ((f*Bg[2])*Bz[3]);((f*Bg[3])*Bz[3]);
(-2)*v.119
(14)*v.120
(-22)*v.119
> ((f*Bg[1])*Bz[4]); ((f*Bg[2])*Bz[4]); ((f*Bg[3])*Bz[4]);
0*v.1
(22)*v.119
(-22)*v.120
\end{verbatim}
La matrice {\`a} {\'e}tudier est
$$\left[
\begin{array}{cccc}
-2\alpha & 22 \beta &-22 \alpha\\
0 & 14\alpha & -22 \beta
\end{array}
\right]
 \cdot$$
Cette matrice est de rang 2 pour tout couple
 $(\alpha,\beta)$ non nul. 
\end{enumerate}
La propri{\'e}t{\'e} $(P)$ est v{\'e}rifi{\'e}e pour cette orbite.\\

\item Caract{\'e}ristique\,:
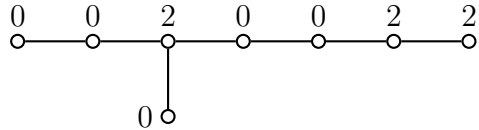
\begin{center}
\begin{pspicture}(-5,-1.2)(10,1)

\pscircle(-2,0){1mm}
\pscircle(-1,0){1mm}
\pscircle(0,0){1mm}
\pscircle(1,0){1mm}
\pscircle(2,0){1mm}
\pscircle(3,0){1mm}
\pscircle(4,0){1mm}

\pscircle(0,-1){1mm}

\psline(-1.9,0)(-1.1,0)
\psline(-0.1,0)(-0.9,0)
\psline(0.1,0)(0.9,0)
\psline(1.1,0)(1.9,0)
\psline(2.1,0)(2.9,0)
\psline(3.1,0)(3.9,0)
\psline(0,-0.1)(0,-0.9)

\rput[b](-2,0.2){$0$}
\rput[b](-1,0.2){$0$}
\rput[b](0,0.2){$2$}
\rput[b](1,0.2){$0$}
\rput[b](2,0.2){$0$}
\rput[b](3,0.2){$2$}
\rput[b](4,0.2){$2$}
\rput[r](-0.2,-1){$0$}
\end{pspicture}
\end{center}
D{\'e}finition du $\mathfrak{sl}_{2}$-triplet\,:
\begin{verbatim}
> e:=x[8]+x[14]+x[16]+x[17]+x[18]+x[19]+x[20]+x[38];
v.8+v.14+v.16+v.17+v.18+v.19+v.20+v.38
> f:=(22)*y[8]+(42)*y[14]+(16)*y[16]+(2)*y[17]+(30)*y[18]+(30)*y[19]
+(2)*y[20]+(16)*y[38];                                                   
(22)*v.128+(42)*v.134+(16)*v.136+(2)*v.137+(30)*v.138+(30)*v.139
+(2)*v.140+(16)*v.158
> e*f;
(32)*v.241+(48)*v.242+(64)*v.243+(96)*v.244+(78)*v.245+(60)*v.246
+(42)*v.247+(22)*v.248
> h:=e*f;;
\end{verbatim}
Calcul de {\tt g}, {\tt Bg}, {\tt z} et {\tt Bz}\,:
\begin{verbatim}
> g:=LieCentralizer(L,Subspace(L,[e]));
<Lie algebra of dimension 22 over Rationals>
> Bg:=BasisVectors(Basis(g));;
> z:=LieCentre(g);
<two-sided ideal in <Lie algebra of dimension 22 over Rationals>, 
  (dimension 5)>
> Bz:=BasisVectors(Basis(z));
[ v.8+v.14+v.16+v.17+v.18+v.19+v.20+v.38, v.99+(-1)*v.107+(-1)*v.108, v.116, 
  v.118, v.120 ]
> h*Bz[2];
(14)*v.99+(-14)*v.107+(-14)*v.108
> h*Bz[3];
(18)*v.116
> h*Bz[4];
(18)*v.118
> h*Bz[5];
(22)*v.120
\end{verbatim}  
Les poids de {\tt z} sont 2,14,18,18,22; d'o{\`u} $m_{r}=22$. Il y a deux matrices {\`a}
{\'e}tudier.

\begin{enumerate}
\item $m_{i_{1}}=14$,  $m_{k_{1}}=10$
\begin{verbatim}
> h*Bg[11];
(10)*v.83+(10)*v.89
> h*Bg[12];
(10)*v.86+(-10)*v.93
> h*Bg[13];
(10)*v.74+(-10)*v.77+(-10)*v.90+(20)*v.97

> ((f*Bg[13])*Bz[2]);
(-70)*v.120
\end{verbatim}

\item $m_{i_{1}}=18$,  $m_{k_{1}}=6$
\begin{verbatim}
> h*Bg[6];
(6)*v.47+(-12)*v.50+(-6)*v.58+(6)*v.59+(12)*v.61+(18)*v.63
> h*Bg[7];
(6)*v.36+(-6)*v.49+(-3)*v.54+(3)*v.64+(-3)*v.66+(-9)*v.69

> ((f*Bg[6])*Bz[3]); ((f*Bg[7])*Bz[3]);
(-12)*v.120
0*v.1
> ((f*Bg[6])*Bz[4]); ((f*Bg[7])*Bz[4]);
0*v.1
(6)*v.120
\end{verbatim}
La matrice {\`a} {\'e}tudier est 
$\left[
\begin{array}{cc}
-12 \alpha & 6\beta 
\end{array}
\right]$; elle est de rang 1 si le couple $(\alpha,\beta)$ est non
nul. 
\end{enumerate}
La propri{\'e}t{\'e} $(P)$ est v{\'e}rifi{\'e}e pour cette orbite.\\

\item Caract{\'e}ristique\,:
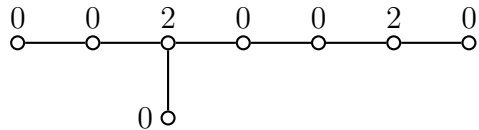
\begin{center}
\begin{pspicture}(-5,-1.2)(10,1)

\pscircle(-2,0){1mm}
\pscircle(-1,0){1mm}
\pscircle(0,0){1mm}
\pscircle(1,0){1mm}
\pscircle(2,0){1mm}
\pscircle(3,0){1mm}
\pscircle(4,0){1mm}

\pscircle(0,-1){1mm}

\psline(-1.9,0)(-1.1,0)
\psline(-0.1,0)(-0.9,0)
\psline(0.1,0)(0.9,0)
\psline(1.1,0)(1.9,0)
\psline(2.1,0)(2.9,0)
\psline(3.1,0)(3.9,0)
\psline(0,-0.1)(0,-0.9)

\rput[b](-2,0.2){$0$}
\rput[b](-1,0.2){$0$}
\rput[b](0,0.2){$2$}
\rput[b](1,0.2){$0$}
\rput[b](2,0.2){$0$}
\rput[b](3,0.2){$2$}
\rput[b](4,0.2){$0$}
\rput[r](-0.2,-1){$0$}
\end{pspicture}
\end{center}
D{\'e}finition du $\mathfrak{sl}_{2}$-triplet\,:
\begin{verbatim}
> e:=x[14]+x[15]+x[16]+x[17]+x[18]+x[19]+x[20]+x[38];                      
v.14+v.15+v.16+v.17+v.18+v.19+v.20+v.38
> f:=(18)*y[14]+(18)*y[15]+(8)*y[16]+(8)*y[17]+(14)*y[18]+(20)*y[19]
+(14)*y[20]+(20)*y[38];                                                  
(18)*v.134+(18)*v.135+(8)*v.136+(8)*v.137+(14)*v.138+(20)*v.139+(14)*v.140
+(20)*v.158
> e*f;
(28)*v.241+(42)*v.242+(56)*v.243+(84)*v.244+(68)*v.245+(52)*v.246
+(36)*v.247+(18)*v.248
> h:=e*f;;
\end{verbatim}
Calcul de {\tt g}, {\tt Bg}, {\tt z} et {\tt Bz}\,:
\begin{verbatim}
> g:=LieCentralizer(L,Subspace(L,[e]));
<Lie algebra of dimension 24 over Rationals>
> Bg:=BasisVectors(Basis(g));;
> z:=LieCentre(g);
<two-sided ideal in <Lie algebra of dimension 24 over Rationals>, 
  (dimension 4)>
> Bz:=BasisVectors(Basis(z));
[ v.14+v.15+v.16+v.17+v.18+v.19+v.20+v.38, v.111+v.112, v.119, v.120 ]
> h*Bz[2];                   
(14)*v.111+(14)*v.112
> h*Bz[3];
(18)*v.119
> h*Bz[4];
(18)*v.120
\end{verbatim}  
Les poids de {\tt z} sont 2,14,18,18; d'o{\`u} $m_{r}=18$. Il y a deux matrices {\`a} {\'e}tudier.

\begin{enumerate}
\item $m_{i_{1}}=14$,  $m_{k_{1}}=6$
\begin{verbatim}
> h*Bg[5];
(6)*v.58+(-6)*v.59+(6)*v.60+(-6)*v.63
> h*Bg[6];
(6)*v.49+(6)*v.64+(6)*v.65+(-6)*v.66+(-6)*v.67+(6)*v.68
> h*Bg[7];
(6)*v.64+(12)*v.65+(-6)*v.66+(-6)*v.67+(12)*v.68+(-6)*v.69
> h*Bg[8];
(6)*v.53+(-6)*v.55+(-6)*v.56+(-6)*v.71+(-6)*v.72+(6)*v.73
> h*Bg[9];
(6)*v.60+(-3)*v.61+(-3)*v.62+(-3)*v.63+(-3)*v.78

> ((f*Bg[5])*Bz[2]);
(14)*v.119
> ((f*Bg[6])*Bz[2]);
0*v.1
> ((f*Bg[7])*Bz[2]);
(14)*v.120
> ((f*Bg[8])*Bz[2]);
0*v.1
> ((f*Bg[9])*Bz[2]);
(7)*v.119
\end{verbatim}
 matrice {\`a} {\'e}tudier est 
$$\left[
\begin{array}{cccccc}
14 & 0 & 0 & 0 & 7 \\
0 & 0 & 14  & 0
\end{array}
\right] \cdot $$ 
C'est clairement une matrice de rang 2.

\item $m_{i_{1}}=18$,  $m_{k_{1}}=2$
\begin{verbatim}
> h*Bg[2];
(2)*v.7+(-10)*v.10+(8)*v.11+(2)*v.12+(-4)*v.29+(2)*v.30+(-4)*v.31+(-4)*v.33
> h*Bg[3];
(2)*v.14+(2)*v.15+(2)*v.16+(2)*v.17+(2)*v.18+(2)*v.19+(2)*v.20+(2)*v.38

> ((f*Bg[1])*Bz[3]); ((f*Bg[2])*Bz[3]);((f*Bg[3])*bz[3]);
(-18/5)*v.120
0*v.1
(-18)*v.119
> ((f*Bg[1])*Bz[4]) ;((f*Bg[2])*Bz[4]); ((f*Bg[3])*Bz[4]);
0*v.1
(18)*v.119
(-18)*v.120
\end{verbatim}
La matrice {\`a} {\'e}tudier est 
$$\left[
\begin{array}{ccc}
0 & 18\beta & -18\alpha \\
-18/5 \alpha & 0 & -18\beta
\end{array}
\right] \cdot $$ 
C'est une matrice de rang 2 pour tout couple $(\alpha,\beta)$ non
nul. 
\end{enumerate}
La propri{\'e}t{\'e} $(P)$ est v{\'e}rifi{\'e}e pour cette orbite.\\

\item Caract{\'e}ristique\,:
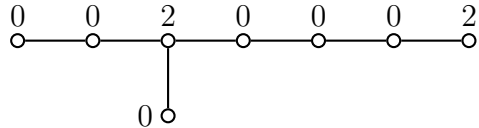
\begin{center}
\begin{pspicture}(-5,-1.2)(10,1)

\pscircle(-2,0){1mm}
\pscircle(-1,0){1mm}
\pscircle(0,0){1mm}
\pscircle(1,0){1mm}
\pscircle(2,0){1mm}
\pscircle(3,0){1mm}
\pscircle(4,0){1mm}

\pscircle(0,-1){1mm}

\psline(-1.9,0)(-1.1,0)
\psline(-0.1,0)(-0.9,0)
\psline(0.1,0)(0.9,0)
\psline(1.1,0)(1.9,0)
\psline(2.1,0)(2.9,0)
\psline(3.1,0)(3.9,0)
\psline(0,-0.1)(0,-0.9)

\rput[b](-2,0.2){$0$}
\rput[b](-1,0.2){$0$}
\rput[b](0,0.2){$2$}
\rput[b](1,0.2){$0$}
\rput[b](2,0.2){$0$}
\rput[b](3,0.2){$0$}
\rput[b](4,0.2){$2$}
\rput[r](-0.2,-1){$0$}
\end{pspicture}
\end{center}
D{\'e}finition du $\mathfrak{sl}_{2}$-triplet\,:
\begin{verbatim}
> e:=x[15]+x[16]+x[17]+x[18]+x[19]+x[20]+x[38]+x[46];                      
v.15+v.16+v.17+v.18+v.19+v.20+v.38+v.46
> f:=(8)*y[8]+(16)*y[15]+(2)*y[16]+(12)*y[17]+(2)*y[18]+(12)*y[19]
+(22)*y[20]+(-14)*y[28]+(8)*y[38]+(14)*y[46];                              
(8)*v.128+(16)*v.135+(2)*v.136+(12)*v.137+(2)*v.138+(12)*v.139
+(22)*v.140+(-14)*v.148+(8)*v.158+(14)*v.166
> e*f;
(24)*v.241+(36)*v.242+(48)*v.243+(72)*v.244+(58)*v.245+(44)*v.246
+(30)*v.247+(16)*v.248
>h:=e*f;;
\end{verbatim}
\rem Il faut pr{\'e}ciser pour cette orbite qu'avec les conventions de
\cite{Triplets}, l'{\'e}l{\'e}ment $X_{47}$ correspond {\`a} l'{\'e}l{\'e}ment {\tt
  x[46]=v.46} du logiciel, c'est pourquoi la d{\'e}finition du
$\mathfrak{sl}_{2}$-triplet est bien en accord avec \cite{Triplets}.\\
\\ 
Calcul de {\tt g}, {\tt Bg}, {\tt z} et {\tt Bz}\,:
\begin{verbatim}
> g:=LieCentralizer(L,Subspace(L,[e]));
<Lie algebra of dimension 28 over Rationals>
> Bg:=BasisVectors(Basis(g));;
> z:=LieCentre(g);
<two-sided ideal in <Lie algebra of dimension 28 over Rationals>, 
  (dimension 3)>
> Bz:=BasisVectors(Basis(z));
[ v.15+v.16+v.17+v.18+v.19+v.20+v.38+v.46, v.119, v.120 ]
> h*Bz[2];
(14)*v.119
> h*Bz[3];
(16)*v.120
\end{verbatim}  
Les poids de {\tt z} sont 2,14,16; d'o{\`u} $m_{r}=16$. Il n'y a qu'une matrice {\`a} {\'e}tudier.

\begin{enumerate}
\item $m_{i_{1}}=14$,  $m_{k_{1}}=4$
\begin{verbatim}
> h*Bg[3];
(4)*v.37
> h*Bg[6];
(4)*v.42+(-8)*v.44+(4)*v.45+(4)*v.48+(-4)*v.53+(4)*v.55
> h*Bg[7];
(4)*v.32+(2)*v.47+(2)*v.51+(2)*v.52+(4)*v.59
> h*Bg[8];
(4)*v.36+(-4)*v.40+(-8)*v.49+(-4)*v.54+(4)*v.57+(-4)*v.64+(4)*v.66

> ((f*Bg[8])*Bz[2]);
(-28)*v.120
\end{verbatim}
\end{enumerate}
La propri{\'e}t{\'e} $(P)$ est v{\'e}rifi{\'e}e pour cette orbite.\\

\item Caract{\'e}ristique\,:
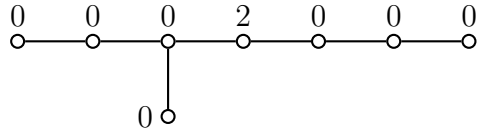
\begin{center}
\begin{pspicture}(-5,-1.2)(10,1)

\pscircle(-2,0){1mm}
\pscircle(-1,0){1mm}
\pscircle(0,0){1mm}
\pscircle(1,0){1mm}
\pscircle(2,0){1mm}
\pscircle(3,0){1mm}
\pscircle(4,0){1mm}

\pscircle(0,-1){1mm}

\psline(-1.9,0)(-1.1,0)
\psline(-0.1,0)(-0.9,0)
\psline(0.1,0)(0.9,0)
\psline(1.1,0)(1.9,0)
\psline(2.1,0)(2.9,0)
\psline(3.1,0)(3.9,0)
\psline(0,-0.1)(0,-0.9)

\rput[b](-2,0.2){$0$}
\rput[b](-1,0.2){$0$}
\rput[b](0,0.2){$0$}
\rput[b](1,0.2){$2$}
\rput[b](2,0.2){$0$}
\rput[b](3,0.2){$0$}
\rput[b](4,0.2){$0$}
\rput[r](-0.2,-1){$0$}
\end{pspicture}
\end{center}
D{\'e}finition du $\mathfrak{sl}_{2}$-triplet\,:
\begin{verbatim}
> e:=x[12]+x[21]+x[30]+x[31]+x[33]+x[42]+x[43]+x[53];                      
v.12+v.21+v.30+v.31+v.33+v.42+v.43+v.53
> f:=(5)*y[12]+y[21]+(5)*y[30]+(2)*y[31]+(8)*y[33]+(2)*y[42]+(8)*y[43]
+(9)*y[53];
(5)*v.132+v.141+(5)*v.150+(2)*v.151+(8)*v.153+(2)*v.162+(8)*v.163+(9)*v.173
> e*f;
(16)*v.241+(24)*v.242+(32)*v.243+(48)*v.244+(40)*v.245+(30)*v.246
+(20)*v.247+(10)*v.248
> h:=e*f;;
\end{verbatim}
\rem Ici encore, il faut pr{\'e}ciser que l'{\'e}l{\'e}ment $X_{32}$ correspond 
{\`a} l'{\'e}l{\'e}ment {\tt x[31]=v.31} du logiciel et la d{\'e}finition du
$\mathfrak{sl}_{2}$-triplet est bien en accord avec \cite{Triplets}.\\
\\ 
Calcul de {\tt g}, {\tt Bg}, {\tt z} et {\tt Bz}\,:
\begin{verbatim}
> g:=LieCentralizer(L,Subspace(L,[e]));
<Lie algebra of dimension 40 over Rationals>
> Bg:=BasisVectors(Basis(g));;
> z:=LieCentre(g);
<two-sided ideal in <Lie algebra of dimension 40 over Rationals>, 
  (dimension 5)>
> Bz:=BasisVectors(Basis(z));;
\end{verbatim}  
Les poids de {\tt z} sont 2,10,10,10,10; d'o{\`u} $m_{r}=10$. Il n'y a qu'une matrice {\`a}
{\'e}tudier.

\begin{enumerate}
\item $m_{i_{1}}=10$,  $m_{k_{1}}=2$
\begin{verbatim}
> ((f*Bg[1])*Bz[2]); ((f*Bg[2])*Bz[2]); ((f*Bg[3])*Bz[2]); ((f*Bg[4])*Bz[2]);
((f*Bg[5])*Bz[2]); ((f*Bg[6])*Bz[2]); ((f*Bg[7])*Bz[2]);
((f*Bg[8])*Bz[2]); ((f*Bg[9])*Bz[2]); ((f*Bg[10])*Bz[2]);
0*v.1
0*v.1
0*v.1
(-5)*v.120
(-10)*v.117
0*v.1
0*v.1
(-5)*v.118
0*v.1
(-10/3)*v.119
> ((f*Bg[1])*Bz[3]); ((f*Bg[2])*Bz[3]); ((f*Bg[3])*Bz[3]); ((f*Bg[4])*Bz[3]);
((f*Bg[5])*Bz[3]); ((f*Bg[6])*Bz[3]); ((f*Bg[7])*Bz[3]);
((f*Bg[8])*Bz[3]); ((f*Bg[9])*Bz[3]); ((f*Bg[10])*Bz[3]);                     
0*v.1
(-2)*v.118
(-2)*v.120
(-2)*v.117
(-8)*v.118
0*v.1
0*v.1
v.120
(4/3)*v.119
0*v.1
> ((f*Bg[1])*Bz[4]); ((f*Bg[2])*Bz[4]); ((f*Bg[3])*Bz[4]); ((f*Bg[4])*Bz[4]);
((f*Bg[5])*Bz[4]); ((f*Bg[6])*Bz[4]); ((f*Bg[7])*Bz[4]);
((f*Bg[8])*Bz[4]); ((f*Bg[9])*Bz[4]); ((f*Bg[10])*Bz[4]);
(-1)*v.119
0*v.1
0*v.1
0*v.1
(-9)*v.119
0*v.1
(-2)*v.118
0*v.1
(-2)*v.120
(2)*v.117
> ((f*Bg[1])*Bz[5]); ((f*Bg[2])*Bz[5]); ((f*Bg[3])*Bz[5); ((f*Bg[4])*Bz[5]);
((f*Bg[5])*Bz[5]); ((f*Bg[6])*Bz[5]); ((f*Bg[7])*Bz[5]);
((f*Bg[8])*Bz[5]); ((f*Bg[9])*Bz[5]); ((f*Bg[10])*Bz[5]);
0*v.1
(-2)*v.120
0*v.1
v.118
(-8)*v.120
(-2)*v.118
(4/3)*v.119
(-2)*v.117
0*v.1
0*v.1
\end{verbatim}
La matrice {\`a} {\'e}tudier est de taille $4 \times 10$\,:
$$\left[
\begin{array}{cccccccccc}
0 & 0 & 0 & -2 \beta &-10\alpha &0 &0 &-2 \delta &0 & 2 \gamma \\
0 & -2\beta & 0 & \delta & -8 \beta & -2 \delta & -2 \gamma  & -5 \alpha & 0
& 0\\
-\gamma &0 & 0& 0&-9\gamma &0 &4/3 \delta & 0& 4/3 \beta & -10/3 \alpha \\
0 &-2 \delta & -2\beta&-5 \alpha &-8 \delta & 0&0 & \beta &-2 \gamma & 0
\end{array}
\right] \cdot
$$
Une {\'e}tude {\'e}l{\'e}mentaire permet de voir que cette matrice est de
rang 4, pour tout 4-uplet $(\alpha,\beta,\gamma,\delta)$ non nul.
\end{enumerate} 
La propri{\'e}t{\'e} $(P)$ est v{\'e}rifi{\'e}e pour cette orbite.\\
\end{enumerate}
{\bf Conclusion\,:} $E_{8}$ v{\'e}rifie $(P)$.
 
\subsection{Cas de $F_{4}$.}
\noindent D{\'e}finition de {\tt L}\,:
\begin{verbatim}
> L:=SimpleLieAlgebra("F",4,Rationals);
<Lie algebra of dimension 52 over Rationals>
> R:=RootSystem(L);
<root system of rank 4>
> P:=PositiveRoots(R);; 
> x:=PositiveRootVectors(R);
[ v.1, v.2, v.3, v.4, v.5, v.6, v.7, v.8, v.9, v.10, v.11, v.12, v.13, v.14, 
  v.15, v.16, v.17, v.18, v.19, v.20, v.21, v.22, v.23, v.24 ]
> y:=NegativeRootVectors(R);
[ v.25, v.26, v.27, v.28, v.29, v.30, v.31, v.32, v.33, v.34, v.35, v.36, 
  v.37, v.38, v.39, v.40, v.41, v.42, v.43, v.44, v.45, v.46, v.47, v.48 ]
> CanonicalGenerators(R)[3];
[ v.49, v.50, v.51, v.52 ]
\end{verbatim}
Dans $F_{4}$, il y a trois orbites nilpotentes distingu{\'e}es non
r{\'e}guli{\`e}res. Pour $F_{4}$, les conventions du logiciel GAP4 sont tr{\`e}s
diff{\'e}rentes de celles adopt{\'e}es dans \cite{Triplets}; dans \cite{Triplets}, le
diagramme de Dynkin est\,:
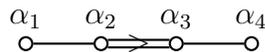
\begin{center}
\begin{pspicture}(-5,-0.2)(10,1)

\pscircle(-2,0){1mm}
\pscircle(-1,0){1mm}
\pscircle(0,0){1mm}
\pscircle(1,0){1mm}

\psline(-1.9,0)(-1.1,0)
\psline(-0.09,0.05)(-0.91,0.05)
\psline(-0.09,-0.05)(-0.91,-0.05)
\psline(0.1,0)(0.9,0)

\rput[b](-2,0.2){$\alpha_{1}$}
\rput[b](-1,0.2){$\alpha_{2}$}
\rput[b](0,0.2){$\alpha_{3}$}
\rput[b](1,0.2){$\alpha_{4}$}

\rput(-0.5,0){$>$}
\end{pspicture}
\end{center}
Il semble que dans GAP4 le diagramme de Dynkin soit plut{\^o}t le suivant\,:
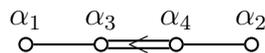
\begin{center}
\begin{pspicture}(-5,-0.2)(10,1)

\pscircle(-2,0){1mm}
\pscircle(-1,0){1mm}
\pscircle(0,0){1mm}
\pscircle(1,0){1mm}

\psline(-1.9,0)(-1.1,0)
\psline(-0.09,0.05)(-0.91,0.05)
\psline(-0.09,-0.05)(-0.91,-0.05)
\psline(0.1,0)(0.9,0)

\rput[b](-2,0.2){$\alpha_{1}$}
\rput[b](-1,0.2){$\alpha_{3}$}
\rput[b](0,0.2){$\alpha_{4}$}
\rput[b](1,0.2){$\alpha_{2}$}

\rput(-0.5,0){$<$}
\end{pspicture}
\end{center}
Par cons{\'e}quent, il est difficile d'utiliser directement les donn{\'e}es de
\cite{Triplets} dans GAP4. On utilise les correspondances suivantes\,:
{\tt x[1]}=$X_{4}$, {\tt x[2]}=$X_{1}$, {\tt x[3]}=$X_{3}$, {\tt x[4]}=$X_{1}$, {\tt x[5]}=$X_{7}$, {\tt
  x[6]}=$X_{5}$, {\tt x[6]}=$X_{5}$, {\tt x[7]}=$X_{6}$, {\tt
  x[8]}=$X_{10}$, {\tt x[9]}=$X_{8}$,  {\tt x[10]}=$X_{9}$,  {\tt
  x[18]}=$X_{18}$. Cependant, m{\^e}me avec ces relations, les
$\mathfrak{sl}_{2}$-triplets de \cite{Triplets} ne conviennent pas. On
utilise la commande {\tt FindSl2} qui permet de chercher une
sous-alg{\`e}bre {\tt s} isomorphe {\`a} $\mathfrak{sl}_{2}$ et contenant {\tt
  e}. On s'assure auparavant que l'{\'e}l{\'e}m{\'e}nt {\tt e} est bien nilpotent {\`a}
l'aide de la commande {\tt IsNilpotentElement} et on v{\'e}rifie aussi que
le $\mathfrak{sl}_{2}$-triplet obtenu correspond bien {\`a} la
caract{\'e}ristique voulue.\\

\begin{enumerate}
\item Caract{\'e}ristique\,:
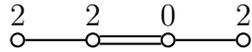
\begin{center}
\begin{pspicture}(-5,-0.2)(10,1)

\pscircle(-2,0){1mm}
\pscircle(-1,0){1mm}
\pscircle(0,0){1mm}
\pscircle(1,0){1mm}

\psline(-1.9,0)(-1.1,0)
\psline(-0.09,0.05)(-0.91,0.05)
\psline(-0.09,-0.05)(-0.91,-0.05)
\psline(0.1,0)(0.9,0)

\rput[b](-2,0.2){$2$}
\rput[b](-1,0.2){$2$}
\rput[b](0,0.2){$0$}
\rput[b](1,0.2){$2$}
\end{pspicture}
\end{center}
D{\'e}finition du $\mathfrak{sl}_{2}$-triplet\,:
\begin{verbatim}
> e:=x[2]+x[4]+x[5]+x[7];
v.2+v.4+v.5+v.7
> IsNilpotentElement(L,e);                   
true
> s:=FindSl2(L,a);                           
<Lie algebra of dimension 3 over Rationals>
> Bs:=BasisVectors(Basis(s));                
[ v.2+v.4+v.5+v.7, v.49+(7/5)*v.50+(9/5)*v.51+(13/5)*v.52, 
  v.25+(7/5)*v.26+v.28+v.29+(4/5)*v.31+(-4/5)*v.34 ]
\end{verbatim}
L'{\'e}l{\'e}ment central de cette base est {\'e}gal au dixi{\`e}me de l'{\'e}lement neutre correspondant {\`a} la caract{\'e}ristique; par suite
en prenant pour {\tt f} dix fois le troisi{\`e}me {\'e}l{\'e}ment de cette
base, on obtient un  $\mathfrak{sl}_{2}$-triplet pour cette caract{\'e}ristique.
\begin{verbatim}
> h:=(14)*H[2]+(26)*H[4]+(18)*H[3]+(10)*H[1];
(10)*v.49+(14)*v.50+(18)*v.51+(26)*v.52
> f:=(10)*y[1]+(14)*y[2]+(10)*y[4]+(10)*y[5]+(8)*y[7]+(-8)*y[10];
(10)*v.25+(14)*v.26+(10)*v.28+(10)*v.29+(8)*v.31+(-8)*v.34
> e*f;
(10)*v.49+(14)*v.50+(18)*v.51+(26)*v.52
\end{verbatim}
Calcul de {\tt g}, {\tt Bg}, {\tt z} et {\tt Bz}\,:
\begin{verbatim}
> g:=LieCentralizer(L,Subspace(L,[e]));
<Lie algebra of dimension 6 over Rationals>
> Bg:=BasisVectors(Basis(g));;
> z:=LieCentre(g);
<two-sided ideal in <Lie algebra of dimension 6 over Rationals>, 
(dimension 3 )>
> Bz:=BasisVectors(Basis(z));
[ v.2+v.4+v.5+v.7, v.20+v.21+(2)*v.22, v.24 ]
> h*Bz[2];
(10)*v.20+(10)*v.21+(20)*v.22
> h*Bz[3];
(14)*v.2
\end{verbatim}  
Les poids de {\tt z} sont 2,10,14; d'o{\`u} $m_{r}=14$. Il n'y a qu'une matrice {\`a} {\'e}tudier.

\begin{enumerate}
\item $m_{i_{1}}=10$,  $m_{k_{1}}=6$
\begin{verbatim}
> h*Bg[3];
(6)*v.11+(6)*v.14+(-12)*v.15+(6)*v.16

> ((f*Bg[3])*Bz[2]);
(30)*v.24
\end{verbatim}
\end{enumerate}
La propri{\'e}t{\'e} $(P)$ est v{\'e}rifi{\'e}e pour cette orbite.\\

\item Caract{\'e}ristique\,:
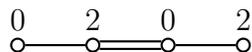
\begin{center}
\begin{pspicture}(-5,-0.2)(10,1)

\pscircle(-2,0){1mm}
\pscircle(-1,0){1mm}
\pscircle(0,0){1mm}
\pscircle(1,0){1mm}

\psline(-1.9,0)(-1.1,0)
\psline(-0.09,0.05)(-0.91,0.05)
\psline(-0.09,-0.05)(-0.91,-0.05)
\psline(0.1,0)(0.9,0)

\rput[b](-2,0.2){$0$}
\rput[b](-1,0.2){$2$}
\rput[b](0,0.2){$0$}
\rput[b](1,0.2){$2$}
\end{pspicture}
\end{center}
D{\'e}finition du $\mathfrak{sl}_{2}$-triplet\,:
\begin{verbatim}
> e:=x[1]+x[6]+x[5]+x[10];
v.1+v.5+v.6+v.10
> IsNilpotentElement(L,e);
true
> s:=FindSl2(L,e);
<Lie algebra of dimension 3 over Rationals>
> Bs:=BasisVectors(Basis(s));
[ v.1+v.5+v.6+v.10, v.49+(5/4)*v.50+(7/4)*v.51+(5/2)*v.52, 
  v.25+v.29+(5/2)*v.30+v.31+v.33+(5/2)*v.34 ]
\end{verbatim}
L'{\'e}l{\'e}ment central de cette base est {\'e}gal au huiti{\`e}me de l'{\'e}l{\'e}ment neutre de
la carat{\'e}ristique. Ici {\tt Bs[1]*Bs[3]=2*Bs[2]}; on prend alors pour
{\tt f} quatre fois le troisi{\`e}me {\'e}l{\'e}ment de cette base. 
\begin{verbatim}
> h:=(10)*H[2]+(20)*H[4]+(14)*H[3]+(8)*H[1];
(8)*v.49+(10)*v.50+(14)*v.51+(20)*v.52
> f:=(4)*y[1]+(4)*y[5]+(10)*y[6]+(4)*y[7]+(4)*y[9]+(10)*y[10];
(4)*v.25+(4)*v.29+(10)*v.30+(4)*v.31+(4)*v.33+(10)*v.34
> e*f;
(8)*v.49+(10)*v.50+(14)*v.51+(20)*v.52
\end{verbatim}
Calcul de {\tt g}, {\tt Bg}, {\tt z} et {\tt Bz}\,:
\begin{verbatim}
> g:=LieCentralizer(L,Subspace(L,[e]));
<Lie algebra of dimension 8 over Rationals>
> Bg:=BasisVectors(Basis(g));;
> z:=LieCentre(g);
<two-sided ideal in <Lie algebra of dimension 8 over Rationals>, 
(dimension 3)>
> Bz:=BasisVectors(Basis(z));
[ v.1+v.5+v.6+v.10, v.23, v.24 ]
> h*Bz[2];
(10)*v.23
> h*Bz[3];
(10)*v.24
\end{verbatim}  
Les poids de {\tt z} sont 2,10,10; d'o{\`u} $m_{r}=10$. Il n'y a qu'une matrice {\`a} {\'e}tudier.

\begin{enumerate}
\item $m_{i_{1}}=10$,  $m_{k_{1}}=2$
\begin{verbatim}
> h*Bg[2];
(2)*v.4+(2)*v.5+(2)*v.6+(-2)*v.7+(2)*v.10
> h*Bg[3];
(2)*v.5+(2)*v.9+(-2)*v.13
> h*Bg[4];
(4)*v.8+(4)*v.11+(-4)*v.12+(-4)*v.14+(-4)*v.16

> ((f*Bg[1])*Bz[2]); ((f*Bg[2])*Bz[2]); ((f*Bg[3])*Bz[2]);
(-8)*v.23
(-2)*v.23
(-2)*v.24
> ((f*Bg[1])*Bz[3]); ((f*Bg[2])*Bz[3]); ((f*Bg[3])*Bz[3]);
(-2)*v.23
(2)*v.23+(-10)*v.24
(-8)*v.24
\end{verbatim}
La matrice {\`a} {\'e}tudier est 
$$\left[
\begin{array}{ccc}
-8 \alpha -2 \beta & -2 \alpha +2 \beta & 0          \\
0 & -10 \beta  & -2 \alpha -8 \beta 
\end{array}
\right] \cdot $$ 
Une br{\`e}ve {\'e}tude de cette matrice montre qu'elle est de rang 2 pour
tout couple $(\alpha,\beta)$ non
nul. 
\end{enumerate}
La propri{\'e}t{\'e} $(P)$ est v{\'e}rifi{\'e}e pour cette orbite.\\

\item Caract{\'e}ristique\,:
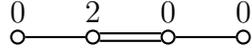
\begin{center}
\begin{pspicture}(-5,-0.2)(10,1)

\pscircle(-2,0){1mm}
\pscircle(-1,0){1mm}
\pscircle(0,0){1mm}
\pscircle(1,0){1mm}

\psline(-1.9,0)(-1.1,0)
\psline(-0.09,0.05)(-0.91,0.05)
\psline(-0.09,-0.05)(-0.91,-0.05)
\psline(0.1,0)(0.9,0)

\rput[b](-2,0.2){$0$}
\rput[b](-1,0.2){$2$}
\rput[b](0,0.2){$0$}
\rput[b](1,0.2){$0$}
\end{pspicture}
\end{center}
D{\'e}finition du $\mathfrak{sl}_{2}$-triplet\,:
\begin{verbatim}
> e:=x[9]+x[10]+x[8]+x[18];   
v.8+v.9+v.10+v.18
> f:=2*y[9]+2*y[10]+2*y[8]+2*y[18];
(2)*v.32+(2)*v.33+(2)*v.34+(2)*v.42
> e*f;
(4)*v.49+(6)*v.50+(8)*v.51+(12)*v.52
> h:=e*f;;
\end{verbatim}
Pour cette orbite, les donn{\'e}es de \cite{Triplets} conviennent, pour
des raisons qui m'{\'e}chappent.\\
\\
Calcul de {\tt g}, {\tt Bg}, {\tt z} et {\tt Bz}\,:
\begin{verbatim}
> g:=LieCentralizer(L,Subspace(L,[e]));
<Lie algebra of dimension 12 over Rationals>
> Bg:=BasisVectors(Basis(g));;
> z:=LieCentre(g);
<two-sided ideal in <Lie algebra of dimension 12 over Rationals>, 
(dimension 3 )>
> Bz:=BasisVectors(Basis(z));
[ v.8+v.9+v.10+v.18, v.23, v.24 ]
> h*Bz[2];
(6)*v.23
> h*Bz[3];
(6)*v.24
\end{verbatim}  
Les poids de {\tt z} sont 2,6,6; d'o{\`u} $m_{r}=6$. Il n'y a qu'une matrice {\`a} {\'e}tudier.

\begin{enumerate}
\item $m_{i_{1}}=16$,  $m_{k_{1}}=2$
\begin{verbatim}
> h*Bg[1]; h*Bg[2]; h*Bg[3]; h*Bg[4]; h*Bg[5];
(2)*v.8+(2)*v.9
(2)*v.4+(-1)*v.11+v.12
(2)*v.12+(-2)*v.13
(2)*v.6+(-1)*v.7+v.14
(2)*v.14+(-2)*v.15
> h*Bg[7];                                                                 
(2)*v.10+(2)*v.18

> ((f*Bg[1])*Bz[2]); ((f*Bg[2])*Bz[2]); ((f*Bg[3])*Bz[2]);
((f*Bg[4])*Bz[2]); ((f*Bg[5])*Bz[2]); ((f*Bg[7])*Bz[2])
(-4)*v.23
(-2)*v.24
(-2)*v.24
0*v.1
0*v.1
(-2)*v.23

> ((f*Bg[1])*Bz[3]); ((f*Bg[2])*Bz[3]); ((f*Bg[3])*Bz[3]);
((f*Bg[4])*Bz[3]); ((f*Bg[5])*Bz[3]); ((f*Bg[7])*Bz[3])
(-4)*v.24
0*v.1
0*v.1
(-2)*v.23
(-2)*v.23
(-2)*v.24
\end{verbatim}
La matrice {\`a} {\'e}tudier est 
$$\left[
\begin{array}{cccccc}
-4 \alpha & 0          & 0           & -2 \beta  & -2 \beta & -2 \alpha \\
-4 \beta  & -2 \alpha  & -2 \alpha  & 0 & 0 & -2 \beta
\end{array}
\right] \cdot $$ 
C'est une matrice de rang 2 pour tout couple $(\alpha,\beta)$ non
nul. 
\end{enumerate}
La propri{\'e}t{\'e} $(P)$ est v{\'e}rifi{\'e}e pour cette orbite.\\
\end{enumerate}
{\bf Conclusion\,:} $F_{4}$ v{\'e}rifie $(P)$.

\nocite{*}
\bibliographystyle{plain}
\bibliography{biblio_indice_normalisateur}

\vspace{3cm}

\normalsize
~\\
{\sc Anne Moreau \\
~\\
Universit{\'e} Paris 7 - Denis Diderot,\\
Institut de Math{\'e}matiques de Jussieu,\\
Th{\'e}orie des groupes,\\
Case 7012\\
2, Place jussieu\\
75251 Paris Cedex 05, France. }\\
~\\
{\em E-mail\,: {\tt moreaua@math.jussieu.fr} }

\end{document}